\documentclass[11pt]{article}
\usepackage[margin=0.75in]{geometry}
\usepackage{setspace} 

\usepackage{graphicx} 
\usepackage{multirow}
\usepackage{ragged2e}
\usepackage{enumitem}

\usepackage{type1cm}        
\usepackage{comment}
\usepackage[lofdepth,lotdepth]{subfig}

\usepackage[table]{xcolor}
\usepackage{makeidx}
\usepackage{graphicx}
\usepackage{multicol}
\usepackage[bottom]{footmisc}
\usepackage{xcolor}
\usepackage{soul}
\usepackage{amsmath}
\usepackage{microtype}

\usepackage{mdwlist}
\usepackage[printonlyused,withpage]{acronym}
\usepackage{makecell}
\usepackage{newtxtext}
\usepackage[varvw]{newtxmath}
\usepackage{float}
\newcommand{\R}{\mathbb R}

\usepackage{hyperref}

\usepackage{amsopn}

\def \A{\mathcal{A}}

\usepackage{lipsum}
\usepackage{amsfonts}
\usepackage{graphicx}
\usepackage{epstopdf}
\usepackage{algorithmic}
\ifpdf
  \DeclareGraphicsExtensions{.eps,.pdf,.png,.jpg}
\else
  \DeclareGraphicsExtensions{.eps}
\fi

\newtheorem{definition}{\noindent $\mathbf{Definition}$}[section]

\title{
Learning the Geometry of Data: \\A Mathematical Review of Shape Space Analysis
}

\author{
Gary P. T. Choi$^{1,\dagger}$, Khanh Dao Duc$^{2,\dagger}$, Shira Faigenbaum-Golovin$^{3,\dagger,\ast}$, Karen Habermann$^{4,\dagger}$, \\Emmanuel Hartman$^{5,\dagger}$, Christoph von Tycowicz$^{6,\dagger}$, Chi Zhang$^{7,\dagger}$, Wenjun Zhao$^{8,\dagger}$, Felix Zhou$^{9,\dagger}$\\
\\
\footnotesize{$^1$Department of Mathematics, The Chinese University of Hong Kong, Hong Kong SAR} \\
\footnotesize{$^2$Department of Mathematics, University of British Columbia, Vancouver, BC V6T 1Z4, Canada} \\
\footnotesize{$^3$Department of Mathematics, Bar-Ilan University, Ramat-Gan, Israel}\\
\footnotesize{$^4$Department of Statistics, University of Warwick, Coventry, CV4 7AL, United Kingdom}\\
\footnotesize{$^5$Department of Mathematics, University of Houston, Houston, USA}\\
\footnotesize{$^6$Zuse Institute Berlin, Takustraße 7, 14195 Berlin, Germany} \\
\footnotesize{$^7$Department of Biomedical Sciences, Texas A\&M University, Dallas, TX 75246, USA}\\
\footnotesize{$^8$Department of Mathematics, Wake Forest University, Winston-Salem, NC 27106, USA}\\
\footnotesize{$^9$Lyda Hill Department of Bioinformatics, University of Texas Southwestern Medical Center, Dallas, TX 75390, USA}\\
\footnotesize{$^\dagger$Author order was determined by alphabetical order; all authors contributed equally to the synthesis of this review.}\\
\footnotesize{$^\ast$Corresponding author, E-mail: shira.golovin@biu.ac.il}
}

\date{ }

\begin{document}

\maketitle

\begin{abstract}
A fundamental objective of machine learning is to analyze structural discrepancies and identify patterns within data. Advancements in acquisition techniques have increasingly produced datasets where data takes shape and has its own complex geometry. These ``shape spaces'' often exhibit subtle variations between objects, presenting unique challenges that traditional machine learning frameworks frequently overlook. This article reviews the emerging field of shape space analysis, providing a unified mathematical and computational foundation for studying collections of geometric data. By standing on the shoulders of classical analysis, this framework introduces specialized tools designed to leverage the structural richness and nuances of shape-based datasets. The analytical pipeline begins with optimal shape parameterization to capture geometric variations, followed by a discussion of the choice of a robust distance metric. Through the lens of manifold learning approaches, we illuminate dynamic trajectories and establish rigorous statistical inferences for fluctuations. By extending classical machine learning into geometry-aware frameworks, we bridge theoretical abstractions with empirical rigor.  This foundation is grounded by two multiscale case studies (subcellular morphology, primate dental evolution) offering a roadmap for researchers as they transition from classical learning to geometry-aware domains. Finally, we evaluate cross-disciplinary opportunities and identify pivotal open questions currently facing the emerging shape space community. This review reveals that, despite the vast diversity of applications spanning neurology, microbiology, anthropology, and computer vision, researchers face similar challenges rooted in unaligned, subtly varying geometries. We identified a significant fragmentation in terminology across disciplines, which this survey aims to resolve by formulating a unified language for shape space analysis. Our findings demonstrate that geometric learning frameworks consistently outperform standard models, particularly when navigating limited labeled datasets and intricate local variations. Ultimately, we show that what appear to be domain-specific hurdles are actually universal geometric properties, requiring a specialized mathematical approach to unlock structural insights across all scales of data. This comprehensive survey provides a unique view by synthesizing perspectives from mathematics, statistics, and computer science to establish a foundational reference for high-dimensional, geometrically constrained data. While traditional studies focus on single-object analysis, this work demonstrates the superior rigour achieved by analyzing entire shape collections to capture structural variation. We offer a practical guide for researchers navigating application-driven challenges. By bridging abstract geometric theory with empirical learning across scales, from microbiology to anthropology, the article offers a unique, unified view of shape modeling. It transforms domain-specific observations into a cross-disciplinary framework for any data possessing physical or structural form.
\end{abstract}


\section{Introduction}

\label{sec:intro}

As we look around us today, we can easily identify the ever-growing diversity of low- and high-dimensional data of shapes emerging across various fields. What has fueled this rapid proliferation is the technological boost (efficient and affordable acquisition methods, the exponential increase in computational power, and the plummeting costs of memory), together with the revolution in machine learning and data analysis through tools from the fields of mathematical, statistical, and computer science analysis, resulting in an explosion of opportunities.

\begin{figure}[t]
\centering
\includegraphics[width=\linewidth]{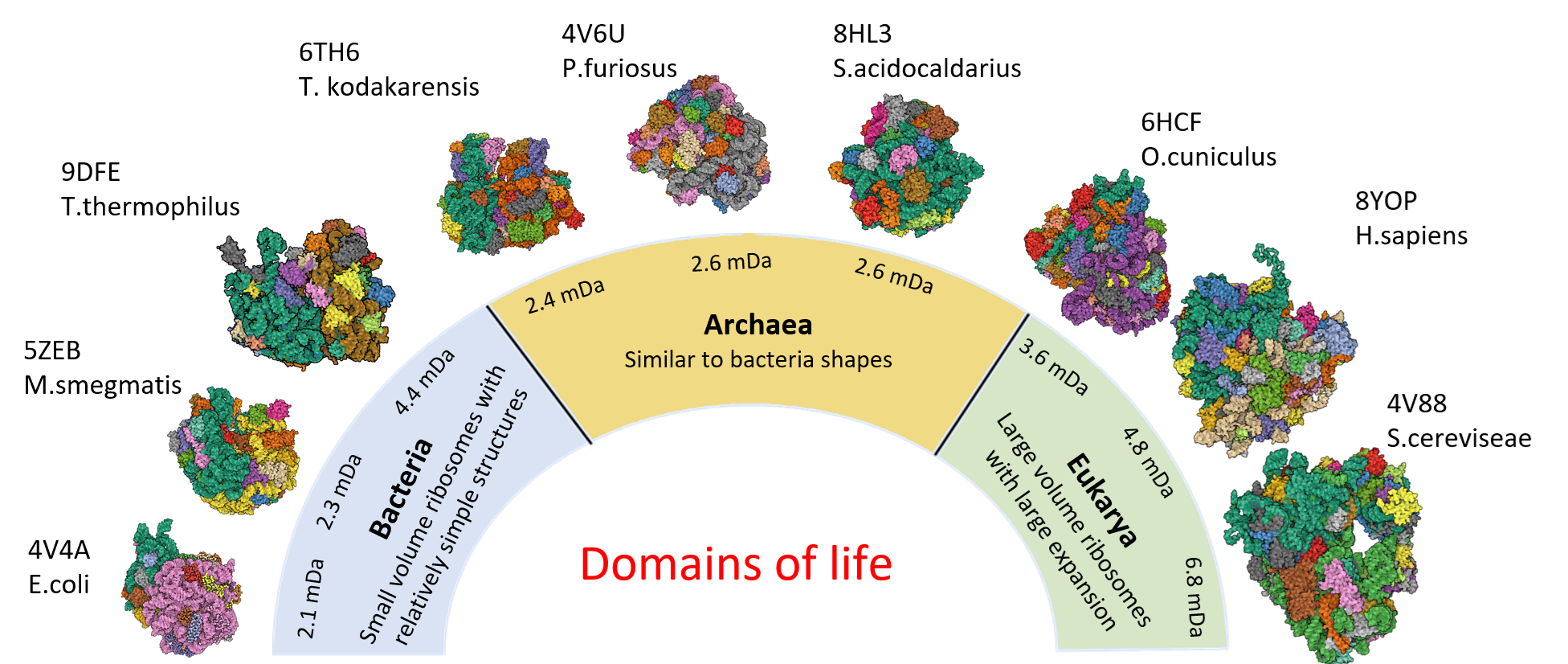} 
\caption{\textbf{Form is function -- the ribosome shape space.} Inspired by the concept of Kendall’s shape space~\cite{kendall2009shape}, our goal is to study the geometric properties of a collection of shapes, analyzing their differences, similarities, and patterns to uncover insights that would otherwise remain hidden. Ribosome volume correlates with functional complexity across domains of life. Each structure occupies a unique position in shape space, with increasing volume and complexity corresponding to expanding functional capabilities (e.g., translational fidelity~\cite{fujii2018decoding}, co-translational folding~\cite{willmund2013cotranslational}, chaperone binding~\cite{kristensen2003chaperone}).
}
\label{fig:shape_space_illustration}
\end{figure}

Developing methods for the challenging tasks of analyzing individual shapes (such as images or 3D point clouds) has matured over the past few decades~\cite{burger2022digital, szeliski2022computer, salomon2007curves}. However, the fast-emerging era, motivated by the new data, poses new mathematical challenges and opportunities. While one can study the intrinsic properties and analyze a single shape (e.g., shape curvature~\cite{shan2019ariadne}), leveraging the power that exists in a collection of shapes (each is a smooth/non-smooth transformation) can enhance the results significantly. For instance, denoising a single image can be a challenge, while denoising a collection of images may be more efficient~\cite{faigenbaum2023manifold, zhang2024robust}. 

The concept of shape space encompasses the analysis of a collection of shapes, aiming to uncover intrinsic patterns and structures within the dataset. As before, the shapes can come in different representations, starting from images, 3D surfaces, or other data structures in low or high dimensions. Shape space analysis is the study of the geometric and statistical properties of shapes, treating them as points within a high-dimensional space designed to capture their key features. The goal is to model and measure the differences, similarities, and patterns among shapes, enabling insights that are otherwise inaccessible. Typically, shapes exhibit only subtle variations from each other, and a primary objective of shape analysis is to identify, quantify, and compare these differences. For example, key questions in this field include shape classification, reconstruction, and alignment, as well as the study of shape variability and evolutionary trends between different shapes. Figure~\ref{fig:shape_space_illustration} illustrates the concept of shape space in the context of ribosomal structures. Examining the geometric variability within ribosomal data is of particular importance, as these structural differences are closely associated with the functional complexity of the ribosome. Questions naturally arise, such as: What is the intrinsic dimensionality of the manifold from which these ribosomes are sampled? Which metric best captures biologically meaningful similarities and differences between ribosomal structures? How can we develop computational tools that preserve both structural and functional relationships within the shape space?

The shape space applications span diverse fields, ranging from medical imaging (starting from the scanned brain \cite{nitzken2014shape, lande_quantitative_1979, guo2022statistical} down to the mycobacteria and cancer cells, e.g.,~\cite{fedorov2023national, cancercell2}), to studying evolutionary anthropology (fruit fly wings \cite{geldenhuys2023deep}, teeth and bones, e.g.,~\cite{faigenbaum2026studying}), archaeology and paleography (flint, bones, and ancient handwritten letters, e.g.,~\cite{faigenbaum2016algorithmic}), as well as computer graphics, and CAGD (faces morphometrics~\cite{egger20203d}, human body pose, e.g.,~\cite{hartman2024varishape, cheng2018parametric}).
These datasets present a range of challenges for shape comparison. The challenges begin with foundational issues (such as data acquisition and denoising) and extend to identifying a suitable parameterization that preserves key geometric and field related features. Naturally, the collected data is rarely aligned, and robust registration across the dataset becomes a critical step. Once alignment is achieved, a central question arises: What is the most appropriate metric for comparing shapes? Such a metric must be invariant under basic transformations (e.g., scaling, rotation, and translation) as well as computationally efficient, and ideally should minimize the need for manual intervention. Importantly, it must also be sensitive enough to detect subtle but meaningful differences. Notably, the challenges in shape space analysis often differ from, and can be more demanding than, those encountered in traditional computer graphics, CAGD, or image processing. This is largely because the differences between shapes in biological datasets are frequently subtle, requiring highly sensitive and robust analytical tools to detect and interpret them meaningfully. While these tasks are nontrivial, having access to a collection of related shapes can help alleviate some of the challenges. As demonstrated in this paper, shape space analysis becomes significantly more powerful when applied to shape collections rather than to isolated specimens.

Shapes are ubiquitous, thus making it impractical to comprehensively survey all existing literature on shape analysis in a single paper. Instead, this survey provides an overview of key concepts in the study of shape space and offers a practical guide on how to approach shape space analysis, highlighting available methods and datasets. Our study contributes to the literature in several ways. To the best of our knowledge, this is the first comprehensive survey addressing shape space problems across diverse applications. These applications generate inherently diverse shapes that require distinct preprocessing and modeling techniques. Despite these differences, we identify common steps and overarching goals that unify the approaches. We also provide a practical “recipe” for successfully studying shape spaces, offering insights into accessible datasets and essential methodologies. Unlike prior studies limited to specific fields, our work integrates methods from mathematics, statistics, and computer science, leveraging the interdisciplinary expertise of the authors to offer a wide-ranging and balanced overview. This paper serves both as an accessible entry point for newcomers and a valuable resource for experienced researchers, highlighting open challenges and emerging trends. Moreover, it lays a foundation for benchmarking new methods by facilitating comparisons to state-of-the-art approaches, ensuring it will remain a vital resource for shaping future advancements in the field.

When discussing shape space analysis, the foundational principles have been extensively covered in several key books, e.g.,~\cite{younes2010shapes, dryden2016statistical, pennec2019book}, which survey the fundamental ideas of shape analysis. In~\cite{younes2010shapes}, Younes describes a large range of methods that have been proposed to represent, detect, or compare shapes (or more generally, deformable objects), together with the necessary mathematical background that they require. In~\cite{dryden2016statistical}, Dryden and Mardia present a comprehensive statistical perspective on shape analysis, emphasizing the joint study of size and shape (form). This edition expands the discussion to include unlabeled size-and-shape analysis, a broader range of three-dimensional applications, and an extended treatment of general Riemannian manifolds. The book~\cite{pennec2019book} by Pennec, Sommer, and Fletcher surveys Riemannian geometric statistics and its use in medical image analysis. These books provide a solid mathematical background on the existing tools. In this paper, we aim to survey the rapid transformation made in the field of shape space over the past six years.

Additional comprehensive reviews on specific topics can be found in the following works. An overview of non-rigid 3D shape analysis is provided by~\cite{laga2018survey}, which covers key aspects such as shape representations, metrics, registration, and statistical modeling. Bauer et~al.~\cite{bauer2014overview} offer an in-depth treatment of shape spaces obtained by quotienting out (re-)parameterizations, while Mitteroecker and colleagues~\cite{mitteroecker2022thirty} review geometric morphometrics with an emphasis on landmark-based approaches. For applications in medical imaging, see reviews on statistical shape models for 3D segmentation~\cite{heimann2009statistical, ambellan2019statistical}.

\textbf{Learning? Yes or No.} With the rise of machine learning, and particularly neural network-based approximation techniques that rely on training with ground truth examples, many challenges are now being addressed through model training. However, shape space analysis presents unique challenges. The main one is data availability, which remains extremely limited due to the rarity of artifacts. Moreover, labeled data for supervised learning is scarce, if available at all, posing further obstacles to the application of such methods. While transfer learning offers a potential solution, identifying a sufficiently similar dataset to the specific problem at hand is often challenging. Consequently, much of the field of shape space analysis continues to rely on classical approaches rooted in mathematics, statistical methods, and computer science, focusing on the intrinsic geometrical properties of shapes themselves. Often, an analytical solution is far more valuable than brute force, for insight can achieve what sheer effort cannot.

\textbf{Critical Learning.} 
Since real-world datasets often contain limited labeled material, evaluating the performance of proposed methods becomes inherently challenging. To address this, it is crucial to design controlled test-case scenarios in which new pipelines or methodologies can be systematically assessed and validated. In addition, while machine learning (ML) and artificial intelligence (AI)  can be immensely powerful, they should be applied with care and a clear understanding of their assumptions and limitations. As noted by Calder et~al.~\cite{calder2022use}, the growing accessibility of ML algorithms-paired with limited expertise in their proper use outside of STEM disciplines-has led to widespread misapplications across the literature.

\textbf{How to read this paper?} To make this paper
accessible to readers with various interests and
backgrounds, it is organized in the following way. In the Introduction (Section~\ref{sec:intro}), we set the stage by exploring what makes shape spaces fundamentally different from other data domains. Section~\ref{sec:data} turns to the data itself: how shapes are acquired, how they take form (implicitly, explicitly, or otherwise), and how one may transition between different representations. In Section~\ref{sec:shape_space_overview}, we review computational approaches for representing shapes, weighing their respective strengths and weaknesses. Section~\ref{sec:methods} then surveys the methodological backbone of the field: preprocessing, alignment, landmarking, and registration, before delving into shape parameterization (Section~\ref{ref:shape_param}) and the choice of shape feature descriptors (Section~\ref{sec:feature_descriptors}). We continue in Section~\ref{sec:distance} with a discussion of distance metrics and similarity measures, highlighting how each incorporates our understanding of geometric variability. The following section, Section~\ref{sec:dynamics}, traces shape dynamics, learning manifolds, constructing temporally invariant coordinate systems, and mapping trajectories through shape space. In Section~\ref{sec:mean}, we define meaningful notions of mean and fluctuation, essential for statistical inference on manifolds. Building on these, Section~\ref{sec:Statistical Analysis} explores machine learning and statistical methods for uncovering hidden geometric patterns, emphasizing the extensions to classical techniques, and geometry-aware learning becomes indispensable. In Section~\ref{sec:practicalGuide}, we bring theory to practice by offering a hands-on guide for researchers working with two real-world shape spaces. Section~\ref{sec:datasets} and Section~\ref{sec:tools} curate accessible datasets and publicly available tools, encouraging reproducibility and further experimentation. Finally, in Section~\ref{sec:Discussion} we close by posing several open questions, some conceptual, others computational, that we hope will spark dialogue, collaboration, and innovation within this emerging community. Together, these discussions aim to bridge data-driven exploration with geometrically meaningful discovery, inviting readers to rethink how data takes form and how form reveals knowledge.

A few reflections on the process of writing this paper. Working on it compelled us to pause, reconsider, and distill the terminology and categorization of the methods we describe. Through many lively and fruitful discussions, we gradually clarified for ourselves-and ultimately in writing-the key ideas that underlie shape space analysis. The diversity of our backgrounds proved to be a true advantage, allowing us to illuminate the concept of shape space from multiple vantage points. At the same time, we became aware of how differently terminology is used across disciplines. Final disclaimer, while we have aimed to provide a broad and representative overview of shape space analysis, it is naturally impossible to capture the full breadth of work in this diverse and rapidly evolving field. Our goal was to highlight key methodologies and conceptual frameworks rather than offer an exhaustive review. We acknowledge that some valuable contributions may have been inadvertently omitted, and we sincerely apologize in advance.

\section{Shape in Data}
\label{sec:data}

\begin{figure}[t]
\centering
\includegraphics[width=1\linewidth]{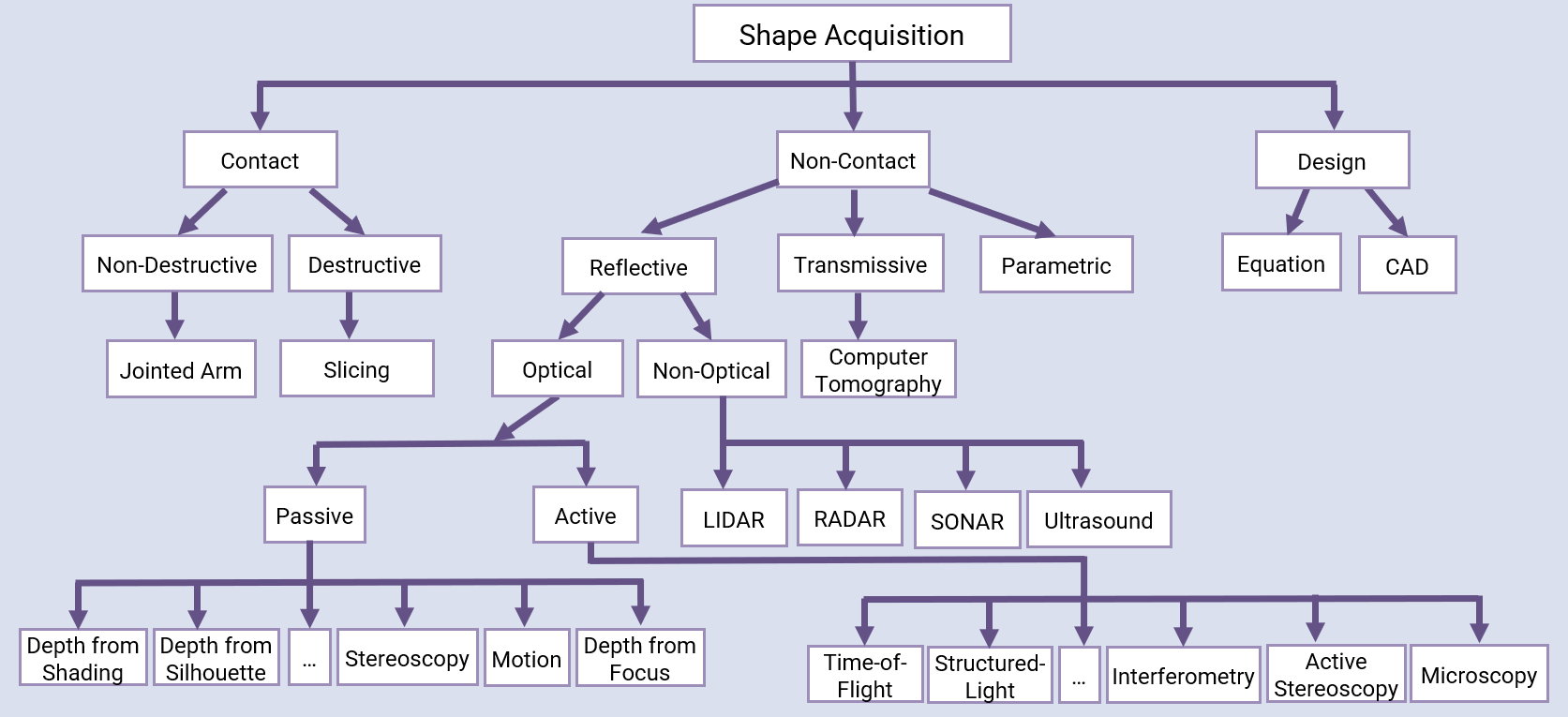} 
\caption{\textbf{Taxonomy of shape acquisition technologies adapted from~\cite{giancola2018survey}.}}
\label{fig:shape_acquisition_technologies}
\end{figure}

The notion of a shape, a circle, a square, or a rectangle is inherently abstract. Various definitions have been developed and used to enable the study of shapes mathematically, see Section~\ref{sec:shape_space_overview}. With respect to computation, shapes are acquired by various technologies and stored in different digital formats. Typically, the acquisition method also determines the storage format. Each format presents advantages and disadvantages for the downstream analyses discussed in Section~\ref{sec:methods}. This section introduces the major acquisition methods, common data formats, and methodologies for interconverting between formats.

\subsection{Shape Data Acquisition}
\label{subsec:shape_data_acquisition}
Shape acquisition is a specialist field of study wider than the scope of this survey paper. 
For more detailed operation principles of individual methods, we refer the reader to dedicated literature~\cite{giancola2018survey}. An overview of different methods is given in Fig.~\ref{fig:shape_acquisition_technologies}. Firstly, acquisition can be split based on physical \textit{contact}, \textit{non-contact}, or \textit{de novo design} techniques. \textit{Contact} techniques can be \textit{destructive}, such as slicing and sectioning an object into 2D shapes successively to be reassembled together. Examples include histology~\cite{kiemen2022coda} and focused ion beam scanning electron microscopy~\cite{knott2008serialfibsem,xu2021fibsem} of 3D biological tissue specimens. Volumetric imaging scans a 3D shape to produce a cube of 2D cross-sectional slices~\cite{cevidanes2006image, yapuncich2019digital, o2024masse}, to be later reconstructed into a 3D object computationally. \textit{Contact} techniques can also be \textit{non-destructive}, such as jointed arm robots that slowly but accurately probe 3D points. \textit{Non-contact} techniques sample regions instead of single points on a target based on measuring the \textit{reflective} or \textit{transmissive} properties of emitted waves or other \textit{parametric} properties. For instance, magnetic resonance imaging is a \textit{parametric} technique that measures the change in magnetization, specifically, the rate of relaxation of nuclear spins following their perturbation by an oscillating magnetic field. Computer tomography is a \textit{transmissive} technique that uses X-ray signals taken from different emitter positions to identify changes in density
within a body. Alternatively, \textit{reflective} techniques focus on analyzing a signal's reflection. \textit{Non-optical} techniques focus on wavelengths outside of the visible and infrared spectrum. Sound Detection and Ranging (SONAR) and Radio Detection and Ranging (RADAR) use radio signals to estimate range maps based on the time the signals run through their environment. Ultrasound measures the intensity of the recollected sound waves. \textit{Optical} techniques exploit the visible and infrared (IR) light to sample a scene or an object. While the visual spectrum aligns with the human vision system, IR wavelengths carry temperature information and are usually more robust to ambient light. \textit{Optical} techniques can be \textit{passive} or \textit{active}. Passive methods use the reflection of natural light on a target to measure its shape. Stereophotogrammetry looks for homogeneous features from multiple cameras to reconstruct a 3D shape, using triangulation and epipolar geometry. Structure from motion (SfM) additionally recovers the camera motion around the object given multiple overlapping viewpoints.
Shape from Silhouette and Shape from Shading measure shape information based on edges and shading theory. Depth of field uses the focus information of the pixels given a sensor focal length to estimate range. \textit{Active} methods enhance shape acquisition by using an external lighting source that provides additional information. \textit{Time-of-flight} systems use the Light Detection And Ranging (LIDAR) principle to estimate depth by emitting light and measuring the time the light goes back-and-forth. \textit{Structured-light} devices project a laser pattern to the target and estimate the depth by triangulation. To cope with depth maps, structured-light cameras project a 2D codified pattern to perform triangulation with. The \textit{active photogrammetry} principle is similar to the passive one but looks for artificially projected features. In contrast with structured-light the projected pattern is not codified and only serves as an additional feature to triangulate with. \textit{Interferometry} projects a series of fringes such as \textit{Moire's} to estimate shapes, requiring iterative spatial refinement in the projected pattern. Microscopy uses cameras to measure the intensity from the reflected visible spectrum (brightfield imaging), from phase shifts of the light source through a sample (phase contrast), or the fluorescence from labeled specimens. Finally, in \textit{de novo} \textit{design}, the shape is explicitly constructed, generally through equations as in the plotting of a circle, or algorithmically, through computer-aided design (CAD) software like Blender, Autodesk 3ds Max, and Unity for applications like 3D animation.

\subsection{Shape Data Formats}
\begin{figure}[t]
\includegraphics[width=1.0\linewidth]{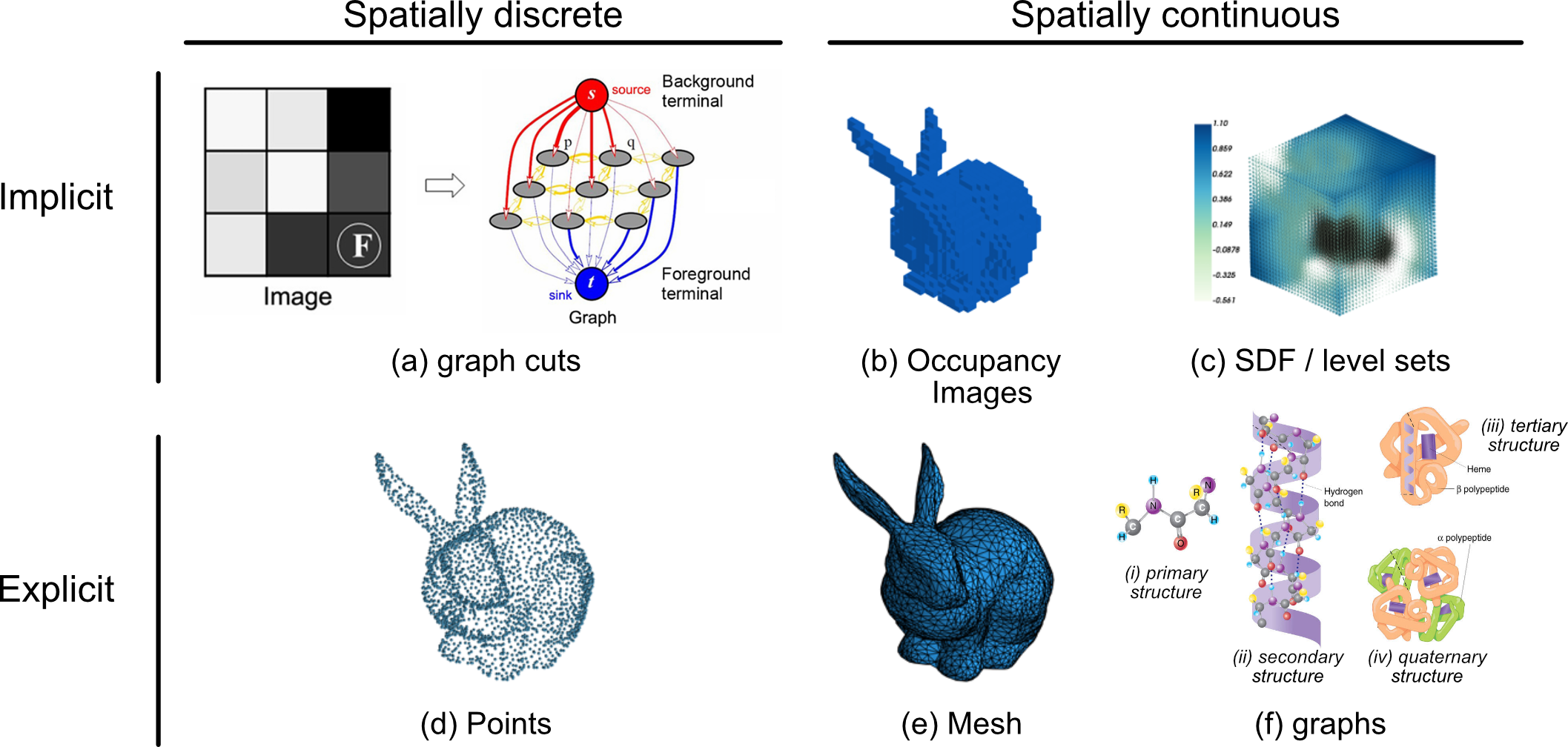} 
\caption{\textbf{Commonly occurring digital formats representing 3D shape implicitly (by interior) or explicitly (by boundary), in a spatially discrete or continuous setting.} (a) Graph cut formulation of an image, connecting pixels to graph nodes with edges representing the strength of affinity or dissimilarity. Illustration adapted from~\cite{boykov2006graph}. (b)--(e) different representations of the Stanford bunny, adapted from~\cite{michalkiewicz2019deep}. (f) Protein structure including primary: amino acid sequence; secondary: how primary structure folds under hydrogen bonds; tertiary: coiling or pleating of the secondary structure through hydrogen, disulphide, ionic bonds and hydrophobic, hydrophilic interactions; quaternary: joining of multiple tertiary structures.}
\label{fig:common_representations_of_3D_shape_in_data}
\end{figure}

Except for \textit{de novo design}, generally, shape acquisition methods do not provide the shape directly as raw data. Instead, computation is necessary to extract individual shapes. This is done by segmentation, the task of labeling all points or voxels belonging to the same object with the same numerical index ID. The segmentation method used is often specific to the object and data type. Point-cloud segmentation techniques~\cite{qi2017pointnet, qi2017pointnet++} are used for point-based acquisitions like those from SfM or time-of-flight systems. Meanwhile, image-based segmentation techniques~\cite{ronneberger2015u,stringer2021cellpose,zhou2025universal} can be used for regularly sampled, grid-like data. From the segmentation, the shape geometry can be stored in digital formats that represent the shape differently. Broadly, one can distinguish between \textit{implicit} and \textit{explicit} representations that capture the shape boundary in a spatially \textit{discrete} or \textit{continuous} manner~\cite{cremers2015image}. Fig.~\ref{fig:common_representations_of_3D_shape_in_data} highlights common shape formats under this categorization.  More details about each category are provided below, alongside with discussion of when it is beneficial to use specific representations.

\subsubsection{Implicit Representations}
Implicit representations specify shape by labeling all points in a physical space as being part of the interior or exterior of the object. 

In the spatially discrete setting, the shape is formulated as a discrete optimization problem, typically in the framework of \textbf{graph cuts}~\cite{boykov2001interactive}. The shape's interior and exterior are specified by finding an optimal partitioning of the adjacency graph representation of an image (Fig.~\ref{fig:common_representations_of_3D_shape_in_data}(a)). In this graph, nodes are individual pixels assigned with a foreground likelihood score, e.g., image intensity, and the edges encode similarities between pixel pairs based on their statistics, e.g., in terms of radial basis function (RBF) kernels.

In the spatially continuous setting, a common representation is \textbf{occupancy images}. The pixel values of an occupancy image reflect the likelihood that the pixel is part of the object interior. The object exterior is typically regarded as background and zero-valued (Fig.~\ref{fig:common_representations_of_3D_shape_in_data}(b)).
Examples of occupancy images are the foreground probability maps predicted by neural networks for segmentation, which score the likelihood of individual pixels being foreground from 0 to 1. A binary segmentation takes a value of 1 for the shape's interior and a value of 0 for the background. An instance image segmentation can also be regarded as an occupancy image where each unique integer ID represents a shape. A classic example is the MNIST handwriting dataset~\cite{lecun1998mnist}, where each binary image represents the shape of a handwritten number. Alternatively, based on a binary segmentation, shape can be specified by a \textbf{signed distance function (SDF)}, also termed \textbf{level sets}, with respect to the boundary, typically defined as distance = 0 (Fig.~\ref{fig:common_representations_of_3D_shape_in_data}(c)). Interior pixels have increasingly positive-valued distances the further they are from the boundary, with a maximum at the shape center. Meanwhile, exterior pixels have increasingly negative-valued distances the further away from the boundary. The opposite sign convention is also used: exterior being positive-valued and interior being negative-valued. Distances are computed by solving partial differential equations~\cite{cremers2015image}.

\subsubsection{Explicit Representations}
Explicit representations parameterize the shape geometry, its boundary or area in 2D, and surface or volume in 3D, specifying a coordinate system that can be used to explicitly index locations on the shape's surface or volume.

In the spatially discrete setting, the predominant approach is to represent shapes by a collection of \textbf{points}, commonly referred to as a point cloud, whose locations are sampled from the surface or volume. The point cloud can cover the shape sparsely or densely, uniformly or non-uniformly (Fig.~\ref{fig:common_representations_of_3D_shape_in_data}(d)). 

In the spatially continuous setting, the connectivity between individual sampled points is also known. \textbf{Meshes} connect surface points in polygons and volume points in polyhedra (Fig.~\ref{fig:common_representations_of_3D_shape_in_data}(e)). Surface polygons are commonly referred to as faces, and volume polyhedra as facets. Strictly speaking, a proper shape mesh should have non-overlapping and non-intersecting faces and facets. Moreover, they should have properties of a surface or volume, namely, be single-component, manifold, and watertight. Meshes that do not have these properties are referred to as ``polygon soups''. The distinction is important. Polygon soups do not necessarily allow for mesh-based processing algorithms available in software libraries, such as Laplacian smoothing. In such cases, mesh repair algorithms can be applied to create a proper mesh, including fixing non-manifold edges and hole patching \cite{attene2010lightweight}, remeshing \cite{valette2008generic, faraj2016multi, botsch2004remeshing}, surface reconstruction \cite{kazhdan2006poisson, huang2018robust}, or shrinkwrapping \cite{portaneri2022alpha,zhou2023surface}. If this is not possible, the mesh can be treated as a point cloud \cite{vasan2025interpretable}. \textbf{Graphs} generalize meshes (Fig.~\ref{fig:common_representations_of_3D_shape_in_data}(f)). Graph nodes represent geometric coordinates, but the connectivity between nodes may be arbitrarily connected, including crossovers and long-range connections, and does not need to correspond to a tessellation of space like a surface or volume. This flexibility allows graphs to represent complex shapes, particularly those with hierarchical structure, notably a protein (Fig.~\ref{fig:common_representations_of_3D_shape_in_data}f). The primary structure is the linear sequence of amino acids comprising the protein. The secondary structure is how the amino acids locally fold. The tertiary structure is the overall 3D shape of the single polypeptide chain, and the quaternary structure is the resulting assembly of multiple polypeptide chains into a single functional protein complex. The primary structure can be captured by simple graphs, but the higher-order shapes require specifying hypergraphs, involving simplicial and cell complexes~\cite{hajij2022topological}. These structures define additional layers of connectivity among subgroupings of points.

\subsubsection{Different Representations for Different Problems} 
The choice of representation dictates the class of shapes that can be modeled and the computational processes that can be performed. 

Generally, it is observed that:
\begin{itemize}
    \item \textbf{Implicit} representations easily generalize to shapes in arbitrary dimensions. Respective algorithms (graph cuts, level set methods) also straightforwardly extend from two to three or more dimensions. For explicit shape representations, the extension to higher dimensions is nontrivial. The arc-length parameterization of curves does not extend to parameterize surfaces.
    \item \textbf{Implicit} representations readily generalize to capture arbitrary topology. The labeling of space as being inside or outside places no topological constraints. In contrast, explicit spatially continuous parameterizations like curves and meshes require sophisticated splitting and merging techniques to model the transition from a sphere to a donut, or the splitting of a single object into multiple objects.  
    \item \textbf{Explicit} representations enables the capturing of point correspondence (see Section~\ref{sec:Landmarking}), and the alignment of semantic shape parts.
    \item \textbf{Explicit} representations often allow for more straightforward and intuitive modeling of shape similarity. We can interpolate between and deform two curves or meshes to obtain intermediate shape meshes, by minimizing appropriate energies \cite{sassen2024repulsive}. The linear interpolation of implicit representations is ill-posed. Convex combinations of
 binary-valued functions are no longer binary-valued. Convex combinations of SDFs are generally no longer SDFs. 
\end{itemize}

Tables~\ref{table:pros_cons_implicit_representations} and~\ref{table:pros_cons_explicit_representations} summarize the advantages and disadvantages of implicit and explicit shape representations, respectively with consideration of the shape analysis methods presented in Section~\ref{sec:methods}.

\begin{table}[t!]
\begin{center}
\begin{tabular}{|p{0.05\textwidth}|p{0.23\textwidth}|p{0.32\textwidth}|p{0.31\textwidth}|}
 \hline
 \rowcolor{gray!30}
  & \centering\textbf{Graph Cuts} & \centering\textbf{Occupancy Images} & \textbf{Signed Distance Function (Level Set)} \\
\hline
\hline
 Pros & \begin{itemize}[leftmargin=*]
     \item Captures both discrete objects and connectivity structure
     \item Robust to topological defects
     \item Natively supports connectivity-aware operations like labelspreading, diffusion, Markov clustering
     \item Easy to extend to higher spatial dimensions
 \end{itemize} & \begin{itemize} [leftmargin=*]
    \item Regular data structure
     \item Geometry: captures surface and volume
     \item Easy for learning
     \item Unique representation
     \item Dense representation: easy to compare similarity and register volume
 \end{itemize} & \begin{itemize}[leftmargin=*]
    \item Regular data structure
     \item Geometry: captures surface and volume
     \item Easy for learning
     \item Unique representation
     \item Customize metric 
     \item Dense representation: compare similarity and register volume
     \item Surface is zero level set
     \item Native support for many shapes in parallel
 \end{itemize}\\ 
\hline
 Cons & \begin{itemize}[leftmargin=*]

    \item No explicit geometry therefore point correspondences and length measurements are not supported
    \item No explicit indexing of locations in shape
    \item Difficult to specify to capture multiple different shapes
 \end{itemize} & \begin{itemize}[leftmargin=*]
    \item Image size determines resolution
     \item Hard to compare similarity and register surfaces
     \item Poorly-suited for computing simultaneously with multiple shapes 
     \item No constraints for computing only on surface
     \item High memory consumption
     \item No explicit traversal of shape
     \end{itemize} & \begin{itemize}[leftmargin=*]
    \item Image size determines resolution
     \item Hard to compare similarity and register surfaces
     \item Expensive distance computations
     \item High memory consumption
     \item No explicit traversal of shape
     \end{itemize}\\ 
 \hline
\end{tabular}
\end{center}
\caption{Advantages and disadvantages of implicit shape representations.}
\label{table:pros_cons_implicit_representations}
\end{table}

\begin{table}[t!]
\begin{center}
\begin{tabular}{|c|p{0.26\textwidth}|p{0.28\textwidth}|p{0.28\textwidth}|}
 \hline
 \rowcolor{gray!30}
 & \centering\textbf{ Points} &\centering \textbf{ Mesh} &\textbf{Graph} \\
\hline
\hline
 Pros & \begin{itemize} [leftmargin=*]
     \item Compact \& memory-efficient
     \item Robust to topological defects
     \item Easy to resample
     \item Supports learning
     \item Supports multiple shapes 
 \end{itemize} &
 \begin{itemize}[leftmargin=*]
     \item Compact and memory-efficient
     \item Accurate, detailed geometry
     \item Explicit neighbor connectivity enables shape deformation
 \end{itemize} & \begin{itemize}[leftmargin=*]
     \item Compact and memory-efficient
     \item Arbitrary specification of geometry
     \item Supports hierarchical shape description
     \item Explicit neighbor connectivity enables shape deformation
 \end{itemize} \\ 
\hline
 Cons & \begin{itemize}[leftmargin=*]
    \item No geometry: no connectivity, no surface normal, no neighbors
     \item Sparse
     \item Hard to compare similarity and register shape
     \item No direct constraint for surface or volume 
     \item not unique
 \end{itemize} &
 \begin{itemize}[leftmargin=*]
    \item Irregular data structure
     \item Non-unique shape description
     \item Harder for learning applications
     \item Hard to compare similarity and register shape
     \item Ill-suited for simultaneously capturing multiple shapes
 \end{itemize} & \begin{itemize}[leftmargin=*]
    \item Irregular data structure
     \item Hard for learning applications
     \item Hard to compare similarity and register shape
     \item Ill-suited for simultaneously capturing multiple shapes 
     \end{itemize}\\ 
 \hline
\end{tabular}
\end{center}
\caption{Advantages and disadvantages of explicit shape representations.}
\label{table:pros_cons_explicit_representations}
\end{table}

Broadly, we identify four key considerations: 

\begin{enumerate}
    \item \textbf{Storage} is the computational footprint required to work with the shape representation in-memory. If a shape occupies only a fraction of an image or requires high-resolution details, it can be efficient to use explicit representations like points, meshes, and graphs that specifically capture the geometry, particularly with modern software libraries like PyG (\url{https://github.com/pyg-team/pytorch_geometric}) and Jraph (\url{https://github.com/google-deepmind/jraph}) supporting their batched processing. However, irregular structure may cause less efficient memory access, more scatter/gather operations, lower GPU utilization, and more overhead compared to points and images that can have a regular structure.     
    \item \textbf{Geometry} is the fidelity of the shape geometry that needs to be captured by the representation. If we are only interested in the extent of shape elongation, we do not require high fidelity, and very sparsely sampled points are sufficient. Meanwhile, if we are analyzing cell surface protrusions, high fidelity is essential, requiring microscopy images to be acquired with minimal pixel size \cite{sapoznik2020versatile}, and the use of meshes with a minimal number of vertices and faces \cite{zhou2023surface}.
    \item \textbf{Flexibility} refers to the ease with which a representation can be transformed into another representation to take advantage of algorithms developed for the other representation. Occupancy images are flexible. The grid connectivity permits a graph cut formulation, and they can be readily meshed to obtain point- and mesh-based representations.   
    \item \textbf{Scalability} is the ability of the representation to support time- and memory-efficient simultaneous processing of a set of shapes. Modern software libraries like PyG (\url{https://github.com/pyg-team/pytorch_geometric}) and Jraph (\url{https://github.com/google-deepmind/jraph}) can batch irregular graph-like representations. However, regular grid-connected representations like images are more efficient, taking advantage of hardware-optimized processing algorithms that exploit regularity like convolutions and Fast Fourier Transforms.
\end{enumerate}

In general, the choice of representation should align with the methods described in Section \ref{sec:methods}, depending on the shape analysis question being addressed.

\subsubsection{Converting between Shape Representations} 
Addressing shape problems in data often involves a series of analyses to be performed (see Section \ref{sec:practicalGuide} such as preprocessing: removing or imputing incomplete shapes due to inaccurate segmentation or partial acquisition; landmarking and alignment to establish correspondences between shapes locally; computing features to compare the correspondent shape regions and statistical testing to find significant differences). Consequently, it can be advantageous to interconvert between representations and maintain multiple equivalent representations to harness the advantages of individual representations for each operation \cite{zhou2023surface}. 

For example:  
\begin{itemize}
    \item Resampling point clouds to a target number of points requires defining a sampling space, either (i) estimated by surface reconstruction~\cite{bernardini1997alphawrap} or nearest neighbour connectivity, or (ii) provided by an occupancy image.
    \item Repairing meshes to simplify topology or for hole patching can be achieved more easily by constructing a volumetric representation either as a binary image~\cite{zeng2020cut,zhou2023surface} or as a tree of partitions~\cite{huang2018robust}.
    \item Reconstructing protein structure from cryo-electron microscopy (Cryo-EM) images involves both occupancy images and graphs. Occupancy images are used to train a neural network to efficiently predict the coarse-grained structure from the thousands of acquired images. Then, to recover the final fine-grained geometry, a structural graph initialized from homology to a similar protein based on amino acid sequence alignment is fitted to the coarse-grained structure using physical force constraints.
    
\end{itemize}

The specific algorithm implementations may vary but the core idea to convert between points, occupancy images, signed distance functions, meshes, and graphs is distilled in Table~\ref{table:converting_between_shape_formats}. The essence is that points correspond to image pixels, mesh vertices, and graph nodes. Occupancy images can be obtained by `densifying' sparse point clouds, e.g., by radially expanding the volume occupied by each point, or by thresholding of a signed distance function. Meshes are obtained by reconstructing the neighborhood connectivity of points by triangulating sparse points~\cite{bernardini1997alphawrap, lee1980two, kazhdan2006poisson, bernardini1999ball} or by meshing the isosurface, i.e., all pixels and voxels with the same value in image data~\cite{newman2006survey, hang2015tetgen}.  Graphs are obtained by constructing the adjacency matrix between mesh vertices or faces~\cite{feng2019meshnet}, grid connectivity of images~\cite{rue2005gaussian, zhou2025universal}, or by using an imposed connectivity structure between points such as a k-nearest neighbors graph~\cite{qi2017pointnet}.  

\begin{table}[t!]
  \centering
\begin{tabular}{ |p{0.12\textwidth}|p{0.1\textwidth}|p{0.12\textwidth}|p{0.15\textwidth}|p{0.2\textwidth}|p{0.15\textwidth}|}
\hline

\multirow{1}{*}{\textbf{From}}
&
 \multicolumn{5}{c|}{\textbf{To}} \\
 \cline{2-6}

   \rowcolor{gray!50}  &  Points & Occupancy Image & Signed Distance Function & Mesh & Graph \\
  \cline{1-6}
  \cline{1-6}
  Points & \centering -- & Raster onto image grid & From occupancy image & Triangulation & construct neighbor adjacency\\
  \hline 
  Occupancy Image & Coordinates & \centering -- & Distance metric & Isosurface meshing & Construct grid connectivity\\
  \hline
  Signed Distance Function & Coordinates & Threshold  & \centering -- & Isosurface meshing & Construct grid connectivity \\
  \hline
  Mesh & Vertices & Voxelization & From occupancy image & \centering -- & Convert triangle and tetrahedral connectivity to vertex adjacency\\
  \hline
  Graph & Coordinates & Raster onto image grid * & From occupancy image & Isosurface occupancy image or reconstruct from points & -- \\ 
  \hline
\end{tabular}
\caption{Converting between different representations of a shape. Conversion is not always exact, but is performed to simplify computations or better visualize particular shape features with a different representation. Results computed in the target representation can always be remapped back to the source representation.}
\label{table:converting_between_shape_formats}
\end{table}

\section{Shape Space Overview}
\label{sec:shape_space_overview}

Across disciplines, \textit{shape} provides a fundamental description of structure and function---whether in anatomical analysis, molecular modeling, or computer vision. The study of \textit{shape spaces} offers a mathematical framework for reasoning about geometric variability. At its core, shape analysis seeks to address three central questions:

\begin{enumerate}
\item \textbf{What is the variation of shapes within a dataset?} \
Addressed via metrics (Section \ref{sec:distance}), feature design (Section \ref{sec:feature_descriptors}), dimensionality reduction, and generative modeling.

\item \textbf{How does one shape relate to another, within or across datasets?} \
Addressed through registration and alignment techniques that establish point correspondences between shapes (Section \ref{sec:alignAndRegister}).

\item \textbf{How can we describe or map specific locations on or within a shape?} \
Addressed through parameterization methods that define intrinsic coordinate systems (Section \ref{sec:surfaceUnwrapping}).
\end{enumerate}

All of these problems require defining a \textit{space of shapes} - a mathematical domain in which each shape corresponds to a point, and distances, averages, and transformations can be meaningfully defined. This notion underlies nearly every computational framework for shape analysis.

Representing shapes in a mathematically consistent way, however, is far from trivial. Meaningful analysis requires a formulation that identifies all representations corresponding to the same intrinsic geometry. Mathematically, this motivates treating shape as an equivalence class under an appropriate family of transformations. In practice, shape analysis may operate either directly on geometric representations or on abstract feature-based encodings tailored for statistical or machine learning-based applications.

\subsection{Representations of Shapes} \label{sec:shape_representation}

Shape representations can be broadly categorized according to how they encode geometric variability. While the quotient-space framework introduced below provides an abstract characterization of shape, practical analysis requires concrete representations that support storage, comparison, and computation. Three complementary families of representations are commonly used:

\begin{enumerate}
\item \textbf{Function Spaces of Parameterizations}, where shapes are represented by maps from a template manifold into Euclidean space, enabling spatially consistent comparisons via explicit coordinates or basis coefficients.

\item \textbf{Measure-Based Representations}, which encode shapes as probability or mass distributions, allowing flexible comparison across differing topologies and sampling densities.

\item \textbf{Feature-Based Representations}, which summarize shapes through geometric, appearance-based, or relational descriptors-hand-crafted or learned.
\end{enumerate}

\begin{figure}[t]
\includegraphics[width=1.0\linewidth]{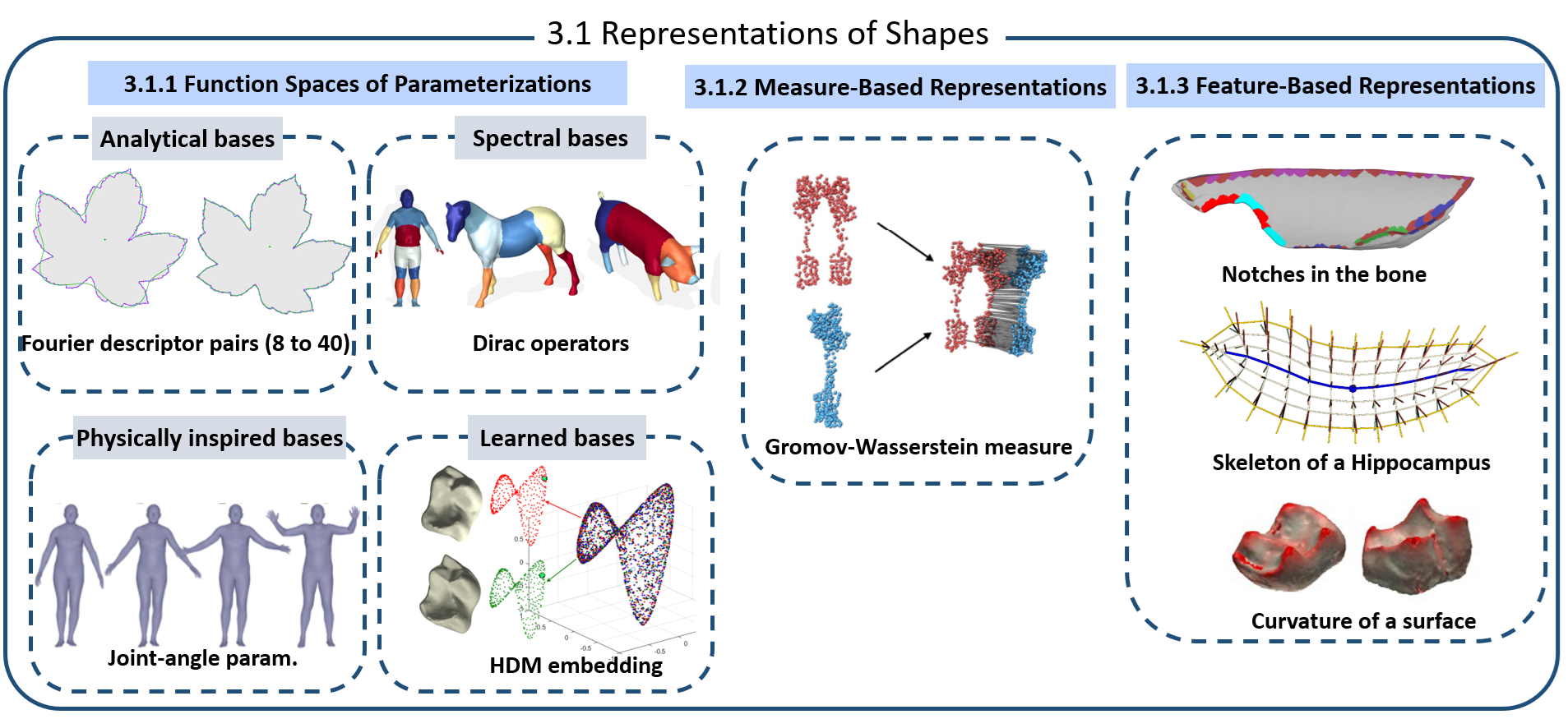} 
\caption{\textbf{Key ideas in Representations of Shapes as discussed in Section~\ref{sec:shape_representation}.} Specifically, we discuss the Function Spaces of Parameterizations, the Measure-Based Representations, and the Feature-Based Representations (images adopted from ~\cite{Burger2013, liu2017dirac, hartman2024varishape, faigenbaum2026studying, tajmir2025alignment, yezzi2022using, pizer2020object, shan2019ariadne}).}
\label{fig:Shape_Representations}
\end{figure}

These paradigms form the foundation for most computational and theoretical approaches to shape analysis. The following subsections summarize their motivations, computational tools, and practical trade-offs (see Fig.~\ref{fig:Shape_Representations} for an overview).

\subsubsection{Function Spaces of Parameterized Shapes}

A natural strategy is to represent shapes as parameterized objects, modeled as subsets of a function space of bijective maps from a template domain into Euclidean space~\cite{floater2005surface,sheffer2006mesh,hormann2007mesh,gu2008computational,crane2017glimpse,gu2020computational}. For example, open surfaces may be represented by embeddings of disk-like domains; closed surfaces or volumes by mappings from a spherical domain; and tubular structures by cylindrical domains with centerline-based coordinates. A key limitation of this approach is the requirement that all shapes share the same topology as the template, making it unsuitable for datasets with topological variability. Nevertheless, explicit parameterizations are widely used in biological and medical settings, where datasets often consist of repeated observations of the same anatomical structure, and precise location indexing everywhere in the shape is desired for comparative analyses. 

An alternative is to represent shapes through coefficients in a chosen basis of the space of parameterizations, providing compact encodings and efficient optimization better suited for machine learning in high-dimensional function spaces. Common choices include:
\begin{itemize}
\item \textbf{Analytical bases} such as Fourier modes or spherical harmonics (e.g.~\cite{Burger2013}),
\item \textbf{Spectral bases} derived from Laplace--Beltrami or Dirac operators (e.g.~\cite{liu2017dirac}),
\item \textbf{Physically inspired bases} such as joint-angle parameterizations for human poses (e.g.~\cite{hartman2023varigrad,hartman2024varishape}),
\item \textbf{Learned bases} obtained from PCA, NMF, Diffusion maps, HDM, sparse coding, or deep generative models (e.g.~\cite{kosambi2017statistics, faigenbaum2026studying}).
\end{itemize}

These basis coefficients capture multiscale geometric variation and support learning-based tasks such as shape-based pattern discovery, clustering, classification, and retrieval. However, individual basis elements may be harder to associate with distinct geometrical features. In cell biology, spectral bases have been used, for example, to project substructures (e.g., mitochondria, surface protrusions, or molecular markers) into low-dimensional coordinates for quantifying spatial organization. Each coefficient captures a different spatial scale, but no single coefficient is associated with a single physical structure.

While parameterizations provide a convenient coordinate system, they are not invariant to pose, scale, or the choice of parameterization itself. We obtain representations that depend only on \emph{shape} by quotienting out appropriate transformation groups. Let $\mathcal{I}$ denote a space of parameterized objects (e.g., landmark configurations, curves, or surfaces). Let $\mathcal{D}$ be the group of reparameterizations acting on $\mathcal{I}$, and let $\operatorname{SE}(d)$ (optionally extended with scalings $\mathbb{R}^{>0}$) represent rigid motions. The associated \textbf{shape space} is the quotient
\begin{equation}
\mathcal{S} = \mathcal{I} \big/ \bigl( \mathcal{D} \rtimes (\operatorname{SE}(d) \times \mathbb{R}^{>0}) \bigr),
\end{equation}
with each shape represented by the equivalence class
\begin{equation}
[q] = \left\{\, s R (q \circ \gamma) \;\middle|\; \gamma \in \mathcal{D},\; R \in \operatorname{SE}(d),\; s \in \mathbb{R}^{>0} \,\right\}.
\end{equation}

Analyzing shapes in this quotient space ensures invariance to reparameterization, rigid motion, and, where desired, scale, allowing application-driven choices of which transformations to remove.

Two primary examples illustrate this model:
\begin{enumerate}
\item \textbf{Landmark shapes.}  
Here $\mathcal{I} = (\mathbb{R}^d)^n$ is the space of ordered landmark configurations, and $\mathcal{D} = S_n$ accounts for landmark permutations~\cite{kendall84,kendall2009shape}. Quotienting by rigid motions and scalings yields Kendall's classical shape spaces.

\item \textbf{Curve and surface shapes.}  
Shapes of parameterized curves or surfaces correspond to
\begin{equation}
\mathcal{I} = \operatorname{Emb}(M,\mathbb{R}^d)
\quad\text{or}\quad
\mathcal{I} = \operatorname{Imm}(M,\mathbb{R}^d),
\end{equation}
with reparameterization group $\mathcal{D} = \operatorname{Diff}(M)$ acting on the domain~\cite{Michor2006,bauer11,SrivastavaKlassen2016}. Quotienting yields spaces of unparameterized shapes, which are central to elastic shape analysis, diffeomorphic registration, and equivariant learning.
\end{enumerate}

\subsubsection{Measure-Based Representations}
Measure-based representations provide a powerful and unifying framework for describing geometric data such as point clouds, curves, surfaces, and volumetric shapes. Instead of treating these objects as explicit parameterizations or discrete approximations of parameterized objects, the idea is to encode them as measures-generalized distributions that capture spatial location, density, and often tangent or orientation information. This perspective originates in geometric measure theory, where geometric objects of arbitrary regularity (smooth, piecewise-smooth, or highly irregular) can be described using Hausdorff measures, Radon measures, or density fields. 

Beyond these basic constructions, richer measure-theoretic models have been developed to encode additional geometric structure. Currents provide an oriented representation of shapes by integrating vector fields over a manifold using differential forms; they are effective for comparing shapes with consistent orientation but fragile when orientation flips occur. Varifolds, by contrast, drop the orientation requirement and instead encode geometry as a measure over positions and tangent spaces. This makes them well-suited for applications such as medical image analysis, shape registration, and data-driven geometric learning where noise, topology changes, and partial observations occur. Both frameworks allow shape comparison through dual norms or kernel embeddings, enabling robust metrics even when the underlying discretization differs.

More recent developments link measure representations to optimal transport (OT) and kernel methods, which provide computationally efficient and differentiable tools for comparing geometric measures. Wasserstein distances, sliced OT, kernelized varifolds, and hybrid current–OT models have become popular in computational anatomy, shape matching, and geometric machine learning. These approaches combine the flexibility of measure-based representations with tractable numerical methods, enabling shape comparison, registration, and statistical modeling in high-dimensional or irregular geometric datasets. Overall, the measure-theoretic viewpoint offers a flexible, mathematically principled, and computationally compatible way to represent and analyze geometric data of widely varying structure.

\subsubsection{Feature-Based Representations}

Feature-based representations describe shapes through vectors of geometric, appearance-based, or relational descriptors rather than through explicit spatial correspondences. This approach is especially useful when handling heterogeneous shapes, when invariance to pose or parameterization is desirable, or when high-level semantics outweigh the need for precise geometry.

Classical methods compute descriptors such as curvature statistics, shape contexts, spin images, appearance histograms, or motion features. Once normalized, these feature vectors support exploration, clustering, dimensionality reduction, and other tasks in which shape variability may be subtle or complex.

Modern approaches use neural networks-such as autoencoders, graph neural networks, contrastive-learning frameworks, or functional-map-based architectures-to learn shape embeddings, and metrics. The embeddings capture geometric and semantic similarity and can be designed to be invariant to shape-preserving transformations. The learnt metrics capture a measure of distance between two embeddings, reflective of their geometric and semantic difference. Generative models further enable sampling or interpolation of shapes conditioned on continuous or categorical variables, supporting applications in cell biology, development, and medical diagnosis where quantifying heterogeneity is essential to establish the statistical significance of a treatment condition.

\section{Methods}
\label{sec:methods}
\subsection{Pre-Processing Steps}

Once the shapes are acquired, analysis often cannot begin immediately. The acquired data typically contains noise, and there may be missing information. Fixing these issues requires pre-processing steps to reconstruct well-formed shapes suitable for computational analysis. Smoothing, local parametric fitting or surface reconstruction often fixes noisily sampled points and patches small holes. These steps can be performed semi-manually (e.g., using off-the-shelf software GeoMagic), automatically (e.g., MeshFix), and by building customized pipelines using software libraries like the Point Cloud Library (PCL) and Open3D. For some shape analysis methods, such as shape parameterization (Section \ref{sec:surfaceUnwrapping}), the shape to be analyzed must have specific properties such as being genus-0. In these cases, pre-processing steps like topological simplification and shrinkwrapping are applied not to clean the data, but to generate a representation of the shape with the required properties.  

The main challenge arises when the shape has large incomplete or damaged patches, or comes in several pieces. In such cases, substantial effort is needed to first reconstruct the missing or broken parts before the shape can be studied. This fragment assembly challenge has been framed as solving a three-dimensional jigsaw puzzle in~\cite{grim2016automatic} or in the vein of multi-piece matching and template comparison in~\cite{zhang20153d}. The challenge commonly arises in practice in the reconstruction of archaeological specimens and in architectural conservation~\cite{schafer2015virtually}. Equally challenging is when the 3D shape can only be partially captured from different angles or low-dimensional views. The former commonly arises in the digitization of 3D scenes in computer vision and the latter, in protein structure determination from 2D cryo-electron microscopy (cryo-EM) images. In these cases, camera parameters and Euclidean geometry are used to map the shape data into a global coordinate system. This often involves also the difficult tasks of determining which data views are associated with the same object, and finding a global optimization that minimizes the discrepancy of the reconstructed shape with all views.

\subsection{Shape Correspondence}

Shape correspondence is a fundamental problem that encompasses key processes such as alignment and registration, often by means of landmarking, to establish meaningful pointwise relationships between shapes. \textit{Alignment} ensures that shapes are positioned in a common coordinate system, minimizing differences due to translation, rotation, or scaling. \textit{Landmarking} involves identifying specific, anatomically or structurally significant points on each shape to serve as reference markers for comparison. \textit{Registration}, often performed using rigid or non-rigid transformations, refines the correspondence by warping one shape onto another while preserving geometric consistency. Together, these techniques enable precise shape analysis and serve as a first step to the comparative analysis that follows. Often, these three come hand in hand together, and enhancing one will improve the others.

\label{sec:alignAndRegister}

\subsubsection{Landmarking}
\label{sec:Landmarking}
Landmarks are often a fundamental step that many methods in the shape analysis pipeline later rely on (starting from registration, comparison, and pattern analysis). 
Landmarks are structurally (or anatomically) correspondent, often termed as ``homologous'' in literature, points on a shape that serve as reference markers for comparison.  Although these are interesting shape features in their own right, they are also frequently an essential step in comparative shape analysis. They provide a structured and interpretable framework for comparing and analyzing shapes across different specimens, time points, or conditions. By serving as key references for shape correspondence, landmarks ensure that meaningful geometric features are consistently aligned. They play a vital role in geometric morphometric analysis (GMM) (both manual and automatic), quantifying shape variations, identifying homologous structures, and enabling statistical shape modeling (see Section~\ref{sec:Statistical Analysis}). Additionally, landmarks facilitate registration and deformation analysis, allowing for comparisons that extend beyond simple geometric transformations. 

\begin{figure}[t]
	\centering
    \includegraphics[width=\textwidth]{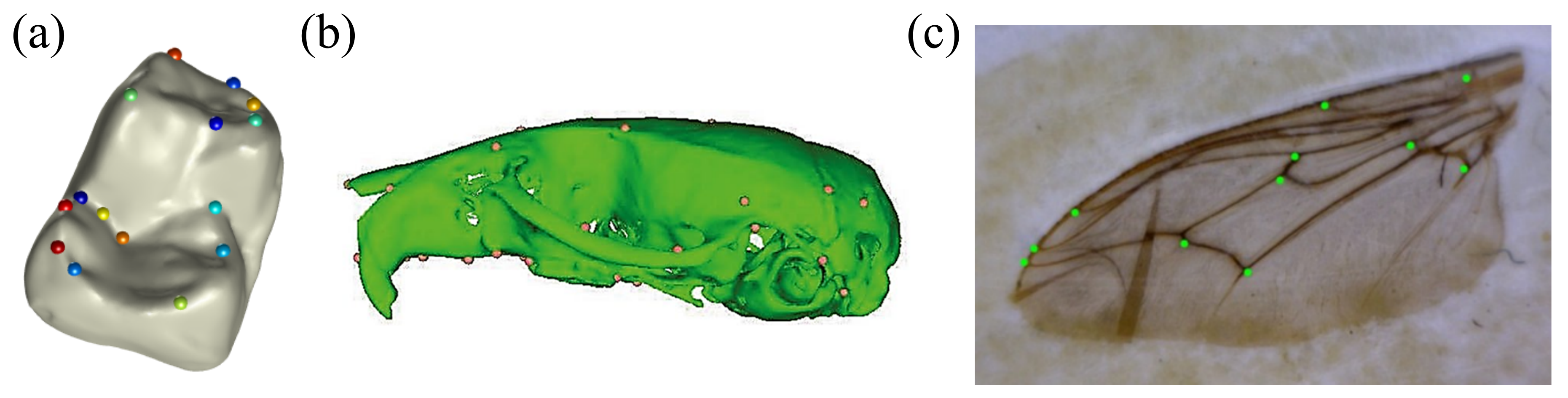}   
	    
	\caption{\textbf{Landmarking in different applications and different data representations (3D surfaces and images).} (a)~Gaussian Process Landmarking of a primate teeth~\cite{gao2019gaussian}. (b)~Transferring the landmark position of a single specimen into another via ALPACA~\cite{porto_alpaca_2021} on a laboratory mouse skull. (c)~Convolutional neural network detection on a fly wing~\cite{geldenhuys2023deep}. }
	\label{fig:landmarking}
\end{figure}

One common approach is to select landmarks that capture key geometric features of the surfaces. In expert-driven landmarking, domain specialists manually identify an equal number of geometrically or semantically meaningful feature points on each surface, typically using computer-based tools, and certify them as being in consistent one-to-one correspondence. While experts are generally effective at identifying salient regions, such as areas of high or low curvature, defining reliable correspondences in intermediate regions between these features is more challenging. Thus, automated procedures for identifying important points were developed. 

Herein, we cover several automatic landmarking methods: Continuous Procrustes landmarking~\cite{al2013continuous} provides pairwise correspondences based on minimizing pointwise distances but lacks global consistency across collections; the Matching Pursuit algorithm~\cite{gao2015hypoelliptic} leverages spectral sparsity in a hypoelliptic Laplacian to extract landmarks across an entire shape collection. In contrast, Gaussian Process Landmarking~\cite{gao2019gaussian, gao2019gaussian_2} selects landmarks sequentially by maximizing prediction uncertainty based on surface geometry, producing biologically meaningful and semantically consistent features (see Fig.~\ref{fig:landmarking}(a)). In \cite{vandaele2018landmark}, the authors use a multi-resolution tree-based approach. Deep learning methods have also been applied to landmark detection. For example, deep neural network has been applied to identify landmarks on Drosophila wings \cite{geldenhuys2023deep} (see Fig.~\ref{fig:landmarking}(c)). Thomas and Maga adapted the MorphVQ framework for automated landmarking~\cite{thomas_leveraging_2024}. Originally developed for landmark-free shape analysis within the functional map framework~\cite{sun_spatially_2023, thomas_automated_2023}, MorphVQ was extended to establish point correspondences among meshes and to select discrete landmarks. In medical imaging, deep learning is frequently used to predict landmark coordinates directly in voxel space, such as in cephalometric analysis~\cite{polizzi2024automatic}.

Another class of approaches derives landmarks through model registration and deformation. ALPACA is a lightweight pipeline that propagates landmarks from a template based on rigid and coherent-point-drift registrations of downsampled point clouds~\cite{porto_alpaca_2021, zhang_automated_2022} (see Fig.~\ref{fig:landmarking}(b)). Dense Correspondence Analysis (DeCA) refines registration-based sparse correspondence using manually placed landmarks and subsequent semilandmark sampling to produce denser representations~\cite{rolfe_deca_2023}.

Last, let us ask: \textit{How many landmarks are sufficient?} This is a one-million-dollar question. In \cite{watanabe2018many}, the authors provide a thorough study of this issue. On the other hand, Mitteroecker and colleagues still recommended starting from denser landmarks for exploration and visualization, even though there is a risk of increasing redundancy and arbitrary weighting of different signals~\cite{mitteroecker2022thirty}. 

Despite their widespread use, robust and reproducible landmark identification remains challenging. Moreover, representing shapes using a finite set of discrete landmarks may fail to capture aspects of continuous geometric variation.  Consequently, landmarking is not strictly required for all shape analysis approaches, and alternative strategies based on direct correspondence estimation are discussed in the following subsections.

\subsubsection{Alignment}

Shape alignment refers to placing shapes into a common coordinate system by removing extrinsic variability such as translation, rotation, and sometimes global scaling, so that shapes can be meaningfully compared without altering their intrinsic geometry. Noteworthy that shape alignment and shape registration are closely related but conceptually distinct tasks in shape analysis. 

Commonly, the first step of finding the best correspondence is alignment. The goal of shape alignment is to establish meaningful correspondences between points, features, or structures across different shapes, ensuring that they are spatially aligned or registered. This process must be both efficient and robust. Given the high number of points on each shape, downsampling is often performed prior to analysis to reduce computational complexity. It is also important to note that finding a common coordinate system is generally more straightforward for similar shapes than for dissimilar ones. This principle is incorporated into the Automated 3D Geometric Morphometrics (Auto3DGM) tool. Auto3DGM begins by using Principal Component Analysis (PCA) to compute the three principal axes for each shape~\cite{boyer2015new}. Other methods also rely on a PCA for this task, e.g.,~\cite{styner2003evaluation}.

\subsubsection{Registration}
\label{sec:Registration}

Registration is a fundamental pre-processing step while working with data, and is often performed after an initial alignment. Without registering the shapes its hard to talk about their comparison. Since registration alone can fill a whole survey paper (e.g. \cite{van2011survey, sahilliouglu2020recent,DynReg}), we will mention only a few methods. Shape registration can be approached through several complementary frameworks. \textit{Rigid registration} aligns shapes using only rotations and translations (and sometimes scaling), as in Iterative Closest Point (ICP). \textit{Non-rigid registration} extends this by allowing deformations such as bending and stretching while preserving intrinsic structure. \textit{Feature-based matching} shifts the focus from raw geometry to descriptors, enabling more robust correspondence. \textit{Spectral or functional methods} move further by representing correspondences as linear operators between function spaces-for example, in Deep Functional Maps-rather than explicit point-to-point matches. \textit{Probabilistic methods} formulate registration as statistical inference, as in Coherent Point Drift (CPD). More recently, \textit{learning-based approaches} leverage data to learn correspondences directly, including point cloud networks (e.g., PointNet and Point Transformer-style architectures), graph neural networks, and functional map pipelines with learned features. Finally, \textit{optimal transport-based methods} model matching as transporting mass between distributions, often using distances such as the Wasserstein metric.

Procrustes superimposition and other least-squares global registration methods, which belong to the class of  \textit{rigid (or similarity) registration}, have become standard procedures in geometric morphometrics. However, several studies have shown that these approaches can influence statistical results and visual interpretations in subjective and sometimes unpredictable ways~\cite{bookstein_reworking_2023, cardini_integration_2019, cardini_procrustes_2022, richtsmeier_promise_2002}. In particular, simulation studies demonstrated that Procrustes alignment may introduce spurious patterns of integration and modularity~\cite{goswami_high-density_2019, cardini_integration_2019}. More generally, global registration methods can obscure localized shape variation through the so-called “Pinocchio effect,” whereby changes concentrated in a small subset of landmarks propagate across the entire configuration~\cite{klingenberg_how_2021, klingenberg_visualizations_2013}. Additionally, arbitrary rotational components introduced during alignment may further complicate visualization and statistical interpretation~\cite{richtsmeier_promise_2002}.

Moving beyond rigid alignment, \textit{non-rigid registration} methods, such as thin-plate splines, aim to model smooth deformations between shapes. These approaches are closely related to classical interpolation and elasticity-based models, and can be extended to more general frameworks such as as-rigid-as-possible deformations and deformation graphs. Nevertheless, they remain sensitive to modeling assumptions. Recently, Fred L. Bookstein questioned several Procrustes-based conventions, including the reliance on homologous landmarks, Procrustes superimposition, and thin-plate spline visualization~\cite{bookstein_reworking_2023}. Using classical cranial growth data, he demonstrated that baseline registration–derived polynomial mappings can yield smoother and more interpretable deformation patterns, emphasizing the critical role of registration choices in exploratory analysis.

Beyond these classical approaches,  \textit{feature-based matching} methods form a central paradigm in shape registration. These methods rely on extracting local or global descriptors from each shape and establishing correspondences based on feature similarity, often followed by robust estimation procedures such as RANSAC. In addition to handcrafted descriptors, recent work also considers learned descriptors, bridging toward data-driven methods. Early work in this direction is well documented in surveys on shape matching and registration~\cite{van2011survey,DBLP:journals/tvcg/TamCLLLMMSR13}, and continues to underpin many modern pipelines.

A complementary line of work adopts a  \textit{probabilistic perspective}, where registration is formulated as statistical inference between point sets or shapes. Classical approaches such as the Coherent Point Drift model one shape as a probability distribution and align it to another by minimizing distributional discrepancy. More generally, Markov Random Field (MRF) formulations model pairwise compatibility between candidate matches and enable the incorporation of constraints such as non-penetration. Originally introduced for problems such as jigsaw puzzle solving~\cite{DBLP:conf/cvpr/ChoBAF08,DBLP:journals/pami/ChoAF10}, these formulations were later applied to structure-from-motion and related tasks~\cite{crandall2012sfm}. Their advantage lies in principled optimization frameworks, such as linear programming relaxations, which often outperform greedy strategies, although they typically require a large number of candidate poses.

A major development in the past decade is the emergence of spectral and functional methods, in particular the  \textit{functional map} framework~\cite{ovsjanikov_functional_2012}. The central idea is to model transformations between shapes as mappings between function spaces, rather than point-to-point correspondences. This representation enables compact formulations and facilitates the incorporation of regularization, such as enforcing commutativity with the Laplace–Beltrami operator. Closely related are diffusion-based and spectral embedding methods, which leverage intrinsic geometry for robust matching. Functional maps also provide a convenient algebraic structure that is well-suited for integration with machine learning. In particular, state-of-the-art isometric shape matching builds on Deep Functional Maps, which combine feature extraction on each shape with a matching module in the functional map space. Training can rely on labeled correspondences or geometric distortion losses, and subsequent work has extended the framework by learning suitable bases, including for point clouds.

Recent progress has increasingly shifted toward \textit{learning-based methods}, where correspondences are inferred directly from data. These include deep functional map pipelines, as well as architectures based on graph neural networks and point cloud networks. More recently, Transformers have been adapted to geometric data, where attention mechanisms incorporate spatial structure and point locations~\cite{DBLP:journals/cgf/RaganatoPM23}. In parallel, foundation models have been leveraged for 3D tasks: methods such as STAR, 3D Highlighter, and follow-up works employ vision-language models, including CLIP, GLIP, and Segment Anything Model, by lifting 2D features to 3D representations. This paradigm enables transferring large-scale pretrained knowledge from images to shape matching and segmentation. Additionally, generative approaches, such as GenCorres, learn distributions over correspondences, while other methods explicitly learn deformation fields or canonical shape spaces.

Finally, moving beyond pairwise matching, recent work emphasizes \textit{global and multi-shape registration}, closely related to ideas from optimal transport and synchronization. In particular, optimal transport formulations model correspondences as soft maps between distributions, enabling robustness to partial matching and noise. The key insight in multi-shape settings is that correspondences should be consistent across a collection of shapes. One approach constructs paths or trajectories between shapes, where composing reliable correspondences between intermediate shapes yields improved results and enforces cycle-consistency constraints~\cite{DBLP:conf/cvpr/ZachKP10,DBLP:journals/cgf/NguyenBWYG11,DBLP:journals/tog/HuangZGHBG12}. Alternative methods optimize over the entire collection, for example by constructing spanning trees or by formulating the problem as low-rank matrix recovery, exploiting the fact that the matrix of pairwise maps is approximately low-rank~\cite{DBLP:journals/cgf/HuangG13}. These formulations can be solved using semidefinite programming, spectral techniques, or non-convex optimization, offering strong empirical performance and theoretical guarantees, although incorporating certain geometric constraints, such as non-penetration, remains challenging.

\subsection{Shape Parameterization}
\label{ref:shape_param}
Shape parameterization typically refers to the construction of a local coordinate system to index any position on a shape’s surface (surface parameterization) or in its volume (volume parameterization)~\cite{sheffer2006mesh,hormann2007mesh}. This is an important topic as the coordinates can generate consistent landmarks across a collection of shapes to enable numerous applications including alignment, comparison, and deformation. A primary use is to integrate out shape variations to compare substructures within the shape. In medicine, a common coordinate framework provided by a template brain enables longitudinal monitoring of disease in brain regions of a single patient or cross-patient comparison for subtyping disease independent of brain shape variation. \\

\subsubsection{Surface Parameterization}
\label{sec:surfaceUnwrapping}

Surface parameterization aims to map a complicated surface onto a simpler domain, thereby reducing the dimensionality and complexity of the shape for subsequent analyses. It has numerous applications in geometry processing, shape registration, and shape analysis. Over the past few decades, many surface parameterization methods have been developed, covering different types of surfaces, choices of the target parameter domains, and the mapping criteria. For overviews and surveys, see~\cite{floater2005surface,sheffer2006mesh,hormann2007mesh,gu2008computational,crane2017glimpse,gu2020computational}. Note that in general, isometric mappings between the given surfaces and the target parameter domains do not exist. Therefore, parameterization methods usually focus on reducing a specific type of geometric distortion while trading off some other distortions. It is also noteworthy that the chosen parameterization is usually closely tied to the geometry of the data. Due to the subtle differences between shapes, a parametrization that works well for one dataset may not be suitable for another, as it is typically designed to capture the specific characteristics of the data under study.

Among different parameterization approaches, conformal parameterizations have been extensively studied and utilized for various applications. Specifically, conformal parameterizations are angle-preserving, and hence the local geometry of a surface as well as the toolset of differential calculus are well-preserved under conformal maps. There are a vast number of conformal parameterization techniques, including linearization~\cite{haker2000conformal}, least-squares conformal map~\cite{levy2002least,desbrun2002intrinsic}, harmonic energy minimization~\cite{gu2004genus}, spectral conformal map~\cite{mullen2008spectral}, Ricci flow~\cite{jin2008discrete,yang2009generalized,zhang2014unified}, curvature prescription~\cite{ben2008conformal}, quasi-conformal composition~\cite{choi2015flash,choi2015fast,choi2021efficient}, boundary first flattening~\cite{sawhney2017boundary}, conformal energy minimization~\cite{yueh2017efficient}, and conformal welding~\cite{marshall2007convergence,choi2020parallelizable}.

Area-preserving parameterizations are another class of wrapping techniques that have also been widely studied. Unlike the conformal parameterizations, the area-preserving parameterizations preserve the proportionate area of the area elements but not the local angles. Some notable existing area-preserving parameterization techniques include Lie advection~\cite{zou2011authalic, zhou2023surface}, optimal mass transportation~\cite{zhao2013area,su2016area,giri2021open,choi2022adaptive}, density-equalizing maps~\cite{choi2018density,choi2020area,zhou2023surface,lyu2024spherical}, and stretch energy minimization~\cite{yueh2019novel,yueh2023theoretical}.

Besides conformal and area-preserving parameterizations, there are also many parameterization methods based on other mapping measures that achieve a balance between conformal and area distortions. Examples include quasi-conformal map~\cite{zeng2009surface,weber2012computing,lam2014landmark,lui2014teichmuller}, bounded distortion map~\cite{lipman2012bounded,aigerman2013injective,chen2015bounded,kovalsky2015large,chien2016bounded}), and other balancing schemes~\cite{nadeem2016spherical,claici2017isometry,choi2025hemispheroidal,zhou2023surface}.

Another important aspect of surface parameterization is the preservation of local injectivity or bijectivity. Specifically, it is desirable that surface parameterization results do not contain any mesh fold-overs or inverted elements. Numerous methods have been proposed for achieving locally injective/bijective parameterization via optimization, such as the advanced MIPS~\cite{fu2015computing}, bijective parameterization with free boundaries~\cite{smith2015bijective}, scalable locally injective mappings~\cite{rabinovich2017scalable}, progressive parameterization~\cite{liu2018progressive},
efficient bijective parameterization~\cite{su2020efficient} (see also the survey paper~\cite{fu2021inversion}). It is also possible to ensure local injectivity/bijectivity via quasi-conformal theory~\cite{choi2015flash,lyu2024bijective} (see also the survey paper~\cite{choi2023recent}).

The topology of the given surface also plays an important role in the computation of surface parameterization. For simply connected open surfaces, one may parameterize them onto a planar disk~\cite{jin2008discrete,yueh2017efficient,choi2018linear}, a rectangle~\cite{lam2014landmark,meng2016tempo}, or some more flexible domains~\cite{levy2002least,mullen2008spectral,sawhney2017boundary}. For multiply connected open surfaces, it is common to parameterize them onto a planar domain with circular holes~\cite{zeng2009generalized,zhu2022parallelizable} or a domain with slits~\cite{yin2008slit}. For genus-0 closed surfaces, a popular choice of the parameter domain is the unit sphere~\cite{haker2000conformal,gu2004genus,zhou2023surface}. For high-genus surfaces, common approaches include partitioning the surfaces and mapping different patches individually~\cite{jin2004optimal} and parameterizing the entire surface onto a planar domain with periodic boundary conditions~\cite{zhang2012canonical,lui2014geometric,yao2026toroidal}.

In recent years, there has been an increasing interest in utilizing neural networks for surface representation and  parameterization~\cite{sinha2016deep,vakalopoulou2018atlasnet,FoldingNet,bednarik2020shape,chen2022auv}. For instance, Neural Surface Maps~\cite{morreale2021neural} considered neural networks as a parametric representation of surfaces and inter-surface maps and used them for surface parameterization. Minimal Neural Atlas~\cite{low2022minimal} used an implicit probabilistic occupancy field to model the parametric domain for parameterization. Nuvo~\cite{srinivasan2024nuvo} used a neural field to represent a continuous UV mapping. The Flatten Anything Model~\cite{zhang2024flatten} achieved surface parameterization using an unsupervised neural architecture. A key advantage of neural network-based parameterizations is their ability to be integrated into end-to-end architectures along with other tasks such as the learning of shape features and classification~\cite{de2025geometric}.

\subsubsection{Volume Parameterization}
\label{sec:volume_parameterization}
Volumes can be parameterized analogous to surfaces. The simplest approach involves drawing parallels to standard coordinate systems, particularly the sphere and cylinder. The sphere can be used for all star-convex shapes, whereby the shape surface coordinates are projected onto the sphere surface by unit length normalization of its displacement from the centroid. Then, linear interpolation is used to sample coordinates along the line between the centroid to each surface coordinate and map these proportionately into the sphere volume~\cite{Viana2023}. The cylinder can be used for tube-like shapes such as vessel fragments and esophagus, where a centerline can be computed to define distance along the centerline, distance from the centerline, and rotation angle. Data fitting splines can also be used to extend the spherical and cylindrical frames~\cite{heemskerk2015tissue,stower2023single}. 

More complex shapes can be represented by tetrahedral meshes and parameterized directly by volumetric mapping approaches, such as harmonic energy minimization~\cite{wang2003volumetric}, optimal mass transportation~\cite{su2016measure,su2017volume}, stretch energy minimization~\cite{jin2015stretch,yueh2019novel2,yueh2020new,tan2025dimensional}, simplicial foliation~\cite{campen2016bijective}, star decomposition~\cite{hinderink2024bijective}, or other techniques with different geometric distortion criteria~\cite{nieser2011cubecover,su2019practical,pan2020volumetric,abulnaga2021volumetric,choi2021volumetric}. A key problem of tetrahedral meshes is the high memory consumption and time involved in solving this much larger set of sparse linear equations. Alternatively, the 3D surface is first parameterized, and then either interpolated, e.g., by harmonic extension to the centroid or a second parameterized 3D surface, or deformed under regularization such as elastic priors to sample the interior volume~\cite{zhou2023surface}. The volume is thus parameterized by a coordinate system comprising the deformed distance and the parameterized surface coordinates. Harmonic extension has notably been applied to map brain volumes. A notable class of such parameterizations in medical imaging for anatomic objects are skeletal models or S-reps~\cite{siddiqi2008medial,pizer2020object}. Deformation need not cover the entire volume. Multi-layer surfaces are used to visualize and measure cell shape parameters proximal to curved embryo surfaces in developmental biology~\cite{heemskerk2015tissue}. The surface to deform can also be judiciously chosen to highlight specific shape features. In cell biology, the deformation of a reference surface was used to construct a topography volume space to specifically isolate cell surface protrusions for analysis~\cite{zhou2023surface}. Temporally indexed volumes can be parameterized with a single set of coordinates of a reference time point after surface-surface or volume-volume alignment. Alternatively, each time point can be parameterized identically with the same coordinate construction.

\subsection{Shape Feature Descriptors}
\label{sec:feature_descriptors}

A shape can also be described by a multidimensional feature descriptor -- a vector concatenating different measurements of a shape's properties. Unlike a shape parameterization, which must faithfully reproduce the physical geometry and may be computable only for particular shape categories and shapes with a single spatial component, a feature could be any scalar measurement and be of different variable types -- categorical, discrete, or real-valued, though using measurement types other than continuous real values may affect the processing of features for downstream applications. Example scalar measurements include area and volume; the number of neighbors, holes, and spatial components; average area of neighboring objects; and histogram of chordal lengths. Thanks to this flexibility, which allows for quantitative description of shapes with arbitrary geometry and topology, and a more compact, information-dense representation compared to geometric coordinates, feature descriptors are the primary workhorse of computer graphics, machine learning, and scientific discovery. Feature descriptors are widely used to analyze and discover patterns in (i) heterogeneous shape collections for example to enable clustering techniques to identify distinct shape classes; (ii) construct latent or low-dimensional spaces to compare shape variation and quantify dynamics; (iii) develop generative shape models; and (iv) multimodal integrative analysis of a shape with other associated data including genetic and chemical characteristics and capturing properties and connectivity with other shapes in a spatial neighborhood.

Shape feature descriptors have been extensively reviewed in earlier works. See surveys on descriptor- and feature-based shape representations in 2D~\cite{zhang2004review, mingqiang2008survey}, in 3D~\cite{kazmi2013survey, valizadeh2022Fourier}, and in data-driven approaches~\cite{rostami2019survey}. Here, we highlight prevalent approaches and document accessible software implementations of these in the dedicated GitHub repository: \url{https://github.com/shirafaigen/ShapeSpaceSurvey.git}.

Any numerical measurement can be used as a feature, and multiple features can be concatenated into a feature descriptor. Since there are countless ways to define and combine features, it is impossible to review every type of shape descriptor. Nevertheless, shape descriptors are designed for three primary purposes: (i) measuring specific shape properties of interest, often for a particular object type such as worms or human posture; (ii) comparing and aligning different shapes or parts of shapes; and (iii) as a discriminative shape signature to find similar objects, classify objects, or combine shape information with other object-associated data. Here, we focus on the general principles behind the most commonly used shape descriptors.

\subsubsection{Global Shape Statistics}
These features capture properties of the whole shape. Common examples are area, convex area, sphericity, eccentricity, and extent. Global features provide insufficient local shape information for fine-grained shape distinction, e.g., a smooth vs a rough sphere. However, they are unmatched in interpretability, aligning with our human perception of shape, and in identifying global trends in shape variation, being robust to errors in upstream processes like image segmentation used to extract shape. As such, global shape statistics are often used to test for statistically significant changes in two shape collections, one from a treatment condition, the other from a control condition. Such comparative problems are central in scientific discovery and arise commonly in the fields of biology and medicine.

\subsubsection{Local Shape Statistics} 
These features characterize the shape characteristics in a local region. At the finest spatial resolution, this involves defining measurements per pixel or voxel for images or per vertex for meshes. These primarily measure the shape fluctuation within a local window. Notable examples include curvature-based measures, such as mean and Gaussian curvature, and roughness measures, such as the distance from a smoother reference shape. Features defined at the finest resolution are suitable for accurate shape matching and registration but are typically too granular for pattern-finding applications such as shape clustering, which also requires the same number of features to be extracted per shape. In these cases, the shape is partitioned into a fixed number of larger regions, and statistics of the fine-grained per pixel/voxel/vertex measurements, such as the minimum, maximum, mean, standard deviation, and kurtosis, are extracted per region and concatenated to form region-based shape feature descriptors. The regions are commonly generated either by partitioning the object surface/volume into approximately equal surface areas/volumes by hard-coded rules, e.g., subdivision by angle; by clustering, accounting for local distance and convexity; or by instance segmentation for shapes with well-defined object parts like protrusions~\cite{driscoll2019robust}. In the latter, the number of parts may vary between individual object shapes, necessitating further summary statistics or histogram-based 'bag-of-words' approaches to ensure a fixed number of features per shape~\cite{wang2014bag}.

\subsubsection{Distributional Statistics}
Considers a shape as a point cloud and measures properties of how points relate spatially to one another. Measurements can involve all pairwise combinations of points. The chordal length histogram measures all pairwise Euclidean point distances to construct a distribution of the lengths \cite{zhou2025development}. However, this is inefficient to compute for complex shapes with high surface curvature features or long, thin protrusions, which require a large number of points to represent \cite{zhou2025development}. Moreover, there can be high redundancy where many points have effectively the same local point distribution. A common approach is to isotropically resample the point cloud to a given number of points \cite{schlomer2011farthest,jacobson2012fast,zhou2025development}. However, if the number of points is too small, the result may be too coarse. Instead, the point cloud distribution can be characterized with respect to a smaller number of ``reference'' points. The typical approach is illustrated by the construction of the shape context descriptor obtained by flattening and concatenating the 2D histogram counting the number of points in each bin from partitioning the space around each reference point by angle and distance~\cite{belongie2000shape}. The distance is usually uniformly sampled in log-space to correct for the larger area covered by bins further away. The main drawback is selecting a fixed number of reference points that sample the whole shape and are located with neighborhoods unique enough to fingerprint the shape but general enough to also be found in other instances of the same shape type.  Certain objects, such as faces, have well-defined and conserved features. When this is not the case, landmark finding techniques (as reviewed in Section~\ref{sec:Landmarking}) are often used. For instance, using surface evolution techniques to identify those points that minimize the shape reconstruction error when the number is fixed~\cite{ramer1972iterative,douglas1973algorithms,garland1997surface,barkowsky2000schematizing}. Using landmark points is often more representative than uniform random sampling, which misses protrusions that often provide the most distinctive shape information. 

\subsubsection{Topological Characteristics}
 Topological data analysis (TDA) aims to construct feature descriptors with more mathematical formalism, based on measuring topological invariants of the shape~\cite{chazal2021introduction}. A notable example is the Euler characteristic for polyhedra, which relates to the genus, the number of ``holes'' in the object. As defined, it is not a useful discriminative descriptor because the Euler characteristic is the same for any shape with the same number of holes, and therefore, too permissive in practice. Instead, the Euler Characteristic Curve (ECC) is computed by scanning the object along a number of angles, partitioning the object into two and measuring the Euler characteristic of the residual object~\cite{smith2021euler,dlotko2023euler}. Shapes with the same Euler Characteristic will have different ECC because the residual object after splitting may have multiple spatial components. The conceptual and computational simplicity of the ECC has led to many variants. Notably, Turner et al.~\cite{turner2014persistent} introduced the Smooth Euler Characteristic Transform (SECT), generalizing the Euler Characteristic to represent each shape or image as function-valued data. Recently, Marsh et al. developed DEtecting Temporal shape changes with the Euler Characteristic Transform (DETECT), a rotationally invariant and dynamic extension of SECT for 2D and 3D shape classification applied to organoid growth~\cite{marsh2024detecting}. The ECC example illustrates the general steps in constructing topological feature descriptors by persistent homology. A filtering operation is used to scan and measure the object’s topological properties at multiple spatial scales. Importantly, the process is used to construct a 2D persistence diagram tracking the birth and death of topological features that also uniquely identifies a shape. The 2D persistence diagram is not directly amenable as a feature vector~\cite{tauzin2020giottotda}. Instead, various functional summaries are used, such as the persistence landscape, Betti curves, and persistence images. 
 Alternatively, topological features like SECT are already equipped with a well-defined inner product structure through reproducing kernel Hilbert spaces, directly allowing for machine learning applications, such as comparison of protein structures~\cite{Crawford2}, and hypothesis testing on populations of shapes from binary~\cite{meng2025randomness} and gray-scale images~\cite{Crawford4}. Key to the effectiveness of topological shape descriptors is the filtering function which ideally is chosen to highlight the most salient shape characteristics. Many filtering functions have been developed, for example, for images, thresholding at increasing grayscale intensities; for branching tree-like shapes, the distance from the root node; for graphs and meshes, the n-ring neighborhood; and for unstructured point clouds, increasing radial distance from individual points.

 \subsubsection{Decomposition Methods}
 These methods describe a shape as a function of a set of basis functions. Unlike the other approaches above, the original shape can still be reconstructed from the coefficients and basis functions. The coefficients are concatenated and used as the shape feature descriptor where they may be further processed, such as by normalization or taking the absolute value, to ensure they provide desirable properties like translation-, scale- and rotational-invariance~\cite{Burger2013}. The basis functions can be directly specified. If so, usually the bases are chosen to be well-characterized orthonormal functions where efficient algorithms have been developed, notably Fourier coefficients~\cite{Burger2013} and Zernike moments for 2D shapes~\cite{nixon2025feature}, and spherical harmonics for 3D shapes~\cite{ruan2019evaluation,Viana2023}. For arbitrary shapes, matrix factorization techniques are used to jointly infer the basis functions and coefficients. The resulting basis functions are specific to individual shapes or shape collections, instead of being universal. The matrix that is factorized is typically a covariance function between pairwise points or a differential operator, notably the Laplace--Beltrami. The choice of matrix factorization influences the properties of the resultant basis function. Principal components analysis (PCA) and eigendecomposition approaches produce orthonormal bases analogous to spectral decomposition~\cite{kosambi2017statistics}. Non-negative matrix factorization (NMF) imposes an additive constraint on the coefficients in the decomposition~\cite{lee2000algorithms}. This leads to sparser coefficients and generates bases that often correspond more to distinct subobjects in a shape.  Considering a human face, PCA finds spatial frequency modes, whereas NMF finds distinct facial features like lips, eyes, and nose~\cite{zhang2008topology}. The former is more efficient to compute and can reconstruct the shape exactly. Therefore, PCA is well-suited for enabling multiscale processing such as hierarchical shape registration. The latter NMF is more computationally expensive, lossy with a finite reconstruction error but more semantic and discriminative, better suited for applications such as identifying distinct landmark features, segmenting the shape into distinctive parts such as protrusions, and shape classification.  

\subsubsection{Data-Learned Descriptors}
\label{sec:features_ML}
These descriptors are generated by fitting neural network models to data using an objective function. The objective can be supervised, based on shape classification; unsupervised, based on data reconstruction using an autoencoder or self-supervised, based on constructing a proxy classification-like task that could be solved using only the data, without semantic annotations. Self-supervised tasks are usually either imputation problems, subdividing the data object into parts and training a model to impute the full object as output using only a subset of parts as input~\cite{he2022masked}; or contrastive learning problems, whereby the model predicts for an object which of two is more similar to it~\cite{costa2025transformed,shu20253d}. 

After model fitting, the feature descriptor is taken as the output of an intermediate network layer. Typically, this is the latent or bottleneck layer.  Data-learnt descriptors can yield several advantages. First, the features may be more specifically tailored to the dataset and the objective. Second, it allows for end-to-end learning whereby the feature-generation can be directly informed by the intended downstream task. Third, the feature computation after model fitting is fast and scalable, as neural networks naturally leverage batch processing and deep learning libraries seamlessly take advantage of parallelized, accelerated GPU compute without user knowledge of GPU coding. Lastly, neural networks allow more flexible integration of data types. Categorical variables allow models to condition on metadata like shape categories~\cite{makhzani2015adversarial}. Perhaps most impressively, neural networks realized the learning of realistic probabilistic models of shape, through sequential nonlinear transformation of simple random variables, in particular of Gaussian random variables, which is not possible by handcrafted design~\cite{vahdat2022lion,zhang20233dshape2vecset}. A notable application is medicine, where generative models are used to sample the larger potential shape space to predict the most likely structure and binding of proteins without existing crystal structures~\cite{abramson2024accurate}, and to design de novo chemical compounds~\cite{alakhdar2024diffusion}. However, fitting the data does not mean that the models have learned geometry or physics~\cite{huang2021arapreg,sassen2024repulsive,yang2025genanalysis}. This leads to unrealistic interpolation or deformation of shapes when interpolating the latent feature space and not mapping shapes in a temporal sequence as nearest neighbors in latent space. There is also no guarantee of performing well on alternative objectives for which the model was not trained on. Notably, using autoencoder-trained features for handwriting digit recognition is significantly worse than supervised classification, but supervised features are poor at interpolation because each digit forms an individual separated cluster in latent space. This is important as it may be difficult to formulate properties as a suitable training objective, and objectives may counter each other. Problems may be open-ended and not have a well-defined objective. Lastly, individual features are often not interpretable without further correlation to the handcrafted features discussed above~\cite{zaritsky2021interpretable,de2025geometric}.

\textbf{Geometric deep learning (GDL)} has emerged as a powerful framework for learning feature-based representations of non-Euclidean data structures such as point clouds, curves, graphs, and meshes. Unlike traditional deep learning models designed for grid-structured data (e.g., images or sequences), geometric methods operate on spaces with irregular geometry and complex topological structure. A unifying principle in this domain is the invariance of learned representations to shape-preserving transformations, such as reparameterizations and rigid motions. A canonical example is learning from point clouds, where the space $\mathcal{I}$ consists of sets of points in $\mathbb{R}^3$ and the symmetric group $\mathcal{D}$ acts by permutation. The PointNet model~\cite{qi2017pointnet} was among the first to design neural architectures invariant to such permutations, effectively modeling functions over unordered point sets. Subsequent extensions~\cite{qi2017pointnet, wang2019dynamic, thomas2019kpconv,FoldingNet} incorporate local neighborhood structures, hierarchical representations, and continuous convolutional filters, enhancing modeling capacity while preserving invariance at the shape level. Similar principles extend to curves and graphs, where objects are represented as parameterized curves in $\mathbb{R}^2$ or $\mathbb{R}^3$, or as graphs encoding adjacency or connectivity. In these cases, the reparameterization group $\mathcal{D}$ may include interval diffeomorphisms for curves or graph automorphisms. Graph-based neural architectures-such as spectral graph convolutions~\cite{bruna2013spectral,defferrard2017convolutionalneuralnetworksgraphs} and spatial graph neural networks like GCNs~\cite{kipf2017semisupervisedclassificationgraphconvolutional}-enable learning of invariant functions over such non-Euclidean domains. For 3D surface meshes, geometric deep learning methods extend graph convolutions to irregular collections of vertices, edges, and faces, while addressing invariance to reparameterization and mesh resampling. Early approaches applied point cloud or graph convolutions directly to mesh vertices and edges but often ignored higher-order connectivity and lacked invariance to surface parameterization. More advanced methods leverage spectral techniques using the mesh Laplacian~\cite{hanocka2019meshcnn,Kostrikov_2018_CVPR} or graph convolutions on primal and dual mesh graphs~\cite{NEURIPS2020_0a656cc1}. Recently, parameterization-invariant measure representations, such as varifolds, have been successfully applied to model complex shapes~\cite{hartman2023varigrad,hartman2024varishape}. GDL based representations of shapes data can be leveraged for various downstream learning tasks. These include serving as inputs to multi-layer perceptron (MLP) models for classification or regression, as well as to autoencoder frameworks for shape reconstruction. In the case of reconstruction, several strategies can be employed to achieve more faithful reconstructions, such as selecting appropriate loss functions and regularization terms~\cite{ARAPReg,LIMP}, or designing specialized decoder network architectures~\cite{groueix2018,hartman2024basis}. 
Feature-based representations are also valuable in dynamic settings, such as modeling and learning trajectories of geometric data for tasks such as surface interpolation and motion extrapolation~\cite{geolatent,hartman2024varishape}.

 \subsubsection{Distance-Based Descriptors} The pairwise distance between the individual shapes $N$ produces a $N\times N$ distance matrix where each row $i$ is the feature descriptor of shape $i$~\cite{johnson2015joint,govek2023cajal,zhao2025fast}. The feature descriptors may be further processed to constrain the numerical values, e.g., by L1, L2 normalization or a kernel function applied to convert the distance matrix to an affinity matrix to improve suitability for downstream data analysis and machine learning applications. The choice of distance directly influences the representation power of the feature descriptor (see Section~\ref{sec:distance}). Another factor is the number of points used to represent the shape, which determines the spatial resolution being compared. Too coarse may miss desired geometric features, and too fine may slow computation and include non-discriminatory or unwanted noisy geometric features~\cite{govek2023cajal,zhou2025development,de2025geometric}. Generally, pairwise distance descriptors are a type of meta-feature descriptor (see subsection below) that may yield excellent performance in clustering and classifying shapes, as they already encode relational information~\cite{govek2023cajal}. However, key drawbacks are the $N\times N$ pairwise matrix computation, which is slow to compute and memory-consuming, making the approach harder to scale to larger datasets. Moreover, the pairwise nature requires recomputing features over the whole dataset whenever new shapes are added. Consequently, much of the research effort is focused on developing and using more efficient low-rank approximations of distance metrics, notably the Nystr\"om method~\cite{nystrom1930praktische,marsh2024detecting,zhou2025development}.

    \subsubsection{Meta-feature descriptors} Meta-feature descriptors are constructed by processing shape features computed for individual shapes using the other described methods. The processing involves the entire shape collection, and therefore meta-feature descriptors cannot be computed for a single object. The rationale for their use is that the individual feature descriptor may be too high-dimensional and noisy. Importantly, its values do not have a point of reference. Meta-feature descriptors take into account feature similarity relative to other shapes in the dataset. There are two main approaches for their construction. The first is capturing a shape’s similarity to all others by computing the feature pairwise distance matrix of objects. An example of such a descriptor is the L2-normalized row vectors of the distance matrix. To capture nonlinearity, a kernel function, such as the RBF kernel, can be further applied to the distance matrix~\cite{marsh2024detecting}. The second approach is to find a small number of prototypical shapes, typically by clustering on the individual shape descriptors and then constructing a meta-feature descriptor reflecting conceptually the description of each shape with respect to the prototypical shapes, such as a vector of feature distance to each prototype. This is also often known as prototype or archetypal analysis~\cite{epifanio2018archetypal}. When prototype analysis is performed on parts of shapes, a histogram descriptor can be constructed for the full shape by allocating each of its parts to the most similar prototype part and counting the number of instances of each prototype part~\cite{wang2014bag}. This is also called a ``bag of features'' or ``bag of words'' approach~\cite{o2011introduction}. Alternatively, Fisher vector uses Gaussian mixture models to capture second-order information beyond first-order counts, usually resulting in a descriptor that is a more compact and dense representation for classification and retrieval~\cite{sanchez2013image}.

\subsection{Distance in the Shape Space}
\label{sec:distance}

The end goal of shape space analysis is to compare shapes and uncover hidden connections between them. A fundamental step in this process is the construction of distances that faithfully capture similarities and differences between geometric objects.

In practice, the design of a shape distance involves balancing several desirable properties. Most fundamentally, a distance should be invariant to irrelevant transformations such as reparameterization, rigid motion, and scaling, so that it reflects intrinsic geometry rather than representation. At the same time, it should be discriminative, assigning small values only to genuinely similar shapes. Additional considerations include stability to noise and discretization, computational efficiency and scalability to large datasets, and the ability to integrate naturally with statistical and learning-based analyses. No single construction optimizes all of these criteria simultaneously, which has led to a variety of complementary approaches in the literature.

The choice of distance directly influences downstream tasks such as clustering, classification, and statistical analysis of shapes. Broadly, we divide the distances into several categories, each offering distinct advantages and challenges depending on the data and application:

\paragraph*{Shape Descriptor-Based Metrics.} We first note that it is possible to define a distance metric between feature descriptors that summarizes the shapes and are equipped with a distance metric, and many options are outlined in Subsection~\ref{sec:feature_descriptors}. This approach is often adopted if specific shape descriptors are of interest based on domain knowledge. 

\begin{figure}[t]
    \centering    
    \includegraphics[width=\linewidth]{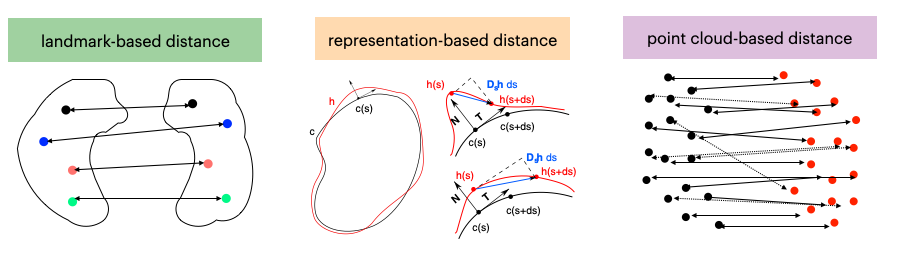}
    \caption{\textbf{Illustration of distances between shapes depending on the data type.} The distances can typically be categorized as (left) landmark-based distance, (middle) representation-based distance, and (right) point cloud-based distance.}
    \label{fig:shape_distance}
\end{figure}

For a general distance, we might consider the input being a raw representation of shapes without any summary, as a set of landmarks or a full point cloud, curves or meshes, areas or volumes. As illustrated in Fig.~\ref{fig:shape_distance}, there are three types of general metrics.

\paragraph*{Landmark Correspondence-Based Metrics.} These distances rely on landmarks, a predefined set of corresponding points across shapes. Metrics such as Procrustes distance measure differences after aligning landmarks such that transformations like rotation, translation, and scaling are not quantified. While these methods are intuitive and computationally efficient, they are sensitive to landmark placement and may struggle with missing or noisy correspondences. In addition, manually placing landmarks requires expert knowledge, is time-consuming, and may also introduce bias. 

In contrast, landmark-free distances, such as Continuous Procrustes Distance, Wasserstein distance, currents-based metrics, and diffeomorphic shape matching, treat shapes as continuous objects and measure differences based on their intrinsic geometric or topological properties. These methods are more flexible, particularly for shapes lacking well-defined landmarks, but often involve higher computational complexity. The choice between these approaches depends on the specific nature of the shapes under study and the desired trade-off between accuracy and efficiency. In the following, we discuss three such categories:

\paragraph*{Representation-Based Metrics.} In many circumstances, the shapes are represented as 2D shape contours/occupancy images, or 3D mesh/volume, and for each class, specific distances have been developed and used, such as the Riemannian metric for curves/surfaces etc. Some of them are also correspondence-based, often through a minimization procedure. We note that this category differs from the landmark-based distance above because it utilizes all the information on the representation instead of just the set of landmarks. For instance, the curve-based metric requires knowledge about the order of points, and the mesh-based metric can take curvature into account.

\paragraph*{Point-Cloud Based Metrics.} When shapes are represented as point clouds, i.e., without information about connectivity as in curves or surfaces-distance metrics are often defined directly on the point sets. Examples include optimal transport-based distances, which do not require explicit correspondence or topological structure and instead rely on geometric proximity across point clouds.

\paragraph*{Data Manifold Metrics.} Another category of distances relies on the data manifold when many shapes are available, and the underlying manifold structure generating data can be inferred. First, the entire dataset is modeled through a diffusion process where the transition probability is often determined based on one of the previously mentioned distances. The diffusion distance is then determined by the aggregation of all possible paths between two shapes.

\subsubsection{Shape Descriptor-Based Distances}

Section~\ref{sec:feature_descriptors} discussed methods for characterizing shapes using feature descriptors rather than their full geometric representations. Examples include global shape statistics (e.g., area, volume), local shape statistics (e.g., surface protrusions), topological characteristics, and data-driven or learned descriptors.
Here, in Table \ref{tab:descriptor-metric}, we summarize several classical distance metrics defined between descriptor representations of two shapes, along with the corresponding input feature types.

\begin{table}[t!]
    \centering
    \begin{tabular}{|p{2.8cm}|p{2.5cm}|p{4cm}|p{3cm}|p{3cm}|}
        \hline
        \rowcolor{gray!30}
        \textbf{Metric} & \textbf{Input} & \textbf{Definition/Concept} & \textbf{Key Features} & \textbf{Challenges} \\\hline
        \hline
        $L^p$ distance and its variants (e.g., Mahalanobis) & Global shape statistics (e.g., area, eccentricity) & \( \|x - y\|_p \) on vectors in \( \mathbb{R}^d \) & Simple, interpretable, fast to compute; Mahalanobis account for correlation & Misses spatial and structural features \\\hline
        
        $p$-Wasserstein distance \cite{Bifurcation} & Distributional features (e.g., histograms) & Optimal transport cost between distributions; equivalent to $L_p$ norm of quantiles & Captures shape variability; invariant to reordering & Computationally expensive in dimensions higher than 1 \\\hline
        Bottleneck \cite{CohenSteiner2006} / 1-Wasserstein \cite{Mileyko2011} & Persistence diagrams & Bottleneck: \( \inf_\gamma \sup_{x} \|x - \gamma(x)\|_\infty \);   1-Wasserstein: $
W_1(D_1, D_2) = $ $\displaystyle \inf_{\gamma: D_1 \to D_2} \sum_{x \in D_1} \| x - \gamma(x) \|_2$

 & Stable to noise; captures topological events & Hard to vectorize; computationally costly \\\hline

        Fourier distance \cite{Zahn1972} & Fourier coefficients of 2D boundary (closed curve) & Euclidean distance on low-frequency Fourier coefficients & Compact; can achieve rotation/scale invariance & Requires consistent parametrization \\\hline
        Zernike distance \cite{HaeKwangKim} & Zernike moments of 2D images in disk & Euclidean distance on Zernike moments & Orthogonal; localized; rotation invariant & Requires disk normalization; limited to 2D \\\hline
        Spherical harmonics (SH) power distance \cite{kazhdan2003rotation} & Functions on a sphere (often in 3D) & \( \sum_\ell \|a_\ell^{(1)} - a_\ell^{(2)}\|^2 \), where \( a_\ell \) are SH coefficients & Encodes global 3D surface; invariant to rotation & Needs spherical parameterization; genus-0 only \\\hline
    \end{tabular}
    \caption{A summary of shape descriptor-based distances.}
    \label{tab:descriptor-metric}
\end{table}

\subsubsection{Landmark-Based Distances} 

\begin{table}[t!]
\begin{center}
\begin{tabular}{|p{2.8cm}|p{1.6cm}|p{4cm}|p{2.5cm}|p{4.3cm}|}
\hline
\rowcolor{gray!30}
\textbf{Metric {\color{red}}} & \textbf{Input} & \textbf{Definition/Concept} & \textbf{Key Features} & \textbf{Challenges} \\\hline
\hline
Procrustes Distance (full/partial) & Landmarks & Measures similarity between shapes by aligning them using translation, scaling, and rotation. & Alignment-based; fast computation & Limited to rigid transformations; requires manually annotated landmarks; not suitable for nonlinear shape changes. \\
\hline
Alpha Procrustes distances \cite{minh2022alpha}& Landmarks & Generalize a family of Procrustes distances in a parameterized manner.& Mean kernel metrics & Allows for additional flexibility as compared to Procrustes Distances. \\
\hline
\end{tabular}
\end{center}
\caption{A summary of landmark-based metrics.}
\label{tab:ldmk-metric}
\end{table}

Next, we consider landmark-based metrics (Table~\ref{tab:ldmk-metric}). First, we describe the distances defined when landmarks are available. Consider two shapes in $d$ dimensions, each represented by $L$ landmark points, chosen to have a one-to-one correspondence.

One well-known approach to model shapes mathematically is due to Kendall~\cite{kendall2009shape}. Formally, the shape of a landmark configuration is all geometrical information about it that remains when location, (isotropic) scale and rotational effects are filtered out. Distance in Kendall's shape space is thus related to the Procrustes analysis that superimposes objects such that their position, orientation, and (optionally) scale match. Depending on the presence of scaling, one distinguishes between full and partial Procrustes superimposition.

\begin{definition}
    Given two shapes $S_1$ and $S_2$ represented by two sequences of homologous (corresponding) landmarks $\{x_\ell ^{(j)} \}_{\ell = 1}^L$ for $j=1,2$. The {\bf partial Procrustes distance} is defined as:
\begin{equation}
    d_P(S_1,S_2) = \inf_{T \in \mathbb{E}(d)} \left( \frac{1}{L} \sum_{\ell=1}^L \| T \left( x_\ell^{(1)} \right) - x_\ell^{(2)} \|^2 \right)^{\frac{1}{2}},
\end{equation}
    where $\mathbb{E}(d)$ is the group of all rigid motion in $d$ dimensions.
\end{definition}

This distance is a metric on size-and-shape space as it is still sensitive to the (absolute as well as relative) scale of the objects. To factor out translation and scale effects, we can further center landmark configurations and resize them to unit size (such that the total variance of landmark coordinates is one). The resulting set of configurations is called the pre-shape space and forms a hyper-sphere. Orbits of pre-shapes under the action of rotations partition the pre-shape space into fibers that correspond one-to-one with shapes. Distances on Kendall's shape space can thus be obtained from the closest distances between fibers in terms of distances on the (pre-shape) sphere.
Restricting the partial Procrustes distance to pre-shapes corresponds to a chordal distance. A geometrically natural distance can be inherited from the projection of the fibers to shapes (i.e.,\ the Riemannian submersion) and corresponds to the closest great circle distance. This Riemannian distance can be expressed in terms of the partial Procustes distance as $\rho(S_1, S_2) = 2\arcsin(d_p(S_1, S_2)/2)$. As remarked by~\cite{al2013continuous}, in some applications, it may be advantageous to consider \textbf{weighted Procrustes distances}, in which each landmark $\ell = 1,...,L$ can be assigned a weight in the computation of the distances. For example, if there are imbalances in the distribution of the landmark points when some regions have higher densities than others.

The distances presented so far are inherited from Euclidean metrics of the spaces $\mathbb{R}^d$ containing the landmarks and do not account for any dependence between landmark pairs. However, in many disciplines, landmarks cannot be assumed to be independent and identically distributed. In biology, for example, we can expect landmarks that are close together (such as the centers of the eyes) to show a more coherent movement than distant ones due to diverse biological processes of regulation. Shape metrics that weight changes in small distances more heavily than changes in larger ones, thus, may be more appropriate for those applications. One approach for such metrics due to Bookstein~\cite{bookstein1989principalwarps} is to derive a metric that gauges shape changes in terms of the induced bending energy modeled via thin plate splines. More recently~\cite{kilian2007geometric,heeren2014exploring}, metrics have been derived based on the objects' surfaces rather than a finite set of landmarks. This line of work defines Riemannian metrics that gauge the severity of shape changes by quantifying the strength of surface deformation. Such approaches are commonly referred to as shell space due to the connection of the metrics to (linearized) bending and stretching energies from continuum mechanics of thin elastic shells. The Riemannian metrics can further be augmented to avoid self-collisions and, hence, nonphysical inter-penetration~\cite{sassen2024repulsive}. 

In general, shell spaces do not admit explicit expressions for geodesics and, hence, require computationally demanding numerical schemes to estimate shortest paths. One strategy to alleviate this challenge is to perform a change of variables by encoding shapes with a collection of features such as triangular surface descriptors~\cite{heeren2016splines}, skeletal representations~\cite{pizer2022skeletons}, and differential coordinates~\cite{tycowicz2018dcm, ambellan2021fcm}.
Geodesic calculus in the according feature spaces can be performed in closed form and, in many cases, accounts for local geometric variability, thus faithfully capturing the nonlinearity inherent to large deformations.

\subsubsection{Representation-Based Distances}

\begin{table}[t!]
\begin{center}
\begin{tabular}{|p{2.2cm}|p{1.5cm}|p{5cm}|p{4cm}|p{2.5cm}|}
\hline
\rowcolor{gray!30}
\textbf{Metric {\color{red}}} & \textbf{Input} & \textbf{Definition/Concept} & \textbf{Key Features} & \textbf{Challenges} \\\hline
\hline
Continuous Procrustes Distance \cite{al2013continuous} & Curves, \mbox{Surfaces} & Extends Procrustes distance to account for continuous deformation between shapes. & Captures both global and local shape changes. & Computationally more expensive; depends on quality of parameterization. \\
\hline
Symmetric distortion energy \cite{Koehl2015} & Curves, Surfaces & Find globally optimal conformal maps between shapes & Metric property for genus zero surfaces, robust to noise & Computationally intensive \\\hline
Elastic Metric \cite{ElasticMetric} & Curves, \mbox{Surfaces} &Measures shape similarity by re-parameterizing curves or surfaces to align optimally under elastic deformation. & Shape deformation captured explicitly; allows for smooth comparison. & Requires optimization over a complex parameter space. \\
\hline
Shell space metric \cite{heeren2014exploring} & Surfaces & Physically-based distance reflecting viscous dissipation required to deform a thin shell. & Riemannian metric (modulo rigid motions) on space of immersions. Captures nonlinear shape variability. & Requires correspondence; computationally expensive. \\\hline
Differential coordinates \cite{tycowicz2018dcm} & Surfaces, \mbox{Volumes} & Employs a differential representation that puts the local geometric variability into focus. & Fast, explicit computation of Riemannian expressions, Captures nonlinear shape variability. & Requires correspondence;\\
\hline
Fundamental coordinates \cite{ambellan2021fcm} & Surfaces & Surface-theoretic approach encoding surfaces via their fundamental forms. & Rigid motion invariant; Fast, explicit computation of Riemannian expressions, captures nonlinear shape variability. & Requires correspondence; Nonlinear out-bedding problem.\\\hline
Varifold Distance \cite{Varifold} & Curves, \mbox{Surfaces} & Measures similarity based on varifold representations of curves and surfaces. & Bi-invariant to reparameterization without optimization, relatively computationally efficient & The distance does not capture large shape deformation.\\
\hline
\end{tabular}
\end{center}
\caption{A summary of representation-based metrics for curves, surfaces, and volumes.}
\label{tab:ldmk-free-metric}
\end{table}

Table~\ref{tab:ldmk-free-metric} provides an overview of representation-based metrics for curves, surfaces, and volumes.

In the absence of landmark points, a reasonable approach is to pick a large number of ``reasonably well distributed'' points on each shape object and perform alignment between different surfaces as an approximation to possible correspondence maps. This motivates the following distance~\cite{al2013continuous}:

\begin{definition}
    Given two surfaces $S_1$ and $S_2$, denote by $\mathcal{A}(S_1,S_2)$ the set of all area-preserving diffeomorphisms, the \textbf{continuous Procrustes distance} between two surfaces is given by:
\begin{equation}
    d_{cP} = \inf_{ \mathcal{C} \in \mathcal{A}(S_1,S_2) } \left[ \inf_{T\in \mathbb{E}(d)} \left( \int_{S_1} \| T(x) - \mathcal{C}(x) \|^2 d \ vol_{S_1}(x) \right)^{1/2} \right].
\end{equation}
\end{definition}

Within the framework of conformal (angle-preserving) maps and focusing on finding globally optimal conformal mapping, the symmetric deformation energy~\cite{Koehl2015} measures the distortion energy:
\begin{definition}
Given two surfaces $S_1$ and $S_2$, the \textbf{symmetric distortion energy} is defined by:
\begin{equation}
d_{sd} = \min_f \sqrt{\int_{S_1} (1-\lambda_f(z))^2\ dA_1} + \sqrt{\int_{S_2} (1-\lambda_{f^{-1}}(z))^2\ dA_2},
\end{equation}
where $\lambda_f(z)$ is a dilation factor that measures the stretching of vectors by $f$ at each point $z$ in $S_1$, and the map $f$ is a conformal map with $\lambda_f(z)=1$ for every $z$. 
\end{definition}

Along the same line of geodesics for the Kendall space, one can view the space of parameterized curves and/or surfaces as an infinite-dimensional manifold and equip it with a Riemannian metric that is invariant under the action of the reparameterization group. If the objects are given by smooth parameterized curves $[0,1]\to \mathbb{R}^2$ with nowhere-vanishing derivative, a square root velocity metric~\cite{ElasticMetric} can be used to quantify the morphological difference between two objects:

\begin{definition}
\label{def:elastic_metric}
    For two parameterized curves $c_1$ and $c_2$ and a regular parameterized path $\alpha\colon [0,1]\to \mathcal{C}$ such that $\alpha(0)=c_1$ and $\alpha(1)=c_2$, the length of $\alpha$, the geodesic distance (\textbf{general elastic metric}) between $c_1$ and $c_2$ is:
\begin{equation}
    d_{a,b}(c_1,c_2) = \inf_{\alpha:[0,1]\to \mathcal{C}, \alpha(0)=c_1,\alpha(1)=c_2} \int_0^1 g_{\alpha(t)}^{a,b} ( \alpha'(t), \alpha'(t) ) ^{1/2} dt,
\end{equation}
    where $g^{a,b}$ denotes an inner product between two curve deformations $h$ and $k$ in the tangent space:
\begin{equation}
    g_c^{a,b}(h,k) = a^2 \int_0^1 \langle D_s h, N \rangle \langle D_sk, N\rangle ds + b^2 \int_0^1 \langle D_sh,T\rangle\langle D_sk, T\rangle ds,
\end{equation}
with $D_s=\frac{1}{\| c'(s) \|}\frac{d}{ds}$ denote a differential operator with respect to the arc length $s$, and $N$ and $T$ respectively are the local unit normal and tangent vectors. In particular, for a certain choice of the parameters of this metric ($a=1/2,b=1$), this family of metrics includes the well-studied \textbf{square root velocity metric ($d_{SRV}$)}. The square root velocity transform~\cite{ElasticMetric} defines the isometry 
\begin{equation}
    \Phi: (\mathcal{I}, d_{SRV}) \to (L^2(M,\R^n),\|\cdot\|_{L^2})
\end{equation}
via
\begin{equation}
  \alpha \mapsto \frac{\alpha'}{\sqrt{\|\alpha'\|}}.
\end{equation}
Furthermore, in~\cite{Bauer2024Elastic}, this simplifying transform was generalized for $d_{a,b}$.
\end{definition}

The Riemannian metric can be generalized to surfaces~\cite{hartman2023elastic} as follows:

\begin{definition}
    For two parameterized surfaces $q_1$ and $q_2$ and a regular parameterized path $q$ with $q(0)=q_1$ and $q(1)=q_2$, the geodesic distance between $q_1$ and $q_2$ is:
\begin{equation}
d_{G}(q_1,q_2) = \inf_{q:[0,1]\to \mathcal{I}} \int_0^1 G_{q(t)} ( \partial_t q(t), \partial_t q(t) )^{1/2} dt,
\end{equation}
where $G$ is a Riemannian metric defined on the space of surfaces. Options for strong Riemannian metrics that induce meaningful distances (first- and second-order Sobolev metrics) can be found in~\cite{hartman2023elastic}.
\end{definition}

Another density-based representation of 3D shape data comes from the study of \textbf{varifolds}~\cite{Stouffer2024}. In this context, data is represented by a density on $\R^3$ times the Grassmannian $G_d(\R^3)$ corresponding to the location of points and the tangent space to the shape respectively. This space of measures, denoted $\mathcal{M}(\R^3,G_d(\R^3))$, is equipped with a Reproducing Kernel Hilbert Space (RKHS) norm.  
\begin{definition}
    Given two varifolds, $\mu,\nu\in\mathcal{M}(\R^3,G_d(\R^3))$ the inner product between $\mu$ and $\nu$ is defined as:
    \begin{equation}
        \langle\mu,\nu\rangle_{V} = \int_{\R^3\times G_d(\R^3)}\int_{\R^3\times G_d(\R^3)} \rho(x,y)\gamma(u,v) d\mu(x,u) d\nu(y,v)
    \end{equation}
    where $\rho$ is a positive kernel on $\R^3$ and $\gamma$ is a positive kernel on $G_d(\R^3)$. The \textbf{varifold norm} is then defined by
    \begin{equation}
        \|\mu-\nu\|^2_V = \langle\mu,\mu\rangle_{V} + \langle\nu,\nu\rangle_{V} -2\langle\mu,\nu\rangle_{V}.
    \end{equation}
\end{definition}
When studying 3D surfaces (or graphs), the varifold representations may be discretized on the faces (edges, respectively), allowing for the varifold norm to be approximated by a simple double sum over the faces (edges) of the two objects. The resulting metric offers a parameterization blind metric between 3D surfaces and graphs; however, the application of this metric requires a careful choice of kernels $\rho,\gamma$ for the specific application to obtain meaningful results.

A category of distances focuses on comparing 2D areas or 3D volumes through registration, that is, finding a mapping that transforms one object into another. In such scenarios, images or volumes are often given as binary or grayscale intensity maps over voxelized domains. Registration between such data is often achieved via minimizing a potential characterizing the discrepancy between the transformed image and its target. Here we present a brief survey on such methods along with their objective functions.

\begin{definition}
Large deformation diffeomorphic metric mapping (LDDMM)~\cite{beg2005computing} is defined as the solution to the following variational problem that involves two images $I_0$ and $I_1$:
    \begin{equation}
\min_{v:\frac{d\phi_t}{dt}=v_t(\phi_t)} \ \left( \int_0^1 \|v_t\|_V^2 dt + \| I_0 \circ \phi_1^{-1} - I_1 \|_{\ell^2}^2\right),
    \end{equation}
where $\|v_t\|_V$ is an appropriate Sobolev norm on the velocity field that controls the smoothness of the registration mapping, and the second term is the squared-error norm that enforces matching of the images. The latter term can be adapted to incorporate any distance metrics between images.
\end{definition}

Under image settings, the choices of distance metrics have been broadened to better reflect underlying geometry, especially under multi-modal settings where images are taken from different imaging devices and/or modalities. For instance, people have used mutual information~\cite{MI} as a new matching criterion. The optimization and interpretation were improved in other distances like normalized gradient fields (NGF)~\cite{Haber2006}, which penalize differences in the orientation of local intensity transitions, thereby enforcing alignment of structural edges instead of raw brightness levels:
\begin{definition}
The normalized gradient field (NGF) distance between two images 
$I_0$ and $I_1$ is defined as
\begin{equation}
D_{\mathrm{NGF}}(I_0,I_1)
=
\int_\Omega 
\| n(I_0,x) \times n(I_1,x) \|^2 \, dx,
\end{equation}
where
\begin{equation}
n(I,x)
=
\frac{\nabla I(x)}
{\sqrt{\|\nabla I(x)\|^2 + \varepsilon^2}},
\end{equation}
and $\varepsilon > 0$ is a small regularization parameter.
\end{definition}

\subsubsection{Point Cloud-Based Distances}

In this section, we introduce metrics that accommodate point clouds as input. Standard geometric baselines for comparing point sets include the \textbf{Hausdorff distance}~\cite{Birsan}, which measures the maximum distance from a point in one set to the nearest point in the other (worst-case error), and the \textbf{Chamfer distance}~\cite{Athitsos}, which computes the average squared distance to the nearest neighbor. While computationally efficient and widely used in deep learning, these metrics rely purely on nearest-neighbor matching and often fail to capture the underlying density or global structural properties of the shape.

To address these limitations, recent approaches employ optimal transport (OT), which treats shapes as probability distributions. Two key metrics derived from OT are the Wasserstein distance and the Gromov--Wasserstein distance, both of which are particularly well-suited for analyzing structured data.

\begin{definition}
    Given two probability distributions \(\mu\) and \(\nu\) defined on a metric space \((\mathcal{X}, d)\), the \textbf{p-Wasserstein distance} is defined as:
\begin{equation}
W_p(\mu, \nu) = \left( \inf_{\gamma \in \Pi(\mu, \nu)} \int_{\mathcal{X} \times \mathcal{X}} d(x, y)^p \, d\gamma(x, y) \right)^{1/p},
\end{equation}
where \(1 \leq p < \infty\), and \(\Pi(\mu, \nu)\) is the set of all joint distributions (couplings) \(\gamma\) on \(\mathcal{X} \times \mathcal{X}\) with marginals \(\mu\) and \(\nu\). Intuitively, \(W_p\) measures the minimal transportation cost to morph one distribution into another under the given ground metric \(d\)~\cite{Villani2009}.
\end{definition}

While the Wasserstein distance compares distributions over a common metric space, many practical scenarios require comparing distributions supported on different metric spaces, such as shapes or graphs with distinct geometries. The \textit{Gromov--Wasserstein distance} generalizes the OT framework to this setting. 

\begin{definition}
Given two metric spaces \((\mathcal{X}, d_\mathcal{X})\) and \((\mathcal{Y}, d_\mathcal{Y})\), and probability measures \(\mu\) on \(\mathcal{X}\) and \(\nu\) on \(\mathcal{Y}\), the \textbf{$p$-Gromov--Wasserstein distance} is defined as:
\begin{equation}
\text{GW}_p(\mu, \nu) = \left(\inf_{\gamma \in \Pi(\mu, \nu)} \int_{\mathcal{X} \times \mathcal{Y} \times \mathcal{X} \times \mathcal{Y}} \left| d_\mathcal{X}(x, x') - d_\mathcal{Y}(y, y') \right|^p \, d\gamma(x, y) \, d\gamma(x', y') \right)^{1/p}.
\end{equation}

\end{definition}

The GW distance measures how well the metric structures of \((\mathcal{X}, d_\mathcal{X})\) and \((\mathcal{Y}, d_\mathcal{Y})\) can be aligned under a probabilistic coupling. It has proved especially useful for tasks involving shape matching, graph alignment, and cross-domain similarity assessments~\cite{Memoli2011}.

Note that more variants of Wasserstein and Gromov--Wasserstein have been introduced to generalize their applications to arbitrary positive measures (e.g., unbalanced or partial Wasserstein or Gromov-Wasserstein metrics), making them suitable for comparing shapes of varying masses~\cite{tajmir2025alignment,koehl2023physicist}. Computing these metrics in practice also presents significant computational challenges. For discrete probability measures in dimensions greater than one, calculating the Wasserstein distance reduces to solving a linear assignment problem, which can be addressed using the Hungarian algorithm with cubic time complexity. The computational burden becomes even more severe for the Gromov--Wasserstein distance, which is widely believed to be NP-hard \cite{GW-NPHard}.
To address these computational limitations, entropic regularization can be employed to relax the original formulations~\cite{cuturi2013sinkhorn,zhang2024gromov}), as well as slicing and quantization techniques \cite{SlicingW, SlicedGW, QuantizedGW,zhao2025fast}, enabling practical application to shape point clouds (e.g., 3D volumes of proteins~\cite{riahi2023alignot,ecoffet2021application}). 

Finally, a hybrid approach that integrates alignment directly into the transport framework is the \textbf{Procrustes--Wasserstein distance}~\cite{Wass-Procustes}:
\begin{definition}
    Given two probability distributions $\mu$ and $\nu$ in Euclidean space, the Procrustes--Wasserstein distance is:
\begin{equation}
    \mathcal{PW}_2(\mu,\nu) = \inf_{T\in \text{ISO}(d)} W_2(T\# \mu, \nu),
\end{equation}
    where $T$ is restricted to the group of rigid isometries (rotations and translations). This metric jointly optimizes for the best spatial alignment and the optimal transport plan.
\end{definition}

\subsubsection{Diffusion-Based Distances}

The previous metrics are defined over the ambient space of data objects. In this section, we discuss a family of distances defined over the intrinsic geometry of the data manifold. The fundamental idea~\cite{Coifman2005} is to build a Markov chain on the data and use the resulting transition probabilities to capture the intrinsic geometry.

\begin{definition}
    Given a collection of objects endowed with a distance metric. Let $(W_t)_{t>0}$ be the diffusion matrices constructed by means of these distances, and let $\psi_1,...,\psi_L$ be the first few dominant non-trivial eigenvectors of $W_\tau$, for an appropriately chosen $\tau>0$ with corresponding eigenvalues $\lambda_1,...,\lambda_L$. For each $\gamma>0$, the diffusion distance is defined as:
\begin{equation}
    D^2_\gamma(k,k') = \sum_{\ell=1}^L \lambda_\ell ^{2\gamma} | \psi_\ell (k) - \psi_\ell (k') |^2.
\end{equation}
\end{definition}

If landmarks and approximative correspondence maps are available as well, the diffusion distance can be adapted to account for discrepancies in correspondences between landmarks via horizontal diffusion distances that model the diffusion process on fiber bundles. We refer the readers to~\cite{gao2021diffusion} for details. 

\subsection{Shape Dynamics}

\label{sec:dynamics}
Analyzing individual object dynamics is an area of intense research with diverse applications to both industry and science. In particular, the evolution history of shape provides indispensable causal information for clustering objects by similar dynamics and for predicting future object behavior. This has proved invaluable in medical imaging for disease monitoring~\cite{namburete2023normative,ahmad2023multifaceted} and in cell biology for functional drug screening~\cite{freckmann2022traject3d,zhou2025identifying}, and in characterization of biological heterogeneity~\cite{zhou2025development,zhou2025identifying}. Due to the differences in how time is acquired and how shapes are derived, the exact shape analysis problem and its solution are often specific to the dataset and subject domain. Nevertheless, there are four common recurring problems.

\begin{enumerate}
    \item \textbf{Discovery of temporal orderings within an arbitrary collection of shapes.}\\
    When it is not possible to track the same object over time or the data represents a single snapshot in time, a core task is reconstructing time based on ordering objects by similarity. This can be particularly difficult if the time structure is unknown a priori, as is common in paleontology, protein structure, and in-situ cell imaging.
    \item \textbf{Alignment and staging of the same common dynamic process between different individual objects.}\\
    The dynamic process of interest is known, for example, the particular disease or biological phenomenon. However, the timing, specifically the onset and endpoint of the process, can vary between individuals. Consequently, staging and alignment are necessary to conduct like-for-like study of the dynamic process.
    \item \textbf{Temporal feature extraction to cluster objects by similar shape and evolution trajectory or for classification.}\\
    The primary goal for medicine and scientific discovery is to generate a categorical classification label to inform treatment in the former and to produce a definitive conclusion in the latter. This necessarily involves transforming the raw shape time series into a feature vector suitable for classification machine learning algorithms.
    \item \textbf{Correlation and causal inference between shape dynamics and other measurements.}\\
    For many scientific disciplines, particularly medicine and biology, shape represents an inexpensive measurement that does not require a specific experimental assay or labeling technique to extract. This makes it useful as a proxy for other signals-of-interest, which is particularly attractive in imaging studies, where the number of microscopy channels is limited by spectral wavelength separability. The ability to correlate molecular signals of interest to specific shape features frees up channels for other imaging. Shape is also a preserved ``signal'' across different spatial scales and is robust to appearance variations. This makes it suitable for integrating measurements, e.g., coupling subcellular signals to tissue features through cell shape or coupling measurements between different imaging assays. Shape is also integral to biological function. In cell biology, a natural killer cell, a type of immune cell, undergoes a sequence of shape changes to engage cancer cells~\cite{sapoznik2020versatile,zhou2023surface}, whereas a T cell, another key type of immune cell, shifts morphology depending on whether it is migrating or ``scanning'' its surroundings~\cite{jenkins2023antigen}. The development of high-resolution microscopy has led to increasing interest in discovering and quantifying these phenomena in an unbiased manner. To this end, non-perturbation causal inference approaches are being developed to find the true causal sequence of changes from background shape variations~\cite{noh2022granger,zhou2023surface}. 
\end{enumerate}

To address these questions, three main analytical approaches are commonly used, balancing implementation simplicity, mathematical rigor, and accuracy against dataset constraints:
\begin{enumerate}[nosep]
    \item \textbf{Manifold learning applied to feature descriptors.} A straightforward, versatile approach suitable for large, shape-heterogeneous datasets. The method can be applied almost universally to shapes from both snapshots and temporal data. Time points can be treated as individual shapes with their own feature descriptors, or a single feature descriptor can be constructed that captures the entire shape sequence. Moreover, non-shape information can be integrated by concatenating shape features with non-shape features into a single feature descriptor. However, the results of the analysis depend on the comprehensiveness and discriminatory power of the features~\cite{zhou2025development}.
    \item \textbf{Constructing a spatiotemporally-invariant coordinate system.} This approach establishes dense pointwise correspondences between all shapes in the dataset, allowing for precise geometrical comparison. Necessarily, it requires accurate registration, which can be challenging or impossible to achieve, particularly if additional properties such as diffeomorphism are desired to be preserved, and shapes have different topological characteristics. As such, this approach is typically used to construct atlases mapping shape variations of the same object. For temporal data, shape dynamics should be sampled at sufficient temporal resolution so that changes are small, so as to be smooth and diffeomorphic between consecutive time points.   
    \item \textbf{Trajectory modeling.} This approach seeks to represent shape dynamics as trajectories in a shape space. Key questions involve finding or fitting a trajectory to discretely observed shapes and comparing the shape trajectories. 
\end{enumerate}

\subsubsection{Manifold Learning Applied to Feature Descriptors}

\label{sec:manifold_analysis_temporal_features}
A common approach to analyzing time-varying shapes is to represent each shape by a feature descriptor (see Section~\ref{sec:feature_descriptors}) and embed these descriptors into a low-dimensional latent space (typically 2D or 3D) using manifold learning. This representation enables easy visualization and comparison of shape variability across all time points and objects~\cite{freckmann2022traject3d,zhou2025identifying,stower2023single}. When explicit time indices are available, it is also common to construct meta-feature descriptors to characterize the entire temporal shape sequence of the same object for clustering and classification of shape dynamics. Typically, direct concatenation of per-timepoint descriptors leads to too high a dimensionality to be discriminative. Additionally, each sequence comprises a different number of time points. Therefore, temporally-aware representations are constructed using summary statistics (e.g., feature means)~\cite{wiggins2023cellphe,wiggins2025cellphepy}; manifold learning such as PCA to create compressed features~\cite{copperman2023morphodynamical}, tSNE and UMAP to embed the pairwise similarity between sequences~\cite{wang2020live,alieva2024behav3d}; or deep learning methods~\cite{soelistyo2022learning,shannon2024cellplato,wiesner2024generative,vasan2025interpretable}. Alternatively, transition dynamics are learnt using state-space models, with model parameters serving as descriptors~\cite{zhong2012unsupervised,gordonov2016time,freckmann2022traject3d}.  

The structure of the latent space depends on the features and the manifold learning method. PCA preserves global metric relationships, providing interpretable axes, and scales efficiently to large datasets, but is limited to linear structure and captures predominantly size and eccentricity in the first two components regardless of dataset, thereby providing less discrimination~\cite{ruan2020image,segal2022vivo,zhou2025development}. Nonlinear methods are therefore more commonly used. Spectral embeddings like Diffusion maps, Laplacian eigenmaps, and potential of heat diffusion for affinity-based transition embedding (PHATE) preserve neighborhood relationships based on diffusion distance and tend to produce branch-like embeddings, making them suitable for trajectory-like organization~\cite{rajpoot2008unsupervised,moon2019visualizing,haghverdi2016diffusion,faigenbaum2026studying}. t-distributed stochastic neighbor embedding (t-SNE) emphasizes local similarity, forming local clusters~\cite{zhou2025development}. Uniform manifold approximation and projection (UMAP) preserves both local and global structure under a Riemannian manifold assumption~\cite{de2025geometric,zhou2025development} whilst PaCMAP,  TriMap, and LocalMAP balance attractive and repulsive forces to embed shapes with similar features together and dissimilar features further away~\cite{wang2021understanding}. 

A limitation of standard manifold learning methods is the assumption that feature vectors are independent and identically distributed. Consequently, embeddings of the same object at different time points do not have to preserve temporal proximity~\cite{zhou2025development,zhou2025identifying}. Temporal coherence can be incorporated through neural network architectures that introduce constraints into the embedding dynamics, such as Koopman-based models that learn linear operators~\cite{brunton2021modern} and partial differential equations in latent space~\cite{rao2023encoding,regazzoni2024learning} or neural architectures (e.g., LSTMs, causal convolutions) that model temporal dependencies~\cite{fortuin2018som,soelistyo2022learning,huijben2023som, wang2024latent}. Alternatively, contrastive learning can enforce partial temporal ordering by embedding neighboring time points closer than distant ones~\cite{yang2019learning,wang2021self}. An advantage of PCA and neural network–based approaches is that new samples can be embedded without recomputing the entire model, unlike t-SNE and UMAP, which embed pairwise similarity relations, facilitating scalable analysis of evolving shape data.

A key application of constructing a lower-dimensional shape space is to simplify the solution to shape analysis questions. Notably, one can infer continuous trajectories between shapes without explicit temporal indices by principal curve-fitting~\cite{gilles2025cryo,yu2025mapping,Magana2026} and fitting of diffusion processes using the latent space coordinates~\cite{haghverdi2016diffusion,lange2022cellrank,weiler2024cellrank,raj2025cell}. These can then be used to define a pseudotime relative to a designated origin point, giving a measure of developmental time for cell differentiation in genetics~\cite{ding2022temporal,saelens2019comparison} and reconstructing a temporal history of disease development from a patient tissue biopsy sampled at a single time point~\cite{van2021revealing}. The space represents a semantic landscape and geometry used to construct a cell population representative fate trajectory robust to individual cell heterogeneity~\cite{zhou2025development,zhou2025identifying}. This was used to visualize the fate of distinct cell populations in a mouse embryo during development~\cite{stower2023single}. In the same study, this space enabled single-cell motion features to be integrated with embryo geometrical features like curvature and molecular intensity features across multiple embryos with different numbers of time points, enabling the definition of complete spatiotemporal morphogenetic cell profiles for distinct cell types~\cite{stower2023single}. 

\subsubsection{Temporally-Invariant Coordinate Systems}
Manifold learning of feature descriptors is powerful for analyzing heterogeneous collections of shapes to discover global patterns. However, features destroy precise location information and therefore is not well suited to single-object morphometric studies that require geometry preservation—such as measuring local shape variation relative to a reference~\cite{li2024unveiling}, removing confounding geometric variation for comparative analysis~\cite{stower2023single}, or comparing surface- or volume-associated signals across time (e.g., molecular intensity~\cite{noh2022granger,zhou2023surface}, subcellular structures~\cite{zhou2023surface,lefebvre2025nellie}, developmental patterns in embryos~\cite{wong20154d,stower2023single,lefebvre2025dynamicatlas}, or disease-related changes in brain anatomy~\cite{zhang2018longitudinal,ahmad2023multifaceted,namburete2023normative}). In all cases, spatially-aware and topologically-correct registration is necessary to establish explicit pointwise correspondence of shapes between time points. For objects with clear anatomical landmarks, such as teeth or faces, sparse landmark matching can estimate geometric transformations, and this may be sufficient to infer shape trajectories, measure shape similarity, and cluster differences~\cite{robinson2001planar,al20143d,gao2019gaussian, driscoll2019robust}. Often, however, particularly in medicine and biology, it is desired to establish dense, point-wise correspondence across the entire surface and volume to factor out shape-induced differences in systematic comparison of signals—such as contrast, molecular expression, and internal structure morphometrics. The dense problem is significantly more challenging, and the optimal registration method is often dataset-dependent. Generally, methods require incorporating smoothness and continuity constraints~\cite{vercauteren2009diffeomorphic,deng2022survey,vishnevskiy2016isotropic}, metric preservation~\cite{hartman2023elastic,sassen2024repulsive}, and deformation priors~\cite{christensen1996deformable, beg2005computing,vercauteren2009diffeomorphic,mosleh2025data,lefebvre2025dynamicatlas} to interpolate transformations across unmatched or ambiguous regions. Volumetric image registration-based methods are popularly used, particularly diffeomorphic demons~\cite{vercauteren2009diffeomorphic}, large deformation diffeomorphic metric mapping (LDDMM)~\cite{beg2005computing}, symmetric diffeomorphic image registration~\cite{avants2008symmetric}, and parametric displacement fields~\cite{vishnevskiy2016isotropic} (see also Section~\ref{sec:alignAndRegister}). Then, to construct a shared coordinate system for a shape collection, either (i) all shapes are pairwise registered to a single template - typically a mean shape (Section~\ref{sec:mean}) or a single shape~\cite{stower2023single,mcdole2018toto}; or (ii) a groupwise approach which aims to jointly minimize the shape variation across all shapes in the collection by dynamically generating the registration template, either by a tree-based shape similarity search~\cite{faigenbaum2026studying} or optimization of the joint loss~\cite{vishnevskiy2016isotropic}. Topological differences between shapes~\cite{younes2010shapes,franccois2022weighted,zhou2023surface}, for example, arising due to segmentation errors or incomplete sampling and high curvature features such as cell surface protrusions~\cite{driscoll2019robust,mazloom2023cellular}, particularly challenge dense registration. For the former, metamorphic registration, mixing diffeomorphic transport with additional intensity terms~\cite{younes2010shapes,franccois2022weighted}; and for the latter, shape parameterization (Section~\ref{sec:surfaceUnwrapping},~\ref{sec:volume_parameterization}) to a canonical geometry such as the disk, plane, or sphere~\cite{postelnicu2008combined,zhou2023surface,yu2025mapping} are good starting points. In particular, shape parameterization of 3D shapes enables mapping to the 2D domain. This mapping naturally serves as a surface-matching constraint. Then, high-curvature features can be readily tracked using well-established 2D algorithms~\cite{zhou2023surface} or registered by treating them as landmarks~\cite{postelnicu2008combined,guo2023automatic}.

\subsubsection{Shape Trajectories}
\label{sec:shape_trajectories}

In many phenomena, such as growth or non-rigid motion, shapes can be assumed to vary smoothly and hence to evolve along a curve in shape space, referred to as a \textit{shape trajectory}.
Usually, we only have access to discrete trajectories given empirically by (sometimes a very limited number of) sampled shapes. Estimating and analyzing shape trajectories often provides a more detailed assessment than investigating individual samples and pairwise differences alone.
In particular, shape trajectories are themselves geometric objects that feature location, size, orientation, and shape. However, standard multivariate statistics only test for differences in location and hence are limited in assessing dynamic variation due to differences in the latter attributes~\cite{collyer_phenotypic_2013}.

As shape trajectories take values in nonlinear shape spaces, processing and analysis tools from multivariate statistics and signal processing are not readily applicable. A common workaround is to linearize the shape data by mapping it to the tangent space at a reference shape~\cite{Adams2013}. This allows the employment of multivariate regression procedures, for example, to investigate size-driven shape change (i.e., allometry) in growth and development, adaptation, and evolution~\cite{adams_ontogenetic_2010, drake_pace_2008, klingenberg_size_2016, klingenberg_methods_2022}. Significance of the slope is tested non-parametrically~\cite{adams_ontogenetic_2010, cardini_larger_2013, drake_pace_2008, klingenberg_morphological_2008, klingenberg_size_2016, klingenberg_methods_2022}. Beyond linear models, principal component analysis found widespread application. In allometry, the dominant principal component (PC1) tends to incline toward the size dimension in form space when there is a significant size-shape relationship, as variability is usually much larger in size than in shape. 
In phenotypic trajectory analysis~\cite{collyer_phenotypic_2013}, the direction of a trajectory is defined as PC1 (fitted to points within that trajectory) and serves as a means to estimate angular differences between subgroups.

\begin{figure}[t]
    \centering
    \includegraphics[width=1\linewidth]{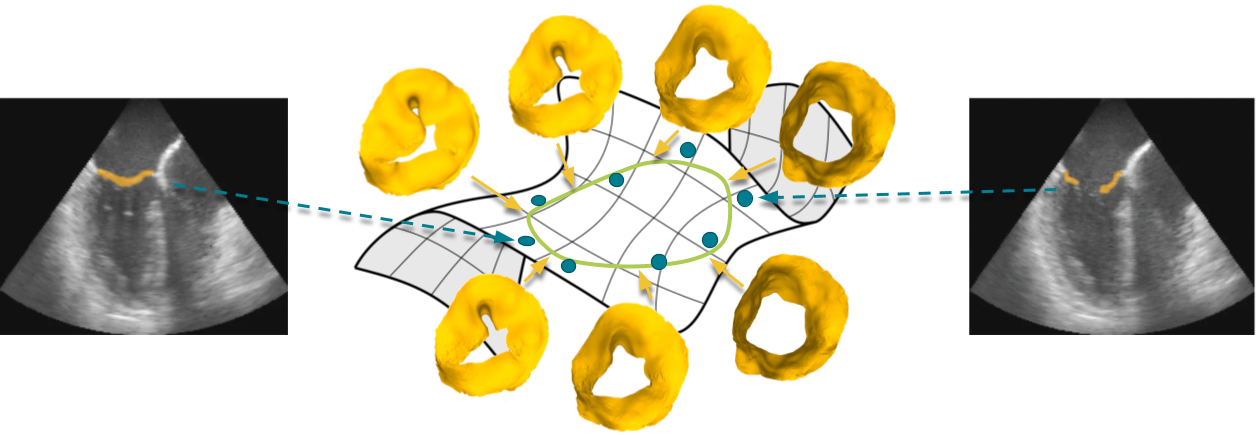}
    \caption{\textbf{Reconstructed meshes from regression of longitudinal mitral valve data covering a full cardiac cycle.} The spline consists of two cubic segments fitted to surfaces extracted from ultrasound images.}
    \label{fig:spline}
\end{figure}

Working exclusively in the tangent space introduces distortions, the severity of which increases with the dispersion of the sample set and in regions of high curvature. This behavior spurred the development of methods that take the nonlinear geometric structure of shape space into account. A fundamental processing task addressed in this field is estimating a smooth trajectory from time-indexed samples that captures important characteristics while leaving out noise or other high-frequency phenomena. To this end, different curve-fitting schemes have been adapted to shape space, such as smoothing splines~\cite{su2012fitting, kim2021smoothing}, geometric flows~\cite{Brandt2016}, and non-parametric regression, including kernel~\cite{davis2010population} and Gaussian process~\cite{mallasto2018wrappedGP} regression. In many cases, we can assume a parametric family of spatiotemporal models to capture the underlying dynamics adequately. Utilization of a parametric model allows us to describe shapes at unobserved times (i.e.,\ shape changes between observation times and -- within certain limits -- also at future times). A common approach is to approximate the time-varying shape data by geodesics in shape space~\cite{Fletcher2013,nava2020geodesic}. Geodesic models are attractive as they feature a compact representation, similar to the slope and intercept term in linear regression, and therefore allow for computationally efficient inference. Non-monotonous shape changes, such as those present in time series of cardiac shape motion (see Fig.~\ref{fig:spline}) or anatomical changes in the human brain over the course of decades, generally do not adhere to constraints of geodesicity. For such phenomena, generalized polynomial~\cite{Hinkle2014} as well as spline regression~\cite{singh2014splines} have been derived based on variational characterizations~\cite{camarinha2023high}. More recently, constructive approaches~\cite{HanikNavayazdanivonTycowicz2024} for intrinsic spline models have been adapted to higher-order regression~\cite{Hanik_ea2020} as they provide exact evaluation without the need for time-discrete approximation schemes.
Beyond geometric characterizations, specialized models encoding transformations inherent to biological or physical shape dynamics have been proposed, e.g., to capture growth~\cite{hsieh2022mechanistic,Charon2023} or secondary motion effects in character animation~\cite{schulz2015animating}, respectively.

Another analysis task is to study distributions of shape trajectories.
When analyzing such observations of shape trajectories, we have to distinguish between morphological differences due to temporal shape evolutions in a single subject, on the one hand, and the geometric variability in a population, on the other. Regression and other techniques intended for cross-sectional data analysis do not take into account the intrinsic correlation between repeated measurements of the same subject, nor do they inform how an instance relates to a population-average trend.
Methods for longitudinal data analysis must therefore be able to capture and disentangle the cross-sectional variability in shape and the temporal variability due to the underlying processes of change.
To this end, hierarchical models have been generalized to the case of shape trajectories.
In their simplest form, hierarchical models describe the data in two stages: at the individual level, subject-specific changes are modeled via smooth, parametric curves. At the group level, the subject-specific trajectories are in turn considered as perturbations from a population-average trend.
In this field, geodesic hierarchical models~\cite{MuFl2012,singh2013hierarchical,nava2020geodesic} that encode both subject-specific and group trends in terms of geodesics have been introduced based on least-squares criteria.
Alternatively, probabilistic formulations~\cite{kim2017riemannian,bone2020learning} were proposed for which subject-specific trends are assumed to be parallel to the average geodesic.
The case of multiple, possibly categorical covariates was addressed in~\cite{hong2019hierarchical} by employing a multi-geodesic approach.
As for regression, more flexible models beyond geodesics were also investigated.
The non-parametric approach in \cite{campbell2018nonparametric} estimates the group average trend via kernel regression on subject-specific geodesics.
In~\cite{HanikHegevonTycowicz2022} and~\cite{NavaYazdaniAmbellanHaniketal2023}, hierarchical models based on intrinsic B\'ezier splines are considered utilizing Riemannian metrics on the B\'ezierfold~\cite{HanikNavayazdanivonTycowicz2024}, i.e.,\ the manifold of B\'ezier splines. Working with parametric models implies that we consider the trajectories to be in a consistent parameterization so that points along the trajectories are in correspondence. However, if the trajectories are temporally misaligned (e.g.,\ due to errors in measurements), we may want to ignore this phase variability entirely and only compare the shapes of the trajectories. To this end, generalizations of the elastic metric (see Definition~\ref{def:elastic_metric}) that apply to manifold-valued curves have been proposed~\cite{su2014mfdSRVF,su2018comparing,celledoni2018homoSRVF}.

\subsection{Mean Shape}
\label{sec:mean}

For shapes and data that do not naturally live in a Euclidean space, it is essential to introduce suitable notions for mean, fluctuations, etc., as the usual definitions in Euclidean space, which heavily rely on affine operations, do not apply.

Fr\'echet~\cite{frechet} generalized the affine notion of the mean of data in a Euclidean space to metric spaces by suitably extending the minimum-square-error property of the ordinary Euclidean mean. This so-called {\em Fr\'echet mean}, also referred to as barycenter in more geometric settings, gives rise to the {\em Riemann center of mass}~\cite{grovekarcher, karcher77, karcher14} on a smooth Riemannian manifold. The notion of mean was further extended to quasi-metric spaces in~\cite{ziezold}. A {\em generalized Fr\'echet mean} defined for continuous functions on topological spaces was introduced in~\cite{huckemanna,huckemannb} by minimizing a squared loss function that links the data space to the descriptor space. In practice, efficient algorithms and software libraries have been developed for computing Fr\'echet means on manifolds and metric spaces. For instane, the open-source \texttt{Geomstats} package \cite{Geomstats} implements a variety of mean-estimation methods across Riemannian, Finsler, and elastic shape spaces, including the Fr\'echet mean under the elastic metric.

When working with the Fr\'echet mean, it is crucial to be aware that properties of the mean and results like the central limit theorem which can feel second nature need no longer apply. To appreciate that the Fr\'echet mean need not be unique one could consider, e.g., two antipodal data points on a sphere whose Fr\'echet mean is the great circle formed by all those points which are equidistance to both data points. Results on the uniqueness of the Fr\'echet mean are discussed in~\cite{karcher77,kendall90,groisser05,afsari11}. On metric spaces which exhibit non-unique Fr\'echet means, there are inevitable statistical challenges in the estimation of Fr\'{e}chet mean, as discussed in~\cite{estimatefrechet}.

Bhattacharya and Patrangenaru in~\cite{CLT-BP} have established a central limit theorem for Fr\'echet means on Riemannian manifolds, which has been followed by many generalizations such as~\cite{bhattacharyabhattacharya, huckemanna, ellingson13, battacharyapatrangenaru14, patrangenaruellingson16,bhattacharyalin17}. A key observation which was made in~\cite{hotzhuckemann15} and was then studied further in~\cite{EltznerHuckemann19,Eltzner22} is that asymptotic fluctuations for the sample mean can occur at scales with exponents strictly less than $\frac{1}{2}$ and can converge to non-normal distributions. This phenomenon, which has been termed {\em smeariness}, can only occur in the non-Euclidean setting. In contrast to smeariness, where a slower convergence to the mean is observed, the Fr\'echet mean can also exhibit {\em stickiness}, where small perturbations of a sample do not change the sample mean yielding faster convergence to the mean, see, e.g.,~\cite{phylo3, phylo1, stickyopenbook,stickyCLT}. A common example which illustrates the phenomenon of stickiness is obtained by taking the 3-spider and placing unit point masses on each leg of the spider the same distance away from the center point. In this case, the sample mean, which lies at the center point, does not change under perturbations of the positioning of the unit point masses. More generally, on stratified spaces the Fr\'echet mean tends to stick to lower-dimensional strata where three or more strata are joined. A central limit theorem for Fr\'{e}chet means in the space of phylogenetic trees has been developed in~\cite{phylo2} and for Fr\'echet means on stratified spaces in the recent series of works~\cite{stratCLT1,stratCLT2,stratCLT3,stratCLT4}. These non-standard asymptotic rates that Fr\'echet means can exhibit in non-Euclidean spaces are inherently related to (synthetic) curvature, with the rough guide that positive curvature leads to smeariness while negative curvature leads to stickiness. It is important though that even within the standard central limit regime curvature effects can alter finite sampling rates up to considerable sample sizes~\cite{curvatureeffects}, with implications for statistical inference~\cite{Ulmer2023}.

Alternatives to the Fr\'{e}chet mean which alleviate issues of smeariness and stickiness are the {\em diffusion means} recently introduced in~\cite{diffusionmean1,diffusionmean2}. Except for some corner cases, diffusion means satisfy a central limit theorem with the standard rate. For $p$ the minimal heat kernel on a complete and connected Riemannian manifold $M$, the {\em diffusion $t$-mean} $E_t(X)$ of a random variable $X$ on $M$ is the set of global minima of a log-likelihood function such that
\begin{equation}
    E_t(X)=\operatorname{argmin}_{y\in M}\mathbb{E}\left[-\ln p(t,X,y)\right].
\end{equation}
As discussed in~\cite{diffusionmean1,diffusionmean2}, this definition can be extended to spaces with non-smooth structures. Notably, diffusion means are defined intrinsically and, in the Euclidean setting, they reduce to the usual mean. Moreover, due to the Varadhan asymptotics $\lim_{t\to0}2t\ln p(t,x,y)=-\mathrm{dist}(x,y)^2$, the diffusion means can be considered as generalizations of the Fr\'echet mean. They can be estimated through bridge sampling,
a field that has been extensively treated in the literature in its own right and with recent extensions to non-Euclidean spaces such as in~\cite{jensen2019simulation,JS21,jensen2022discrete,buiInferencePartiallyObserved2023,corstanjeSimulatingConditionedDiffusions2024,sRbridge24}.
It should be stressed that~\cite{diffusionmean1,diffusionmean2} make the assumption that the underlying space is both geodesically and stochastically complete. While geodesic completeness has been studied and understood for a variety of distances and shape spaces, see, e.g.,~\cite{bauer2014overview,Bruveris17,maor23,maor24,discretesobolev,HPS25}, the study of their stochastic completeness property still requires attention going forward.
Recently, a broad characterization for spaces of discrete regular curves and for landmark spaces has been obtained in~\cite{curve_sto_complete} and~\cite{HHS24,HS26}, respectively.

Recently, an efficient algorithm for estimating mean and infinitesimal covariance for data on a Riemannian manifold has been developed in~\cite{mostprobable1} through the means of most probable paths for anisotropic Brownian motions. Most probable paths are paths on the Riemannian manifold that maximize a certain path probability, though with the path probability being measured with respect to the driving Euclidean Brownian motion that is mapped by stochastic development to the Riemannian manifold. The algorithm developed in~\cite{mostprobable1} obtains an approximate solution very efficiently owing to the use of constrained optimization. Most probable paths for stochastic flows with Eulerian noise and deterministic drifts, which, e.g., appear in fluid dynamics and shape analysis, as well as for developed processes are discussed in~\cite{mostprobable2,mostprobable3}. These works have applications in shape analysis, e.g., by providing ways of optimally transforming landmark configurations.

\section{Machine Learning to Uncover Morphological Patterns in Shape Space}

\label{sec:Statistical Analysis}

Up to this point, we explored various shape representations and parameterizations, alongside methods for improving correspondence and measuring similarity between shapes. These steps form the foundation of shape space analysis, offering domain experts quantitative tools to represent and investigate the data. Selecting an optimal shape representation, often closely tied to both the dataset and the research question, constitutes a critical step in any shape analysis workflow. As we will see in this section, such representations enable ML and statistical methods to more effectively uncover patterns in the data.

Although shape representations and other pre-processing steps stand on their own, they also provide the basis for applying machine learning and statistical methods to uncover underlying structures in the data. Tasks such as classification and modeling shape variation are essential for gaining deeper insight into species differences and other intrinsic characteristics of the dataset. In this section, we will review key methods in area. While Machine Learning (ML) methods are developed beyond shape analysis, the goal of this section is to highlight methods that have been developed and adapted for biological shape analysis. Comprehensive treatments for diving deeper into the subjects are available elsewhere (e.g.,~\cite{calder2025neural}).

Before delving into the technical details of specific algorithms, it is important to emphasize careful experimental design, including selecting appropriate methods, testing/training data, and interpreting the results. A recent paper highlighted a series of common missteps, errors, and violations of proper ML protocols that appear with disconcerting frequency in the growing body of literature~\cite{calder2022use}.

\subsection{How Far Off Are We? Accuracy Analysis}

In the following sections, we delve into the details of ML and statistical analysis. A crucial point to emphasize is that every new algorithm, method, or application to new data must be carefully evaluated for accuracy and robustness. It is important to estimate the extent to which a chosen method is effective.
In some cases, accuracy may be inherently limited. For example, insufficient data limits the ability to recognize patterns, which in turn compromises the identification of key properties. It is therefore advisable to test methods under multiple settings to obtain a clearer understanding of method performance. The remaining challenge lies in designing and carrying out this evaluation rigorously.

When evaluating error, multiple aspects that influence reliability must be carefully considered: First, the characteristics of the data itself are critical. Evaluation must be performed on real labeled data or on simulated data designed to approximate real scenarios. Sample size, distribution, and potential class imbalance also strongly affect the reliability of accuracy estimates. When pooling data from multiple sources to increase sample size, additional sources of variability must also be considered. Differences in imaging and digitization modalities (e.g., CT, laser scanning, photogrammetry), as well as interobserver variation, may introduce systematic error, particularly when comparing species or subspecies~\cite{zhang2023open, robinson2017error}.

Second, the validation strategy must be carefully chosen. This includes decisions such as whether to use a train/test split, cross-validation, leave-one-out, or bootstrap, and how to ensure that all clusters in the data are represented in a balanced manner. It is equally important to prevent leakage from the test data into the model, which could artificially inflate accuracy.

Third, error evaluation should be aligned with the task. Classification tasks are often assessed with confusion matrices, accuracy, precision, recall, or F1-score, whereas regression tasks may rely on metrics such as mean squared error (MSE), root mean squared error (RMSE), or mean absolute error (MAE). In shape and morphology studies, domain-specific metrics, such as geometric distances or similarity measures, are often most appropriate (see Section~\ref{sec:distance}).

Fourth, statistical significance and uncertainty should be assessed. Confidence intervals for error rates provide a measure of uncertainty, and hypothesis testing can help determine whether observed differences in performance are meaningful or merely due to chance. 

Finally, the cost of different types of errors must be evaluated. Depending on the research question, false positives and false negatives may carry unequal importance, and some mistakes may carry higher consequences than others.

Taken together, these considerations provide a structured and rigorous framework for evaluating algorithm performance. Properly addressing them ensures that reported results are both reliable and interpretable, reducing misinterpretation and overconfidence. For further discussion of common pitfalls in ML evaluation, see~\cite{liao2021we, calder2022use}.

\subsection{Shape Classification and Clustering}
The application of ML in shape space analysis focuses on extracting patterns and clustering based on overall similarity. These approaches address questions related to taxonomy, biological property identification, species recognition, evolutionary and ecological relationships, and disease phenotype classification. 

What distinguishes shape space analysis from classical ML applications is the flexibility in data representation. Classification can be performed directly on the shapes, on intrinsic features derived from them, or on feature representations embedded into low-dimensional spaces (see Section~\ref{sec:features_ML}). The choice of embedding level depends on several factors: the representation of shape (Section~\ref{sec:shape_space_overview}), the preprocessing applied, and the distance metric employed (see Section~\ref{sec:distance}).

In shape space analysis, the classification methods themselves do not differ from standard clustering or classification approaches commonly applied elsewhere~\cite{calder2025neural}. What differs is the nature of the data and research questions, which are inherently tied to shape representation and geometry. In addition, the assumption that shapes are sampled from a manifold can be leveraged to obtain improved solutions. Once shapes are represented in either high- or low-dimensional form, classification can be effectively carried out using well-established ML algorithms, such as random forests, support vector machines, K-means, or even phylogenetic tree–based approaches (e.g.,~\cite{yezzi2022using}). Cluster separation may subsequently be assessed using statistical procedures to compare group structure, such as one-way Procrustes ANOVA~\cite{zelditch_geometric_2012, adams_geomorph_2021, klingenberg_geometric_1998, baab2012shape}. Numerous applications illustrate this approach, including classifying and predicting obstructive sleep apnea from facial semilandmarks, classifying fossil human skull and wing shapes, and identifying bone fragments from archaeological records~\cite{yezzi2022using, macleod2017use, bellin2021geometric, monna2022machine, baab2021assessing}. 

\begin{figure}[t]
  \centering
  \includegraphics[width=.33\textwidth]{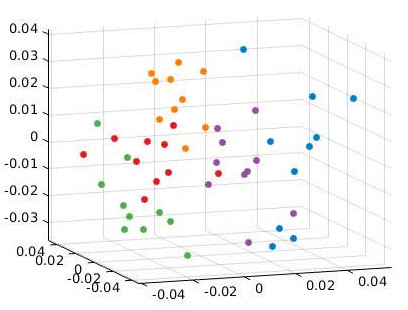}
  \includegraphics[width=.35\textwidth]{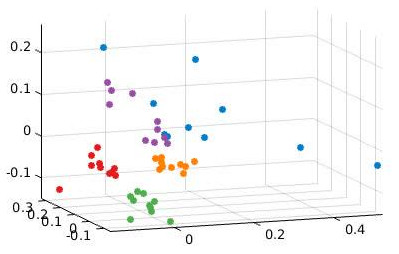}
  \includegraphics[width=.31\textwidth]{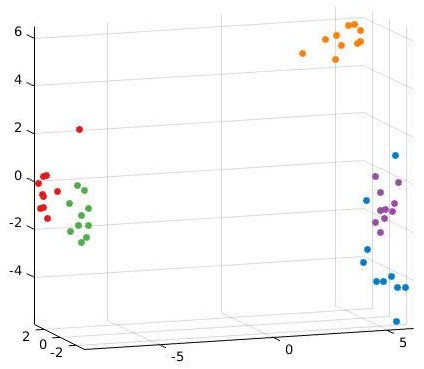}
  \caption{\textbf{Impact of data embedding on clustering performance.} Geometric Similarity of Anatomical Surfaces using different distance calculations.   Left: Embeddings into $\mathbb{R}^3$ using \emph{Multi-dimensional Scaling} of continuous Procrustes distance (Table~\ref{tab:ldmk-free-metric}); Middle: DM calculated after the diffusion map embedding; Right: and HDM distance~\cite{gao2021diffusion}.}
  \label{fig:MDS_comparison}
\end{figure}

What about the levels of data representation? As discussed in Sections~\ref{sec:shape_space_overview} and \ref{ref:shape_param}, an essential step is deciding how to represent the data itself. Once a representation is selected, along with an appropriate metric between representations, the analysis is carried out at that level. The chosen representation can strongly influence the outcomes, from obscuring meaningful structures to successfully revealing biologically meaningful groups. A clear example of the importance of choosing the right representation for clustering is presented in the summary paper~\cite{faigenbaum2026studying}. For example, at the shape level, clustering based on continuous Procrustes distance failed to yield meaningful group separation (Fig.~\ref{fig:MDS_comparison}, left). Diffusion distances provided modest clustering improvement, though still only weakly indicative (Fig.~\ref{fig:MDS_comparison}, middle). Once geometric information of the shapes was incorporated into the data embedding, the clusters immediately became well defined (see Fig.~\ref{fig:MDS_comparison}).

Although hierarchical clustering is a general clustering method, we emphasize its particular value in evolutionary studies. In evolutionary and shape analyses, a central question is how phylogeny is constructed and how it influences morphological variation~\cite{adams_phylogenetic_2019}. Once data are parameterized, phylogenetic trees can be inferred, often by analyzing covariance patterns among taxa means, which serve as the operational units. Phylogenies are typically derived from distance or covariance matrices that summarize similarities among taxa. 

To date, most landmark-based approaches focus on detecting and isolating the phylogenetic component of overall morphological variation represented by landmarks, given an existing Brownian motion-based phylogenetic tree for taxa within a sample \cite{adams_phylogenetic_2019}. These methods are largely adapted from univariate analyses. One commonly used method is the phylogenetic transformation approach based on the eigendecomposition of the phylogenetic covariance matrix for estimating variance associated with phylogeny in linear modeling~\cite{adams_phylogenetic_2019, adams_phylogenetic_2018}. Extensions of this approach can be used for evolutionary rate estimation and Phylogenetic Partial Least Squares (PLS) to assess trends of integration in phylogeny directly from landmark data~\cite{adams_phylogenetic_2019, adams_quantifying_2014}. Other methods include using a multivariate version of K statistics, $K_{mult}$, to estimate the strength of signals versus expectation~\cite{adams_phylogenetic_2019} and phylogeny-based PCA for visualizing overall variations adjusted for or associated with phylogeny~\cite{polly_phylogenetic_2013, adams_phylogenetic_2019}.

On the other hand, directly using landmarks or other shape features to reconstruct a phylogenetic tree remains uncommon due to the difficulty of defining discrete and localized morphological characters from global shape representations (see discussions in~\cite{adams_morphometrics_2011}).

\subsection{Hypothesis Testing within Shape Space}
At times, classification or machine learning algorithms alone cannot fully answer the questions posed by the data; in such cases, robust statistical analysis becomes essential for obtaining reliable and statistically significant results. Hypothesis testing in shape space provides a framework for evaluating whether observed morphological differences exceed what would be expected by chance, given the geometry of the data. 

To do this, the research question is expressed as a working hypothesis, and a corresponding null hypothesis ($H_0$), the opposite statement, is constructed. An experiment is then designed to test the data against $H_0$ and to compute the probability of observing the results at least as extreme as those obtained under the assumption of $H_0$. If this probability (the p-value) is sufficiently small, $H_0$ is rejected to support the working hypothesis; otherwise, the evidence is considered insufficient to reject $H_0$. In short, hypothesis testing is a systematic way of using data to determine whether patterns reflect real effects or are simply due to chance.

The first challenge is to translate the research question into a precise statistical working hypothesis. This often involves converting the question posed by domain experts (e.g., biologists) into a mathematical-statistical form. Once the optimal representation of the data is selected (Section~\ref{sec:shape_space_overview}), the research question-or a portion of it-can be formulated as a working hypothesis. For example, in the shape space of letters discussed in~\cite{faigenbaum2016algorithmic}, the research question was to determine how many different writers were represented in the dataset. This was approached by testing pairs of inscriptions: the working hypothesis stated that ``the two distributions of features belong to different classes,'' while the null hypothesis $H_0$ assumed that “both distributions belong to the same class.”

Another key challenge in shape analysis is the high dimensionality of the data, often coupled with limited sample sizes. In many studies of shape analysis, the number of variables far exceeds the number of observations, leading to the ``curse of dimensionality''. Addressing this requires careful dimensionality reduction or regularization to ensure robust statistical inference (see discussions in~\cite{adams_morphometrics_2011, rohlf_applications_1998, zelditch_discovery_2000, cryoemchallenge}). In addition, randomization-based procedures, such as permutation tests or bootstrap methods, are frequently used to generate empirical null distributions from the observed sample~\cite{collyer_method_2015, collyer_rrpp_2018, zelditch_geometric_2012, zhao2025fast}. These approaches allow statistical inference even when classical parametric assumptions are difficult to justify.

Recent studies have pushed statistical analysis of shapes beyond traditional Euclidean methods based on tangent projections, enabling rigorous inference in nonlinear shape spaces using geodesic distances. Meng et al.\ introduced a framework for random shape inference via the smooth Euler characteristic transform (SECT), mapping shapes into a reproducing kernel Hilbert space without relying on landmarks or diffeomorphisms, and applied functional ANOVA with a chi-square statistic from a Karhunen–Loéve expansion for two-sample tests~\cite{meng2025randomness}. Zhang et al.\ developed a non-parametric k-sample energy test on nonlinear shape spaces, including Kendall and elastic shape spaces, demonstrating its use for mitochondrial shapes without tangent projections~\cite{zhang2022nonparametric}. Guo et al.\ represented brain arterial networks as elastic graphs in a Riemannian quotient space, enabling computation of geodesics, means, and covariances, and illustrated PCA-based visualization of shape changes with age~\cite{guo2022statistical}. Robinson and Turner proposed a nonparametric framework using persistence diagrams from topological data analysis to test population differences, applied to brain fMRI data in ADHD patients~\cite{robinson2016hypothesis}. Together, these approaches establish a new frontier for statistical analysis of shape variations in curved shape spaces, particularly when Euclidean projections distort distances.

\subsection{Learning Biological Patterns from Shape Space} 
Discovering and understanding how patterns, or regularities, of morphological variations encode evolutionary, functional, genetic, and ecological information is a central goal of biological phenotypic analysis. Shape space analysis provides an unprecedented suite of geometric and topological features that allows detection of both global and localized patterns of variations that may lead to the discovery of novel underlying biological factors.

The shape features of interest range from (a) global summaries of shapes such as compactness, sphericity, or aspect ratios to (b) geometric descriptors like curvatures, (c) spectral features derived from differential operators, as well as (d) statistical scores encoding how a shape relates to the dominant modes of variation within a population. Commonly, features are morphometric invariants that enable analysis independent of scaling, rotation, or translation.

By utilizing these features, one can quantify shape similarities and differences that are not immediately apparent in raw measurements, including the detection of subtle morphological trends, the grouping of related forms, and the identification of underlying constraints shaping the data. For instance, in evolutionary biology, pattern analysis of primate teeth can reveal the functional significance of cusps across different diets and track evolutionary changes in tooth morphology over time. Similarly, studying the patterns of cells can highlight developmental or functional processes, as well as abnormalities associated with disease or environmental influences. 

In order to learn the connection between shape features and biological factors, one can turn to multivariate regression of these features, which can reveal how shape encodes biological signals, for instance, by linking morphology to age, sex, ecological or phylogenetic traits, functional performance, or disease severity. Similar to the classification tasks, methods themselves are not necessarily new; rather, it is their application within shape space that is novel. Among the commonly employed approaches are multivariate regression models, Procrustes ANOVA accompanied by canonical correlation analysis, kernel regression, Gaussian process regression, multivariate analysis of variance (MANOVA), and measuring the Homogeneity of feature space~\cite{klingenberg_methods_2022, collyer_method_2015, zelditch_geometric_2012, mitteroecker2022thirty}. At times, it is better not to mine the entire shape space, rather than the solution lies in pairwise comparison of the data and deducing the information from it~\cite{faigenbaum2016algorithmic}. 

Across these examples, biological pattern analysis in shape space provides a rigorous framework for translating complex geometric data of phenotypic variation into interpretable insights, bridging the gap between raw measurements and scientific understanding. Below, we review major areas of application beyond classification. Though landmark-based analysis is still most popular in biological research, similar techniques can be adapted to shape space features, as reviewed in previous sections.

\subsubsection{Major Areas of Biological Applications}
A first major area of application is the study of size-shape relationships, such as how multivariate allometry drives shape changes in growth and evolution~\cite{klingenberg_methods_2022, mitteroecker2022thirty}. Multivariate regression and ANOVA with landmark data are most frequently used \cite{klingenberg_methods_2022, adams_ontogenetic_2010, mitteroecker2022thirty}. Notably, one can incorporate size as an additional dimension in the data to construct a form space to represent biological integration of size and shapes~\cite{klingenberg_methods_2022, mitteroecker2022thirty, richtsmeier_promise_2002}.

A second major area of application is the analysis of biological asymmetry in shapes. Shape asymmetry can be broadly classified into two categories: matching (between bilateral structures, e.g., left and right wings) and object (between connected halves of the same structure, e.g., within a wing). We can differentiate between directional asymmetry, which reflects consistent mechanisms at the population level, and fluctuating asymmetry, which captures individual deviations that may indicate developmental disturbances~\cite{klingenberg_analyzing_2015, rolfe2018associations}. Patterns in shape asymmetry are most commonly quantified using Procrustes ANOVA to parse and assess variance associated with symmetric and asymmetric components~\cite{klingenberg_analyzing_2015}. Other methods include using the Mantel test to assess matrix correlation between symmetric and asymmetric components~\cite{klingenberg_geometric_1998} and using Dense Correspondence Analysis (DeCA) to visualize and score the degree of facial asymmetry~\cite{rolfe_deca_2023, rolfe2025streamlining}.

Another major application is the study of integration and modularity~\cite{klingenberg_studying_2014, zelditch_what_2021, mitteroecker_conceptual_2007}. Integration refers to strong global covariation among traits driven by shared genetic, developmental, or functional mechanisms, whereas modularity describes stronger covariation within subsets of traits that reflect relatively independent processes. These concepts differ in degree rather than kind. For example, humans and great apes share similar global integration between the face and neurocranium, but humans show greater dissociation among facial submodules (e.g., the brow) during development and evolution~\cite{mitteroecker_ontogenetic_2009}.

Most empirical studies evaluate predefined modules based on developmental hypotheses~\cite{adams_comparing_2019, zelditch_what_2021}. Partial Least Squares (PLS) is the primary tool for quantifying integration between landmark blocks by identifying orthogonal axes of maximal covariation~\cite{zelditch_what_2021, zelditch_geometric_2012}. The strength of integration across modules can be summarized using RV coefficients or the covariance ratio (CR), and differences among groups are commonly assessed through CR-derived Z-scores~\cite{adams2016comparison, zelditch_what_2021}.

Different module partitions can be compared using the $Z_{12}$ statistic, which contrasts their CR scores~\cite{adams_comparing_2019}. Global integration can be quantified with relative eigenvalue variance ($V_{rel}$), which evaluates how the empirical distribution of trait variation deviates from expectations under perfect integration~\cite{conaway_effect_2022, watanabe_statistics_2022}. Beyond these summary measures, more complex frameworks have been developed. Hierarchical landmark partitioning examines how integration changes across nested levels of trait organization~\cite{claes_genome-wide_2018, claes_improved_2012, matthews_using_2023}. Network-based approaches such as MINT (Modularity and Integration by Node-based Testing) identify modules from covariation structure among landmarks~\cite{marquez_mint_2014, goswami_emmli_2016, adams_comparing_2019}, while EMMLi (Evaluating Modularity with Maximum Likelihood) statistically compares alternative modular hypotheses using a likelihood framework. However, both MINT and EMMLi require further validation before broad application in biological studies~\cite{zelditch_what_2021, adams_comparing_2019}.

Landmark-based shape data have also been applied in quantitative genetics to investigate the genetic basis of morphological variation. A common approach is to incorporate landmark configurations into the breeder’s equation, allowing estimation and visualization of predicted selection responses directly in shape space. For example, deviation between actual change in craniofacial outline from the expected selection response may indicate developmental constraint~\cite{klingenberg_quantitative_2001, klingenberg_evolution_2010}. Another major application is the use of landmark data in quantitative trait loci (QTL) mapping, where shape coordinates or derived axes are treated as phenotypic traits and tested for association with genomic regions, thereby linking specific loci to morphological modules~\cite{klingenberg_evolution_2010, navarro_does_2016, katz_facial_2020}. In genome-wide association studies (GWAS), either principal component (PC) scores or raw shape coordinates are used as traits to map against single nucleotide polymorphisms (SNPs). This enables the identification of genetic variants that contribute to morphological differences, such as facial variation in healthy populations or structural features associated with disease phenotypes~\cite{weinberg_3d_2016, white_insights_2021}.

\section{Practical Guidance to Shape Space Analysis, through the Lens of Two Examples} 
\label{sec:practicalGuide}
In this chapter, we provide practical guidance for conducting shape-space analysis, drawing on two concrete examples that span different scales and contexts. These examples illustrate key principles and methodological choices that can inform future studies. The first example is of microscale data, focusing on subcellular morphology, while the other one is a macroscale analysis of primate teeth. We decided to include these examples as a guide for future projects, clarifying the sequence of steps and the core analytical tools involved. This chapter brings together all the ideas presented so far into a coherent workflow that can be readily applied to other cases. We also include a more philosophical discussion of practical questions that arise when working within shape space.

\subsection{Micro-Scale -- Subcellular Morphology and Signal Dynamics across Space and Time}

At the microscale, shape analysis is often used to study biological cells, where shape is typically described as cell morphology. The motivation traces back to the notion ``form is function'': whereby cells with a characteristic function also appear to adopt a characteristic shape. Classic examples include neurons, hematopoietic cells, muscle cells, adipocytes, sperm, and egg cells. As a result, characterizing cell shape has become essential in areas such as drug screening and disease prognosis, where the goal is often to identify or predict specific biological phenotypes.

\textbf{An unexpected turn.} Traditionally, studies of cell morphology follow a standard pipeline: extract handcrafted or learned shape descriptors, optionally combine them with molecular signatures, and then use clustering or predictive modeling to identify phenotypes or infer biological conditions~\cite{zhou2025development,zhou2025identifying,de2025geometric}. This framework rests on the assumption that shape is largely the outcome of cellular signaling pathways. However, advances in fluorescent biosensors and high-resolution microscopy have begun to challenge this assumption. These technologies make it possible to measure molecular activity across the cell surface with unprecedented spatial detail. Strikingly, this has led to findings-such as those by Weems et al.~\cite{weems2023blebs} that reshape our understanding. Their study showed that metastatic melanoma cells generate dynamic hemispherical protrusions, or blebs, not merely to aid in migration, as conventionally thought, but also to activate pro-survival signaling pathways to evade cell death checkpoints and promote metastasis. This discovery raises deeper questions: ``What are the causal mechanisms by which shape drives signaling either in conjunction with or independent of transcription and translation; what context are these mechanisms activated; and how significant are these mechanisms in affecting a functional cell response?'' Answering these questions requires quantitative methods that can track a cell’s highly irregular 3D surface, decorated with transient, high-curvature protrusions, while simultaneously capturing the molecular activity associated with these regions. This, in turn, enables correlation and causal-inference analyses that do not rely on experimental perturbations that may not reflect behavior in the unperturbed system~\cite{vilela2011s, welf2014using, noh2022granger}.\\

\begin{figure}[t!]
\centering
\includegraphics[width=\linewidth]{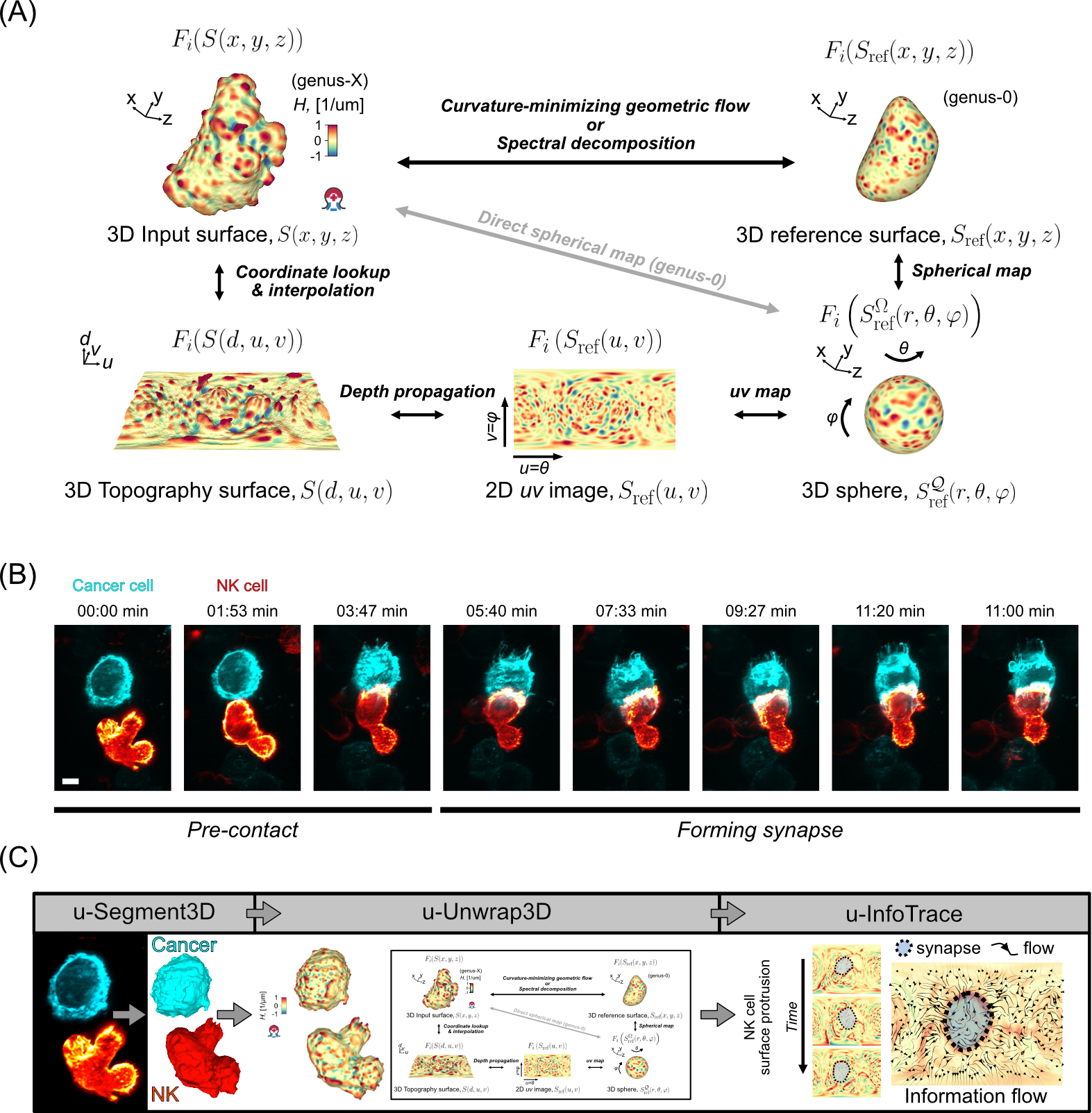} 
\caption{\textbf{Transformation between 3D surface representations to investigate the causal mechanisms underlying cell surface protrusion dynamics.} (A) Overview of the u-Unwrap3D~\cite{zhou2023surface} workflow, which enables bidirectional transformation between a Cartesian 3D surface and multiple surface representations. (B) Snapshots illustrating morphological changes as a natural killer cell recognizes and engages a cancer cell, forming a physical interface (the immunological synapse). (C) Analysis of protrusion dynamics during synapse formation. Complete 3D cell surfaces are obtained using u-Segment3D~\cite{zhou2025universal}; surfaces were then registered spatiotemporally with u-Unwrap3D~\cite{zhou2023surface}, and protrusions, quantified by mean curvature, are mapped into 2D unwrapped images. Consensus flow in transient protrusion dynamics is identified using multiscale information flow computed by u-InfoTrace~\cite{zhou2022multiscale}, equipped with a dynamic differential covariance causal measure~\cite{chen2022ddc}.
}
\label{fig:u_Unwrap3D_causality_micro_flow}
\end{figure}

\textbf{A General Framework for 3D Cell Surface Analysis.} Towards enabling spatiotemporal causal inference of morphology and signaling in dynamic 3D cells, Zhou et al. developed u-Unwrap3D~\cite{zhou2023surface} - a first comprehensive framework that addresses both correspondence and sampling problems introduced by segmentation errors and topological changes for 3D cell surfaces imaged by high-resolution confocal or light-sheet microscopy. u-Unwrap3D is an open-source Python package that maps a segmented 3D cell surface mesh into several complementary representations (Fig. \ref{fig:u_Unwrap3D_causality_micro_flow}). A central idea in u-Unwrap3D is the construction of a genus-0 reference surface, $S_{\text{ref}}$-a smoothed, protrusion-free version of the cell. This reference surface, generated using curvature-minimizing geometric flows such as conformalized mean curvature flow (cMCF)~\cite{kazhdan2012can}, active contours~\cite{kass1988snakes}, or by spectral decomposition~\cite{liu2017dirac}, enables definition and isolation of protrusions on the original surface. The reference is made genus-0 by shrinkwrapping and subsequently mapped to canonical spherical and rectangular $(u,v)$ coordinates. The spherical map is initially quasi-conformal, then the vertices are advected to linearly diffuse the area distortion to obtain a quasi-equiareal parameterization of $S_{\text{ref}}$. The $(u,v)$ parameterization is constructed to (i) maximize the mapping of surface features-of-interest to the equator, which has the lowest metric distortion, through specification of a vertex-based weight function, and solving a weighted eigendecomposition of the sphere geometry; and (ii) optimize aspect ratio by minimizing the Beltrami coefficient to obtain the lowest global distortion. Lastly, the $(u,v)$ parameterized $S_{\text{ref}}$ is propagated along the gradient of the distance transform of $S_{\text{ref}}$, to map the proximal surface volume (or even to map the entire inner volume). This generates a topographic surface representation that selectively reconstructs the 3D geometry of protrusions. Each representation offers unique advantages: the original input surface for geometrical measurements; $S_{\text{ref}}$ for defining protrusion-of-interest and simplifying the dense registration problem; canonical sphere and $(u,v)$ parameterizations for performing tracking and registration, and topography for segmenting individual protrusions and performing protrusion-specific operations, e.g., measuring protrusion-specific geometrical parameters like height. By taking advantage of each representation and the bidirectional mapping between them, Zhou et al. demonstrated how u-Unwrap3D facilitates the analysis of unconstrained cell geometries and enables quantitative analysis of previously challenging spatiotemporal data. Notably, they establish a retrograde flow of protrusions on natural killer cells in immunological synapse formation with cancer; measure information flows~\cite{zhou2022multiscale} to causally link K14+ cells with invasive branching in breast organoids; perform unsupervised instance segmentation of protrusions regardless of morphological complexity; quantify bleb-mediated recruitment event of the Septin molecule involved in melanoma cancer cell survival~\cite{weems2023blebs}; and measure the speed of membrane ruffling.  

\subsection{Macro-Scale -- Morphological Evolution across Primate Species through Tooth Shape Analysis}

\begin{figure}[t]
\centering
\includegraphics[width=\linewidth]{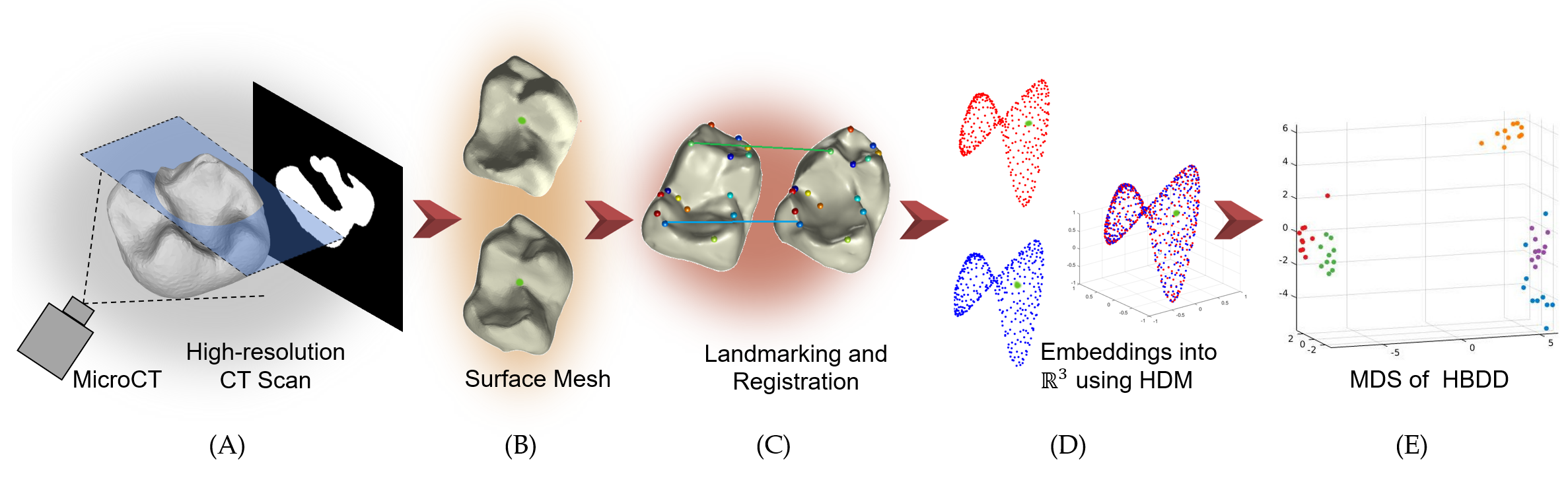} 
\caption{\textbf{Illustration of the main steps in the workflow for analyzing molar tooth morphology in an evolutionary context: from physical specimens to quantitative shape comparison.} (A) CT scanning of molar teeth~\cite{yapuncich2019digital}; (B)~Surface reconstruction from scan data to produce 3D meshes; (C) Alignment, registration, and automatic detection of biologically meaningful landmarks, with correspondence identification using Gaussian process landmarking~\cite{gao2019gaussian}; (D) Definition of geometric distances and embedding of surfaces into a common template space; (E) Quantification of anatomical surface similarity using the HBDD~\cite{gao2021diffusion}, followed by embedding into $\mathbb{R}^3$ via MDS.
}
\label{fig:shape_space_macro_Flow}
\end{figure}

At the macro scale, shape space analysis is frequently used to investigate morphological evolution across species~\cite{choi2018planar, porto2009evolution, faigenbaum2026studying}. A particularly compelling example investigates the morphology of primate teeth to understand how form diversifies under different dietary pressures. This end-to-end shape analysis project was carried out at Duke University by Ingrid Daubechies and her research group, in collaboration with evolutionary biologist Doug M. Boyer, and is summarized in~\cite{faigenbaum2026studying}. The process begins with the acquisition of high-resolution 3D models, derived from surface scanning or computed tomography (CT) (see Fig.~\ref{fig:shape_space_macro_Flow}A). These scans are processed to generate point clouds and surface meshes, which are then cleaned and completed through several post-processing steps (Fig.~\ref{fig:shape_space_macro_Flow}B).

Principal Component Analysis (PCA) is first applied to the aligned shapes to extract major axes of variation. This is followed by optimization of the Continuous Procrustes Distance, which identifies the best point correspondences and geometric distance between each pair of surfaces. Gaussian process landmarking~\cite{gao2019gaussian} further refines the analysis (Fig.~\ref{fig:shape_space_macro_Flow}C). With accurate shape comparisons established, the research turns to statistical evaluation of morphological variation and evolutionary hypotheses, using the Continuous Procrustes framework~\cite{al2013continuous}. Tools from differential geometry are central to this phase, including the Horizontal Diffusion Map (HDM)~\cite{gao2021diffusion}, which embeds the surface coordinates into a low-dimensional space based on eigenvectors of the Horizontal Laplace-Beltrami operator (Fig.~\ref{fig:shape_space_macro_Flow}D). Here, colors (e.g., red and blue) represent points from different teeth projected into diffusion space. Finally, Multidimensional Scaling (MDS) is applied to the horizontal-based diffusion distance (HBDD) matrix, producing a shape space that reflects known dietary groupings: the folivores \textit{Alouatta} (red) and \textit{Brachyteles} (green) cluster together, as do the frugivores \textit{Ateles} (blue) and \textit{Callicebus} (purple); the insectivore \textit{Saimiri} (orange), by contrast, lies far from these herbivorous groups. Beyond clustering, the pipeline enables extraction of additional geometric and biological features, including curvature estimation~\cite{shan2019ariadne} and consistent segmentation of shapes across the collection.

\section{Shape Space Datasets}

\label{sec:datasets}

There are numerous challenges that naturally fall under the umbrella of shape space analysis. In fact, many of the core problems addressed in data science, ranging from clustering and classification to representation learning, can be reinterpreted through the lens of shapes and the spaces they inhabit. This perspective highlights the geometric and structural aspects of data that are often critical for meaningful analysis. In this section, we aim to present a diverse set of examples drawn from various research domains. We begin with applications in biological and medical imaging, where the geometry of anatomical structures plays a key role, continue with questions in evolutionary biology that involve comparing morphologies across species, and conclude with datasets arising in computer graphics, where shape representation and manipulation are central. These examples illustrate the wide applicability of shape space thinking and its potential to unify seemingly distinct problems through a shared mathematical framework.

Collecting a dataset of sufficient size for a machine learning task is a scrupulous process, requiring significant time, effort, and financial support. We would therefore like to sincerely thank all those who have chosen to make their data publicly available. Sharing these datasets not only helps others develop better algorithms but also advances the field by enabling reproducibility and new insights.

\subsection{Established Benchmark Datasets} In the realm of shape analysis and recognition, a collection of standard benchmark datasets has emerged as resources for evaluating algorithms and methodologies. For 2D image data, the well-known MNIST~\cite{lecun1998mnist} and its extension EMNIST~\cite{cohen2017emnist} are classic benchmarks featuring handwritten digits, commonly used for image recognition models, as well as ancient handwriting~\cite{faigenbaum2016algorithmic}. Similarly, benchmarks are available for face/pose recognition (e.g., CMU Pose, Illumination, and Expression (PIE) database~\cite{sim2002cmu}).
In the 3D domain, 3D polygonal models are available for segmentation, classification, and deep learning, such as the Stanford 3D scanning repository~\cite{turk1994zippered}, PrincetonSB~\cite{shilane2004princeton}, ShapeNet~\cite{chang2015shapenet}, and ABC~\cite{koch2019abc}. Collectively, these datasets form a robust foundation for 2D and 3D shape analysis, enabling standardized benchmarking across various applications and methodologies.

\subsection{Data from the Real World} In addition to benchmarking datasets, real-world datasets play a crucial role in advancing the community of shape analysis, as they provide the complexity and variability found in actual environments. It is not only that these datasets originate from diverse fields, but each also presents unique computational challenges and raises distinct research questions shaped by the specific context of its field. The differences in data modality and structure, scale, and measurement modalities often require tailored analytical approaches, highlighting the need for versatile methods capable of addressing both technical and domain-specific complexities. Data tell a story, and exploring its geometric features can uncover narratives about the underlying structure, patterns, and relationships within the system, revealing insights that might be obscured by traditional manual analyses. A non-exhaustive list of real-world datasets from various fields is provided in Table~\ref{tab:bio_med_data} and Table~\ref{tab:other_data}, along with references and links to data. 

We begin by surveying biological and medical imaging datasets. Biological applications constitute only a small portion of shape space studies, spanning scales from the micro to the macro (see discussion in Section~\ref{sec:practicalGuide}).

\begin{table}[t!]
\centering
\begin{tabular}{|p{3.5cm}|p{7.8cm}|p{1.7cm}|p{2.8cm}|}
\hline
\rowcolor{gray!30}
Dataset & Short description & Ref. & Link \\\hline
\hline
Cell contours via fluorescence microscopy & 2D coordinates of cancer cell contour shapes under different treatments & \cite{Alizadeh2019-ua,cancercell2} & \href{https://github.com/geomstats/geomstats/tree/main/geomstats/datasets/data/cells}{Bone cancer}, \href{https://github.com/wxli0/dyn/tree/main%4092c7a58/dyn/datasets/breast_cancer}{Breast cancer}\\\hline
2D cell images & Large-scale benchmark data for live cell images from microscopy with manual annotation of cell types & \cite{Edlund2021} & \href{https://paperswithcode.com/dataset/livecell}{LIVECell}\\\hline
Atomic structure of ribosome & Visualization and analysis for 3D ribosome structures & \cite{kushner2023riboxyz} & \href{https://pdb101.rcsb.org}{Protein Data Bank}, \href{https://www.emdataresource.org}{EMDataResource}, \href{https://ribosome.xyz}{RiboXYZ}\\\hline
Cryo-EM density maps for a community challenge & Synthetic density maps deforming along a global motion from molecule dynamics & \cite{cryoemchallenge,Astore2025} & \href{https://osf.io/8h6fz/}{Dataset}\\\hline
3D cell images & Cell images in 3D with more than 200,000 live cells & \cite{Viana2023, zhou2023surface} & \href{https://open.quiltdata.com/b/allencell/packages/aics/hipsc_single_cell_image_dataset}{Dataset 1},  \href{https://zenodo.org/records/8166238}{Dataset~2}\\\hline
3D neuron cell images & Neurons reconstructed from Patch-seq technology & \cite{govek2023cajal} & Data availability section of \cite{govek2023cajal} \\\hline
Brain imaging data & Brain images via MRI, PET, MEG, EEG, and iEEG & \cite{OpenNeuro, AllenBrain} & \href{https://openneuro.org}{OpenNeuro}, \href{https://portal.brain-map.org}{Allen Brain Map}\\\hline
Several Human Neuroimaging data & Large and comprehensive set of correlated clinical data, Manually Labeled MRI Brain Scan Database, etc. & \cite{Jack2008, kennedy2016nitrc, luo2009neuroimaging} & \href{https://www.nitrc.org/}{NITRC}, \href{https://adni.loni.usc.edu/data-samples/adni-data/neuroimaging/mri/mri-image-data-sets/}{ADNI}\\\hline
BIOCARD & Predictors of Cognitive Decline Among Normal Individuals; MRI, PET-PiB, Neuropathology data, Diagnoses & \cite{miller2014amygdalar} & \href{https://www.gaaindata.org/partner/BIOCARD}{BIOCARD}\\\hline
Colon10K & Colonoscopic images sampled from a video & \cite{ma2021colon10k} & \href{https://endoscopography.web.unc.edu/place-recognition-in-colonoscopy/}{Colon10K}\\\hline
MedShapeNet & Large-scale dataset of 3D medical shapes & \cite{li2025medshapenet} & \href{https://medshapenet.ikim.nrw/}{MedShapeNet}\\\hline

\end{tabular}
\caption{Biological and medical imaging datasets.}
\label{tab:bio_med_data}
\end{table}

\begin{table}[t!]
\centering
\begin{tabular}{|p{3.6cm}|p{7.5cm}|p{1.8cm}|p{2.8cm}|}
\hline
\rowcolor{gray!30}
Dataset & Short description & Ref. & Link \\\hline
MorphoSource & 3D meshes and CT-images of biological skeletal material & \cite{boyer2016mophosource} & \href{https://www.morphosource.org/users/sign_in?locale=en}{MorphoSource}\\\hline
Second mandibular surfaces & 116 Second mandibular molars of prosimian primates and non-primate, as well as observer-determined landmark & \cite{boyer2011algorithms} & \href{https://www.wisdom.weizmann.ac.il/~ylipman/CPsurfcomp/}{Dataset}\\\hline
Fruit fly wings & Wing images taken on microscopes with magnification & \cite{FruitFlyWing, edwards2007database} & \href{http://gigadb.org/dataset/100141}{Dataset}\\\hline
Leaf shapes & 2D coordinates of leaf shape across plant groups & \cite{LeafData} & \href{https://figshare.com/articles/dataset/Leaf_coordinates_zip/5056441/1?file=8559085}{Dataset}\\\hline
Cloud images & 2D cloud image data with annotation of cloud types & \cite{CloudShape, Rasp2020} & \href{https://github.com/raspstephan/sugar-flower-fish-or-gravel}{Dataset}\\\hline
Human Body Shapes & Time Series of 3D Meshes Representing Human Bodies/Faces & \cite{dfaust:CVPR:2017,COMA:ECCV18} & \href{https://dfaust.is.tue.mpg.de}{DFAUST}, \href{https://coma.is.tue.mpg.de}{COMA}\\\hline
3D morphable models & 3D data for faces, ears, bodies, and others & \cite{egger20203d,CAPOD,TOSCA} & \href{https://github.com/3d-morphable-models/curated-list-of-awesome-3D-Morphable-Model-software-and-data}{Curated List}, \href{https://sites.google.com/site/pgpapadakis/CAPOD}{CAPOD}, \href{https://cvg.cit.tum.de/data/datasets/partial}{Partial}\\\hline
Shark feeding apparatus & 3D landmark configurations of basking shark gill arches & \cite{Paskin2022data} & \href{https://github.com/morphomatics/3Dfrom2DLandmarks/tree/main}{Data}\\\hline
Blender database & databases provided by Blender, including human based meshes, poses etc. & NA & \href{https://www.blender.org/download/demo-files/}{Blender}\\\hline
Stone tools & Use-related shape trajectories of Neolithic stone axes from Central Europe & \cite{mayer2025} & \href{https://zenodo.org/communities/wear/}{WEAR}\\\hline
Biomedical collection & Collection of small real data ranging from ribosome tunnel data to cell imaging time series  & NA & \href{https://osf.io/56v42/?view_only=dfd518264b20403bbf72fc9b9d9a6df9}{Khanh Dao Duc Lab} \\\hline
Imaging Data Commons & National Cancer Institute open imaging repository with web-based 3D viewer powered by deep learning segmentation and analytical tools. & \cite{fedorov2023national} & \href{https://datacommons.cancer.gov/repository/imaging-data-commons}{NCI Imaging Data Commons}\\\hline
\end{tabular}
\caption{Other datasets.}
\label{tab:other_data}
\end{table}

\section{Software Tools}
\label{sec:tools}

A wide range of problems in shape space have been studied - each motivated by the specific challenges posed by particular datasets (see Section~\ref{sec:datasets}). We aimed to provide an extensive list of the existing tools for others to apply them easily. However, the breadth of available methodological tools is far too extensive to be exhaustively covered in this paper. In particular, the rapidly growing body of machine learning–based approaches has produced a long and continually expanding list of relevant contributions. Rather than attempting an incomplete survey, we decided to maintain a more comprehensive and regularly updated list of related methods and references in a dedicated GitHub repository: \url{https://github.com/shirafaigen/ShapeSpaceSurvey.git}.

\section{Discussion and Future Research Direction}
\label{sec:Discussion}

In this paper, we discussed a wide range of methods pertaining to the shape space. We reviewed a representative - though necessarily non-exhaustive-set of methods for analyzing shape space. Since both the scientific question and the nature of the data strongly influence methodological choices, our goal was not only to catalog existing approaches but also to provide guidance on how specific data characteristics motivate the use of particular methods. We began by discussing different types of shape representations, including implicit and explicit formulations, and outlined the types of problems for which each is most appropriate. We then reviewed several major representational paradigms, including function spaces of parameterizations, measure-based representations, and feature-based representations. A substantial portion of the paper is devoted to methodological foundations and preprocessing steps, such as shape correspondence, landmarking, and registration. In particular, we emphasize that landmarking remains a critical component of comparative shape analysis, as the identification of anatomically or structurally meaningful points often plays a decisive role throughout the analysis pipeline. 

We further discussed the presentation and organization of shape feature descriptors, distinguishing between global and local statistics, distributional and topological characteristics, decomposition-based methods, data-driven descriptors, and meta-feature descriptors constructed from existing shape features. We also highlighted the central importance of metric choice for capturing meaningful discrepancies between shapes, including concepts such as shape trajectories and mean shapes. Then, once shapes and metrics were defined, we addressed the range of questions that can be tackled using machine learning to uncover morphological patterns in shape space, with an emphasis on methods that are specifically designed for shape-structured data rather than generic machine learning techniques. To illustrate how the full pipeline operates in practice, we provided practical guidance for conducting shape-space analysis, drawing on two concrete examples that span different scales and contexts, from micro- to macroscopic settings. Finally, we presented several representative shape-space datasets and discussed existing tools, accompanied by a GitHub repository that is continuously updated.

Before we close this chapter and conclude that all the work has been done, this paper is by no means a concluding chapter of shape analysis. Rather than closing the book, each study we reviewed opens new avenues for research - some inspired by emerging datasets, others by unresolved theoretical questions. Below, we outline several key challenges that continue to shape the field of shape space analysis. 
\begingroup
\renewcommand\labelenumi{(\theenumi)}
\begin{enumerate}

\item \textbf{AI revolution.} With the rapid dissemination of what is often referred to as the ``AI revolution,'' how does it affect the shape space analysis? Advances in machine learning-particularly in representation learning, geometric deep learning, and data-driven modeling-have expanded the range of questions that can be addressed. So is the shape space part of this revolution, and how will the shape space research be affected by it?

\item \textbf{Out-of-the-box solution?} Existing algorithms can rarely be applied as is, and typically must be carefully adapted to the specific data and scientific question at hand. Differences in data origin, data topology, scale, noise characteristics, and sampling often require tailoring each component of the analysis. Consequently, a central question that arises is not merely which individual methods to use, but how to coherently design and apply the entire shape-space analysis pipeline, in a way that is appropriate for the data and the phenomena being studied.

\item \textbf{Are landmarks essential?} In evolutionary and developmental biology, shape analysis has been dominated by landmark-based geometric morphometrics, partly because landmark data integrate well with established multivariate statistical frameworks. As a result, landmark-based approaches are not well suited to questions that require identifying and comparing individual morphological characters across samples. This misalignment helps explain why shape data are rarely used directly in phylogenetic reconstruction and character-based systematics, where the explicit representation of localized traits remains central to the discovery of novelties and the reconstruction of phylogenetic relationships~\cite{adams_morphometrics_2011, schwartz2004getting, mongiardino2017discrete, peterson2016phenotypic, palci2019geometric}. Therefore, in addition to integrating shape features and shape space attributes with statistical analysis, complementary shape representations are needed to study localized morphological features that may carry distinct developmental, functional, or evolutionary significance.

\item \textbf{Shape registration: close, but not quite there.} Despite decades of research, shape registration remains a challenging problem, with existing methods often struggling in the presence of large deformations, topological variability, noise, or weakly defined correspondences. Moreover, how can efficient shape registration be achieved without relying on initial correspondences? 
What is the solution for incomplete or partially corrupted shapes?

\item \textbf{Evaluation of the developed methods.} A systematic study based on reliable ground truth (GT) data is still largely missing, along with standardized quantitative measures for evaluating how well different methods perform in regular scenarios, as well as addressing how space analysis can be made robust against noise, missing data, and artifacts introduced during acquisition. Are there mathematical guarantees for the robustness and stability of the algorithms in noisy datasets? In particular, there is a lack of agreed-upon benchmarks and evaluation protocols for assessing the quality of key components of the pipeline, such as registration accuracy and shape deformation fidelity. This gap makes it difficult to compare methods objectively, understand their failure modes, and determine which approaches are most appropriate for specific data regimes or scientific questions.

    \item \textbf{Multi-modal and hybrid representations.} How can heterogeneous data modalities-such as 2D images, point clouds, surface meshes, and volumetric representations-be integrated into a unified and coherent shape-space framework? Addressing this question requires representations and metrics that can bridge differing dimensionalities, resolutions, and noise characteristics, while preserving meaningful geometric and semantic information across modalities.

\item \textbf{Multiscale and hierarchical shape representations.} While today each data is addressed separately, we anticipate creating an integrative representation across scales, e.g., describing both the protein structure, protein distribution in the cell, cell morphology, and cell environment, will be beneficial. 

\item \textbf{Model-based learning.} Integrating shape-space knowledge with data-driven systems-particularly deep learning-raises the question of how geometric constraints, invariances, and physically or biologically motivated priors can be incorporated into learning frameworks in a principled way. Combining explicit shape models with learned representations offers the potential to improve interpretability, data efficiency, and generalization, while ensuring that learned solutions remain consistent with the underlying geometric structure of shapes.

\end{enumerate}
\endgroup

\section*{Acknowledgments}
Gary P. T. Choi thanks the Croucher Tak Wah Mak Innovation Award for supporting his research. Shira Faigenbaum-Golovin is grateful to the Eric and Wendy Schmidt Fund for Strategic Innovation, the Zuckerman-CHE STEM Program, Duke University, and in part the Simons Foundation under Grant Math+X 400837 for supporting her research. Karen Habermann gratefully acknowledges support from the EPSRC-funded Prob\textunderscore AI Hub, with grant reference EP/Y028783/1. Wenjun Zhao gratefully acknowledges that this research was supported in part by the Pacific Institute for the Mathematical Sciences and the Simons Foundation.

\bibliographystyle{abbrv}
\bibliography{references.bib}

\end{document}